%% file: all.tex
\documentclass[12pt,reqno]{amsart}

\usepackage{a4wide}
\usepackage{xcolor}
\definecolor{darkgreen}{rgb}{0,0.45,0}
\definecolor{darkred}{rgb}{0.75,0,0}
\definecolor{darkblue}{rgb}{0,0,0.6}
\usepackage[pdfborder=0,colorlinks, anchorcolor=darkred,citecolor=darkgreen,linkcolor=darkgreen,urlcolor=darkblue, unicode, final]{hyperref}
\usepackage[T1]{fontenc}    
\usepackage{enumitem}
\usepackage{graphicx}

\usepackage{amsmath}
\usepackage{amssymb}
\usepackage{amsxtra}
\usepackage{amsthm}

\newif\ifnewadd
\newif\ifnewcut
\input{imports}

\setlist{}
\setenumerate{leftmargin=*,labelindent=\parindent}
\setitemize{leftmargin=*,labelindent=0.5\parindent}
\setdescription{leftmargin=1em}

\swapnumbers

\theoremstyle{plain}

\newtheorem{thm}{Theorem}[subsection]
\newtheorem{lem}[thm]{Lemma}
\newtheorem{cor}[thm]{Corollary}

\newtheorem{prop}[thm]{Proposition}

\theoremstyle{definition}
\newtheorem{defn}[thm]{Definition}
\newtheorem{ex}[thm]{Example}
\newtheorem{ntn}[thm]{Notation}

\theoremstyle{remark}
\newtheorem{obs}[thm]{Observation}
\newtheorem{rec}[thm]{Recall}
\newtheorem{rmk}[thm]{Remark}

\makeatletter
\let\c@equation\c@thm
\makeatother
\numberwithin{equation}{subsection}

\title{Cartesian exponentiation and monadicity}

\author[Riehl]{Emily Riehl}
\address{
  Department of Mathematics \\
Johns Hopkins University \\
Baltimore, MD 21218\\
  USA
}
\email{eriehl@jhu.edu}

\author[Verity]{Dominic Verity}
\address{
  Mathematical Sciences Institute\\
  The Australian National University\\
  ACT 2601\\
  Australia
}
\email{dominic.verity@anu.edu.au}

\date{\today}


\begin{document}

  \ifpdf
  \DeclareGraphicsExtensions{.pdf, .jpg, .tif}
  \else
  \DeclareGraphicsExtensions{.eps, .jpg}
  \fi

  \begin{abstract}
An important result in quasi-category theory due to Lurie is that the cocartesian fibrations are \emph{exponentiable}, in the sense that pullback along a cocartesian fibration admits a right Quillen right adjoint that moreover preserves cartesian fibrations; the same is true with the cartesian and cocartesian fibrations interchanged. To explicate this classical result, we prove that the pullback along a cocartesian fibration between quasi-categories forms the oplax colimit of its ``straightening,'' a homotopy coherent diagram valued in quasi-categories, recovering a result first observed by Gepner, Haugseng, and Nikolaus. As an application of the exponentiation operation of a cartesian fibration by a cocartesian one, we use the Yoneda lemma to construct left and right adjoints to the forgetful functor that carries a cartesian fibration over $\qB$ to its $\ob\qB$-indexed family of fibers, and prove that this forgetful functor is monadic and comonadic. This monadicity is then applied to construct the reflection of a cartesian fibration into a \emph{groupoidal} cartesian fibration, whose fibers are Kan complexes rather than quasi-categories.
  \end{abstract}

  \maketitle
  \tableofcontents

  \input{intro}
  \ifnewadd
  \input{background}
  \fi
  \ifnewcut\else
  \input{weighted}
  \input{cartesian}
 \input{comprehension-review}
\fi
 \input{pullback}
  \input{exponentiation}
  \input{comonadicity}

\input{reflection}

  \bibliographystyle{abbrv}
  \bibliography{../../common/index}


\end{document}


%% file: imports.tex
\usepackage[silent]{fixme}
\fxsetup{
  status=draft,
  layout={marginclue,noinline}
}
\fxsetface{inline}{\itshape}

\input{macros}

\usepackage{mathtools}

\makeatletter
\newcommand{\qc}[1]{{\mathord{\text{\normalfont{\textsf{#1}}}}}}
\newcommand{\qop}[1]{{\mathord{\qc{#1}}}}
\def\ec@#1#2<.>{\mathcal{#1}\mkern-2mu\text{\normalfont{\textsf{\slshape #2}}}}
\newcommand{\ec}[1]{{\mathord{\ec@#1<.>}}}
\newcommand{\eop}[1]{{\mathord{\ec{#1}}}}
\makeatother

\newcommand{\qA}{\qc{A}}
\newcommand{\qB}{\qc{B}}
\newcommand{\qC}{\qc{C}}
\newcommand{\qD}{\qc{D}}
\newcommand{\qE}{\qc{E}}
\newcommand{\qF}{\qc{F}}
\newcommand{\qG}{\qc{G}}
\newcommand{\qK}{\qc{K}}
\newcommand{\qL}{\qc{L}}

\newcommand{\eK}{\ec{K}}
\newcommand{\eL}{\ec{L}}

\newcommand{\eA}{\ec{A}}

\newcommand{\eC}{\ec{C}}
\newcommand{\eD}{\ec{D}}
\newcommand{\eE}{\ec{E}}

\newcommand{\Hom}{\qop{Hom}}
\newcommand{\Fun}{\qop{Fun}}

\renewcommand{\coCart}{\qop{coCart}}
\newcommand{\Cart}{\qop{Cart}}

\newcommand{\sCptd}{\sSet\text{-}\category{Cptd}}

\newcommand{\SSet}{\eop{SSet}}
\renewcommand{\qCat}{\eop{QCat}}
\newcommand{\qqCat}{\qop{qCat}}
\newcommand{\qKan}{\qop{Kan}}
\renewcommand{\Kan}{\eop{Kan}}

\renewcommand{\msSet}{\sSet_\bullet}

\newcommand{\gr}{^{\mathrm{gr}}}

\newcommand{\Cube}{\mathord{\sqcap\llapm\sqcup}}
\newcommand{\CHorn}{\mathord{\sqcap\llapm\sqcap}}

\newmop{core}
\newmop{Fam}
\newmop{last}
\newmop{rv}
\newmop{ob}

\newcommand{\sslice}{{\mathord{\mathord{/}\!\!\mathord{/}}}}

\usepackage{xr}

\newcommand{\extRef}[3]{%
  {\protect\IfBeginWith{#3}{itm:}{}{#2.}}\ref*{#1:#3}}

\newcommand{\refI}{\extRef{found}{I}}
\newcommand{\refII}{\extRef{cohadj}{II}}
\newcommand{\refIII}{\extRef{complete}{III}}
\newcommand{\refIV}{\extRef{yoneda}{IV}}

\newcommand{\refVI}{\extRef{comprehend}{VI}}
\newcommand{\refVII}{\extRef{holimits}{VII}}
\newcommand{\refVIII}{\extRef{construction}{VIII}}

%% file: macros.tex
\usepackage[mathscr]{eucal} 
\usepackage{mathrsfs}       
\usepackage{slashed}        

\usepackage{mathtools}
\usepackage{extpfeil}
\usepackage{proof}
\usepackage{ifnextok}
\usepackage{xstring}

\usepackage[all]{xy}
\SelectTips{cm}{}
\SilentMatrices
\ifx\pdfoutput\undefined
  \xyoption{dvips}
\fi
\usepackage{tikz-cd}

\newcommand{\ifnotdef}[1]{\expandafter\ifx\csname#1\endcsname\relax}%

\DeclareMathAlphabet{\mathbbe}{U}{bbold}{m}{n}

\makeatletter

\def\re@DeclareMathSymbol#1#2#3#4{%
    \let#1=\undefined
    \DeclareMathSymbol{#1}{#2}{#3}{#4}}

\ifnotdef{Top}
  \DeclareSymbolFont{tcSyC}{U}{txsyc}{m}{n}
  \SetSymbolFont{tcSyC}{bold}{U}{txsyc}{bx}{n}
  \DeclareFontSubstitution{U}{txsyc}{m}{n}

  \re@DeclareMathSymbol{\Top}{\mathord}{tcSyC}{120}
  \re@DeclareMathSymbol{\Bot}{\mathord}{tcSyC}{121}
\fi

\ifnotdef{righthalfcup}
  \DeclareFontFamily{U}{MnSymbolC}{}
  \DeclareSymbolFont{mnSyC}{U}{MnSymbolC}{m}{n}
  \SetSymbolFont{mnSyC}{bold}{U}{MnSymbolC}{b}{n}
  \DeclareFontShape{U}{MnSymbolC}{m}{n}{
      <-6>  MnSymbolC5
     <6-7>  MnSymbolC6
     <7-8>  MnSymbolC7
     <8-9>  MnSymbolC8
     <9-10> MnSymbolC9
    <10-12> MnSymbolC10
    <12->   MnSymbolC12}{}
  \DeclareFontShape{U}{MnSymbolC}{b}{n}{
      <-6>  MnSymbolC-Bold5
     <6-7>  MnSymbolC-Bold6
     <7-8>  MnSymbolC-Bold7
     <8-9>  MnSymbolC-Bold8
     <9-10> MnSymbolC-Bold9
    <10-12> MnSymbolC-Bold10
    <12->   MnSymbolC-Bold12}{}

  \re@DeclareMathSymbol{\righthalfcup}{\mathord}{mnSyC}{184}
  \re@DeclareMathSymbol{\lefthalfcap}{\mathord}{mnSyC}{185}
\fi

\DeclareFontFamily{U}{MnSymbolA}{}
\DeclareSymbolFont{mnSyA}{U}{MnSymbolA}{m}{n}
\SetSymbolFont{mnSyA}{bold}{U}{MnSymbolA}{b}{n}
\DeclareFontShape{U}{MnSymbolA}{m}{n}{
    <-6>  MnSymbolA5
   <6-7>  MnSymbolA6
   <7-8>  MnSymbolA7
   <8-9>  MnSymbolA8
   <9-10> MnSymbolA9
  <10-12> MnSymbolA10
  <12->   MnSymbolA12}{}
\DeclareFontShape{U}{MnSymbolA}{b}{n}{
    <-6>  MnSymbolA-Bold5
   <6-7>  MnSymbolA-Bold6
   <7-8>  MnSymbolA-Bold7
   <8-9>  MnSymbolA-Bold8
   <9-10> MnSymbolA-Bold9
  <10-12> MnSymbolA-Bold10
  <12->   MnSymbolA-Bold12}{}

\re@DeclareMathSymbol{\twoheadedswarrow}{\mathord}{mnSyA}{30}

\makeatother

\newcommand{\mlaux}[3]{\setbox0=\hbox{$\mathsurround=0pt #2{#3}$}%
  \dimen0=\dp0\advance\dimen0 by \ht0\lower#1\dimen0\box0}

\newcommand{\makellapm}[2]{\hbox to 0pt{\hss$\mathsurround=0pt #1{#2}$}}
\newcommand{\llapm}{\relax\mathpalette\makellapm}
\newcommand{\makerlapm}[2]{\hbox to 0pt{$\mathsurround=0pt #1{#2}$\hss}}

\newcommand{\makelapm}[2]{\hbox to 0pt{\hss$\mathsurround=0pt #1{#2}$\hss}}

\newcommand{\makeushort}[3]{%
  \setbox0=\hbox{$\mathsurround=0pt #2{#3}$}%
  \hbox to 1\wd0{\hss\underbar{\hbox to #1\wd0{\hss\box0\hss}}\hss}}

\def\makebigger#1#2#3{\scalebox{#1}{$\mathsurround=0pt #2{#3}$}}
\def\bigger#1#2{{\relax\mathpalette{\makebigger{#1}}{#2}}}

\def\scaleuphalf{1.0954}
\def\scaleupone{1.2}

\newcommand{\op}{^{\mathord{\text{\rm op}}}}

\newcommand{\defeq}{\mathrel{:=}}

\makeatletter

\def\newmop{\@ifstar{\@newmop m}{\@newmop o}}
\def\@newmop#1{\@ifnextchar[{\@@newmop #1}{\@@@newmop #1}}
\def\@@newmop#1[#2]{\@declmathop #1#2}
\def\@@@newmop#1#2{\expandafter\@declmathop\expandafter #1\csname #2\endcsname{#2}}

\makeatother

\newdir{ |}{{}*!/-5pt/@{|}}

\newmop{im}
\newmop{coim}
\newmop{dom}
\newmop{cod}
\newmop{id}
\newmop{Map}

\newmop{obj}
\newmop{arr}
\newmop{sq}
\newmop{norm}

\newmop{el}

\newmop{ev}

\newmop{fun}

\newmop{Ext}
\newmop{icon}
\newmop{pbk}
\newmop{icom}

\newmop{cell}
\newmop{cof}
\newmop{fib}

\newmop{diag}

\newmop*{colim}
\newmop*{holim}

\newmop[\lan]{lan}
\newmop[\ran]{ran}

\newmop*{cone}
\newmop*{cocone}
\newmop*{base}

\newcommand{\comma}{\mathbin{\downarrow}}

\makeatletter
\newcommand{\rotatemath}[2]{\rotatebox[origin=c]{180}{$\m@th #1{#2}$}}
\makeatother

\newmop[\pr]{pr}
\newmop[\dm]{d}

\newmop[\cpr]{in}

\newcommand{\pocorner}{\hbox to 8pt{{\vrule height8pt depth0pt width0.5pt}%
    \vbox to 8pt{{\hrule height0.5pt width7.5pt depth0pt}\vfill}}}
\newcommand{\poexcursion}{\save[]-<15pt,-15pt>*{\pocorner}\restore}
\newcommand{\pbcorner}{\vbox to 0pt{\kern 4pt\hbox to 0pt{\kern 4pt%
      \vbox{{\hrule height0.5pt width7.5pt depth0pt}}%
      {\vrule height8pt depth0pt width0.5pt}\hss}\vss}}
\newcommand{\pbexcursion}{\save[]+<5pt,-5pt>*{\pbcorner}\restore}

\newmop{cls}

\newcommand{\leib}[1]{\mathbin{\widehat{#1}}}

\newcommand{\category}[1]{\underline{\smash[b]{\text{\rm{#1}}}}}

\newcommand{\cattwo}{{\bigger{1.12}{\mathbbe{2}}}}
\newcommand{\catone}{{\bigger{1.16}{\mathbbe{1}}}}
\newcommand{\iso}{{\mathbb{I}}}

\makeatletter

\def\Del@Sym{{\bigger\scaleuphalf{\mathbbe{\Delta}}}}

\def\del@fn{\futurelet\del@next}
\def\del@dn{\def\del@next}

\def\parsedel@{%
  \ifx +\del@next \del@dn+{\Del@Sym_{\mathord{+}}}%
  \else \del@dn {\del@fn\parsedel@@}%
  \fi\del@next}

\def\parsedel@@{%
  \ifx\space@\del@next \expandafter\del@dn\space{\del@fn\parsedel@@}%
  \else\ifx [\del@next \del@dn[{\del@fn\parsedel@@@}%
  \else\ifx _\del@next \del@dn{\Delta}%
  \else\ifx ^\del@next \del@dn{\Delta}%
  \else \del@dn{\Del@Sym}%
  \fi\fi\fi\fi\del@next}

\def\parsedel@@@{%
  \ifx\space@\del@next \expandafter\del@dn\space{\del@fn\parsedel@@@}%
  \else\ifx t\del@next \del@dn t{\Del@Sym_\infty\del@fn\parsedel@@@@}%
  \else\ifx b\del@next \del@dn b{\Del@Sym_{-\infty}\del@fn\parsedel@@@@}%
  \else \del@dn{\errmessage{unexpected modifier}}%
  \fi\fi\fi\del@next}

\def\parsedel@@@@{%
  \ifx\space@\del@next \expandafter\del@dn\space{\del@fn\parsedel@@@@}%
  \else\ifx ]\del@next \del@dn]{}%
  \else \del@dn{\errmessage{expecting close of option block}}%
  \fi\fi\del@next}

\def\Del{\del@fn\parsedel@}

\makeatother

\newcommand{\Horn}{\Lambda}

\newcommand{\Cat}{\category{Cat}}

\newcommand{\sCat}{\sSet\text{-}\Cat}
\newcommand{\sSet}{\category{sSet}}
\newcommand{\qCat}{\category{qCat}}

\newcommand{\Adj}{\category{Adj}}
\newcommand{\Mnd}{\category{Mnd}}

\newcommand{\msSet}{\category{msSet}} 

\newcommand{\coCart}{\category{coCart}}

\newcommand{\dmod}[3]{\xymatrix@=1.25em{{#2} \ar[r]|\mid^{ {#1}} & {#3}}}

\newcommand{\pbshape}{{\mathord{\bigger\scaleupone\righthalfcup}}}

\newcommand{\face}{\delta}

\newcommand{\degen}{\sigma}

\newcommand{\fbv}[1]{\{{#1}\}}

\newcommand{\join}{\mathbin\star}

\newmop{dec}
\newmop{fatdec}
\newmop{slc}
\newmop{fatslc}

\newcommand{\slice}{/}
\newcommand{\slicel}[2]{\vphantom{#2}^{{#1}\slice{}}\mkern-2mu{#2}}
\newcommand{\slicer}[2]{{#1\mkern-1mu}_{{}\slice{#2}}}

\newmop{ir} 
\newmop{incl}

\newcommand{\laxslicer}[2]{{#1\mkern-1mu}_{{}{\sslice}{#2}}}

\newcommand{\nrv}{N}

\newcommand{\nrvhc}{\nrv}
\newcommand{\hN}{\nrv}
\newcommand{\gC}{\mathfrak{C}}

\newcommand{\boundary}{\partial}

\def\reedyfilt#1_#2{#1_{\leq #2}}
\newmop{sk}
\newmop{cosk}
\newmop{res}

\newmop{map}
\newmop{Ho}

\newmop[\pth]{path}
\newmop{cyl}

\newmop[\const]{c}

\newcommand{\Kan}{\category{Kan}}

\newmop{coll}
\newmop{wgt}

\newdir{ >}{{}*!/-7pt/@{>}}
\newdir{u(}{{}*!/-4pt/@^{(}}
\newdir{d(}{{}*!/-4pt/@_{(}}
\newdir{|>}{%
  !/4.5pt/@{|}*:(1,-.2)@^{>}*:(1,+.2)@_{>}*+@{}}

\makeatletter
\def\makeslashed#1#2#3#4#5{#1{\mathpalette{\sla@{#2}{#3}{#4}}{#5}}}

\def\@mathlower#1#2#3{\setbox0=\hbox{$\m@th#2#3$}\lower#1\ht0\box0}
\def\mathlower#1#2{\mathpalette{\@mathlower{#1}}{#2}}
\makeatother

\newcommand{\inc}{\hookrightarrow}

\newcommand{\tfib}{\twoheadrightarrow}

\newcommand{\xtfib}[1]{\xtwoheadrightarrow{#1}}
\newcommand{\xfibt}[1]{\xtwoheadleftarrow{#1}}

\newcommand{\longtwoheadrightarrow}{\mathrel{\mathord{-}\mkern-3mu\mathord\twoheadrightarrow}}

\newcommand{\we}{\xrightarrow{\mkern10mu{\smash{\mathlower{0.6}{\sim}}}\mkern10mu}}

\newcommand{\trvfib}{\stackrel{\smash{\mkern-2mu\mathlower{1.5}{\sim}}}\longtwoheadrightarrow}

\newcommand{\vsim}{\mathrel{\rotatebox{270}{$\sim$}}}

\newcommand{\To}{\Rightarrow}

\makeatletter

\def\tens@fn{\futurelet\tens@next}
\def\tens@dn{\def\tens@nextcont}
\newtoks\tens@toks
\def\addtotens@toks#1{\tens@toks=\expandafter{\the\tens@toks#1}}

\def\parsetens@@{%
    \ifx\space@\tens@next \expandafter\tens@dn\space{\tens@fn\parsetens@@}%
    \else\ifx ^\tens@next \tens@dn ^##1{\parsetens@procsep^\addtotens@toks{##1}%
      \tens@fn\parsetens@@}%
    \else\ifx _\tens@next \tens@dn _##1{\parsetens@procsep_\addtotens@toks{##1}%
      \tens@fn\parsetens@@}%
    \else\tens@dn{\ifx *\tens@last \else\addtotens@toks\egroup\fi\the\tens@toks}%
    \fi\fi\fi\tens@nextcont}

\def\parsetens@procsep#1{%
  \ifx *\tens@last \addtotens@toks{#1}\addtotens@toks\bgroup%
  \else\ifx \tens@last\tens@next \addtotens@toks,%
  \else \addtotens@toks\egroup\addtotens@toks\bgroup%
    \addtotens@toks\egroup\addtotens@toks{#1}\addtotens@toks\bgroup%
  \fi\fi\let\tens@last\tens@next}

\newcommand{\tn}[1]{\let\tens@last=*\tens@toks={#1}\tens@fn\parsetens@@}

\makeatother

\def\adjdisplay#1-|#2:#3->#4.{{%
    \xymatrix@R=0em@!C=2.5em{%
      *+[l]{#3} \ar@/_0.55pc/[rr]_-{#2} & {\bot} &
      *+[r]{#4}\ar@/_0.55pc/[ll]_-{#1}}}}

\def\adjdisplaytwo#1-|#2:#3->#4.{{%
\xymatrix@=1.2em{
      {#3}\ar@/_1.5ex/[rr]_-{#2}^-{}="one"
      & & {#4}
      \ar@/_1.5ex/[ll]_-{#1}^-{}="two"
      \ar@{}"one";"two"|{\bot}
    }}}

\def\tripleadjdisplay#1-|#2-|#3:#4->#5.{{%
\xymatrix@=2.4em{
{#4}\ar[r]|{#2} &
{#5} \ar@/_3ex/[l]_{#1}^{\bot} \ar@/^3ex/[l]_{\bot}^{#3}}
}}

\def\adjinline#1-|#2:#3->#4.{{#1}\dashv{#2}:#3\to #4}

\newcommand{\pent}[1]{
  \xybox{
    \POS (0,-15)*+{\a}="0",
         (-14,-5)*+{\b}="1",
         (-9,12)*+{\c}="2",
         (9,12)*+{\d}="3",
         (14,-5)*+{\e}="4"
    \POS"0" \ar "1"^{\labelstyle \ab}|{}="01"
    \POS"1" \ar "2"^{\labelstyle \bc}|{}="12"
    \POS"2" \ar "3"^{\labelstyle \cd}|{}="23"
    \POS"3" \ar "4"^{\labelstyle \de}|{}="34"
    \POS"0" \ar "4"_{\labelstyle \ae}|{}="04"
    \ifcase #1
    \POS"0" \ar "2"|{\labelstyle \ac}="02"
    \POS"0" \ar "3"|{\labelstyle \ad}="03"
    \POS"02";"1"**{}, ?(0.3) \ar@{=>} ?(0.7)^{\labelstyle \abc}
    \POS"03";"2"**{}, ?(0.25) \ar@{=>} ?(0.5)_{\labelstyle \acd}
    \POS"04";"3"**{}, ?(0.2) \ar@{=>} ?(0.4)_{\labelstyle \ade}
    \or
    \POS"1" \ar "3"|{\labelstyle \bd}="13"
    \POS"1" \ar "4"|{\labelstyle \be}="14"
    \POS"13";"2"**{}, ?(0.3) \ar@{=>} ?(0.7)_{\labelstyle \bcd}
    \POS"14";"3"**{}, ?(0.25) \ar@{=>} ?(0.5)_{\labelstyle \bde}
    \POS"04";"1"**{}, ?(0.25) \ar@{=>} ?(0.5)_{\labelstyle \abe}
    \or
    \POS"2" \ar "4"|{\labelstyle \ce}="24"
    \POS"0" \ar "2"|{\labelstyle \ac}="02"
    \POS"02";"1"**{}, ?(0.3) \ar@{=>} ?(0.7)^{\labelstyle \abc}
    \POS"04";"2"**{}, ?(0.2) \ar@{=>} ?(0.35)_{\labelstyle \ace}
    \POS"24";"3"**{}, ?(0.2) \ar@{=>} ?(0.6)^{\labelstyle \cde}
    \or
    \POS"1" \ar "3"|{\labelstyle \bd}="13"
    \POS"0" \ar "3"|{\labelstyle \ad}="03"
    \POS"04";"3"**{}, ?(0.2) \ar@{=>} ?(0.4)_{\labelstyle \ade}
    \POS"13";"2"**{}, ?(0.3) \ar@{=>} ?(0.7)_{\labelstyle \bcd}
    \POS"03";"1"**{}, ?(0.25) \ar@{=>} ?(0.5)^{\labelstyle \abd}
    \or
    \POS"2" \ar "4"|{\labelstyle \ce}="24"
    \POS"1" \ar "4"|{\labelstyle \be}="14"
    \POS"24";"3"**{}, ?(0.2) \ar@{=>} ?(0.6)^{\labelstyle \cde}
    \POS"04";"1"**{}, ?(0.25) \ar@{=>} ?(0.5)_{\labelstyle \abe}
    \POS"14";"2"**{}, ?(0.25) \ar@{=>} ?(0.5)^{\labelstyle \bce}
    \else\fi
  }
}

\newcommand{\pentofpent}[1]{
  \def\baselen{#1}
  \begin{xy}
    0;<\baselen,0mm>:
    *{\xybox{
        \POS(0,-4)*[o]{\pent 0}="zero"
        \POS(16,40)*[o]{\pent 3}="three"
        \POS(72,40)*[o]{\pent 1}="one"
        \POS(88,-4)*[o]{\pent 4}="four"
        \POS(44,-36)*[o]{\pent 2}="two"
        \ar@<1ex>"zero";"three"^-{\objectstyle\abcd}
        \ar@<1ex>"three";"one"^-{\objectstyle\abde}
        \ar@<1ex>"one";"four"^-{\objectstyle\bcde}
        \ar@<-1ex>"zero";"two"_-{\objectstyle\acde}
        \ar@<-1ex>"two";"four"_-{\objectstyle\abce}
        \ar@{=>}(44,-5);(44,+15)^{\objectstyle\abcde}
     }}
  \end{xy}
}


%% file: intro.tex

\section{Introduction}

Famously the category $\Cat$ of small categories is not a topos because, among other things, it fails to be locally cartesian closed. A finitely complete category $\eE$ is \emph{locally cartesian closed} just when each slice category $\eE_{/B}$ is cartesian closed, or equivalently, which the pullback functor associated to any morphism $f \colon A \to B$ admits a right adjoint (as well as a left adjoint given by composition with $f$):
\[
\xymatrix@C=3em{ \eE_{/B} \ar[r]|{f^*} & \eE_{/A} \ar@/_3ex/[l]_{\Sigma_f}^\perp \ar@/^3ex/[l]^{\Pi_f}_\perp}\]

In the case $\eE=\Cat$, those functors $f$ for which the pullback functor $f^*$ does admit a right adjoint $\Pi_f$ are called \emph{exponentiable} and have been characterized by Conduch\'{e} \cite{Conduche:1972as}. Famously, 
  \begin{enumerate}[label=(\roman*)]
  \item\label{itm:desiderata-i} All cartesian and cocartesian fibrations $p\colon E\to B$ of 1-categories are
    exponentiable.
  \item\label{itm:desiderata-ii} If $p\colon E\to B$ is a cocartesian fibration and $q\colon F\to E$
    is a cartesian fibration then the pushforward $\Pi_p(q\colon F\to E)$
    is also a cartesian fibration, and the dual result holds when the cartesian and cocartesian fibrations are interchanged. 
  \end{enumerate}

  In \cite{Lurie:2017ha} Lurie established $\infty$-categorical analogues of
these results for quasi-categories.\footnote{A general characterization of the exponentiable functors between quasi-categories, while not the focus of our interest here, can be found in  \cite[\S B.3]{Lurie:2017ha} or \cite{AyalaFrancis:2017fo}.}
 Subsequent authors, for instance Barwick
  and Shah~\cite{Barwick:2016aa}, have stressed the importance of these results to the theory and practice of $\infty$-categories, and we have further applications in mind. In \cite[\S12.3]{RiehlVerity:2022eo}, we use this result to prove that \emph{modules} between quasi-categories admit all right and left extensions. It follows that the question of existence of pointwise right and left Kan extensions can be reduced to the existence of certain limits and colimits. For those results, it is useful to have a somewhat more refined version of these results than is easily found in the literature---see especially Theorem \ref{thm:cosmological-pushforward} and Corollary \ref{cor:cosmological-pushforward}---which is the motivation for the present exposition


A \emph{cocartesian fibration} of quasi-categories is an isofibration\footnote{In the Joyal model structure on simplicial sets, we refer to the fibrations between fibrant objects (the quasi-categories) as \emph{isofibrations} because they have a lifting property for isomorphisms analogous to that for the isofibrations in classical category theory.} $p \colon \qE \tfib \qB$ whose fibers depend covariantly functorially on $\qB$. In the simplest non-trivial case, when $\qB = \Del^1$, the data is given by a pair of quasi-categories $\qE_0$ and $\qE_1$ together with a functor $\qE_0 \to \qE_1$.  In general, the \emph{comprehension construction} of \cite{RiehlVerity:2017cc} ``straightens'' $p \colon \qE \tfib \qB$ into a simplicial functor $c_p\colon\gC\qB \to \qCat$ that sends each vertex $b \in \qB$ to the fiber $\qE_b$. The domain category appearing here is the \emph{homotopy coherent realization} of the quasi-category $\qB$, a cofibrant simplicial category\footnote{The simplicial categories that are cofibrant in the Bergner model structure are precisely the \emph{simplicial computads} that are freely generated by their non-degenerate ``atomic'' $n$-arrows for each $n \geq 0$, admitting no non-trivial factorizations; see Definition \ref{defn:simplicial-computad}.} that indexes $\qB$-shaped homotopy coherent diagrams. At the level of objects and 1-arrows $f \colon a \to b$ in $\qB$, the comprehension construction is defined     by lifting the 1-arrow $f$ to a $p$-\emph{cocartesian 1-arrow} with codomain $\qE$:
\[
        \begin{xy}
      0;<1.4cm,0cm>:<0cm,0.75cm>::
      *{\xybox{
          \POS(1,0)*+{1}="one"
          \POS(0,1)*+{1}="two"
          \POS(3,0.5)*+{\qB}="three"
          \ar@{=} "one";"two"
          \ar@/_5pt/ "one";"three"_{b}^(0.1){}="otm"
          \ar@/^10pt/ "two";"three"^{a}_(0.5){}="ttm"|(0.325){\hole}
          \ar@{=>} "ttm"-<0pt,7pt> ; "otm"+<0pt,10pt> ^(0.3){f}
          \POS(1,2.5)*+{\qE_{b}}="one'"
          \POS(1,2.5)*{\pbcorner}
          \POS(0,3.5)*+{\qE_{a}}="two'"
          \POS(0,3.6)*{\pbcorner}
          \POS(3,3)*+{\qE}="three'"
          \ar@/_5pt/ "one'";"three'"_{\ell^\qE_{{b}}}
          \ar@/^10pt/ "two'";"three'"^{\ell^\qE_{{a}}}_(0.55){}="ttm'"
          \ar@{->>} "one'";"one"_(0.325){p_{b}}
          \ar@{->>} "two'";"two"_{p_{a}}
          \ar@{->>} "three'";"three"^{p}
          \ar@{..>} "two'";"one'"_*!/^2pt/{\scriptstyle \qE_f}
          \ar@{..>}@/_10pt/ "two'";"three'"^(0.44){}="otm'"
          \ar@{=>} "ttm'"-<0pt,4pt> ; "otm'"+<0pt,4pt> ^(0.3){\ell^\qE_{{f}}}
        }}
    \end{xy}  \]
Together, this data defines a lax cocone $\ell^\qE$  under the comprehension functor $c_p$ with nadir $\qE$ the data of which is given by a functor $\ell^\qE\colon \gC[\qB \join \Del^0] \to \qCat$ that restricts along $\gC\qB \inc \gC[\qB\join \Del^0]$ to $c_p$. In fact, $\ell^\qE$ is a colimit cocone:

{
\renewcommand{\thethm}{\ref{cor:total-space-as-oplax-colimit}}
\begin{cor}
The domain of a cocartesian fibration $p \colon \qE \tfib \qB$ is equivalent to the oplax colimit of the associated comprehension functor $c_p \colon \gC\qB \to \qCat$, with colimit cocone:
            \begin{equation*}
    \xymatrix@=1.5em{
      & {\catone+\gC{\qB}}\ar@{^(->}[dl]\ar[dr]^{\langle \qE,c_p \rangle} & \\
      {\gC[\qB\join\Del^0]}\ar[rr]_-{\ell^\qE} && {\qCat}
    }
  \end{equation*}
  \end{cor}
  \addtocounter{thm}{-1}
  }

In particular the domain $\qE$ of the cocartesian fibration $p$ can be recovered up to equivalence as the oplax colimit of a the comprehension functor $c_p \colon \gC{\qB} \to \qCat$. Gepner, Haugseng, and Nikolaus, who obtain a similar result to Corollary \ref{cor:total-space-as-oplax-colimit} as one of the main theorems of \cite{GepnerHaugsengNikolaus:2017lc}, interpret this result as a proof that ``Lurie's unstraightening functor is a model for the $\infty$-categorical analogue of the Grothendieck construction.''\footnote{Unfortunately, the assignment of the terms ``oplax colimit'' and ``lax colimit'' given in \cite[2.8]{GepnerHaugsengNikolaus:2017lc} is opposite to the one used here. The standard convention in 2-category theory is that the 2-cell component of an oplax natural transformation is parallel to its 1-cell components, while these 2-cells are reversed in a lax natural transformation. A lax cocone is then a lax natural transformation whose codomain is a constant diagram. Confusingly, due to the principle that a $W$-weighted colimit in an enriched category coincides with a $W$-weighted limit in the opposite category, oplax colimits represent lax cocones under a diagram.
}
Their methodology is quite different from ours, constructing oplax colimits directly at the quasi-categorical level, whereas our comprehension construction enables us to work at the level of simplicial categories and functors. The comprehension functor $c_p \colon \gC{\qB} \to \qCat$ can be used to define a ``straightening'' of the pullback of $p$ along any generalized element $b \colon X \to \qB$, even in the case where $X$ is not a quasi-category simply by restricting the comprehension functor (and its lax cocone) along $b$. We derive Corollary \ref{cor:total-space-as-oplax-colimit} as a special case of our first main theorem, which proves that the fiber $E_b$ is equivalent to the oplax colimit of this straightened diagram.

{
\renewcommand{\thethm}{\ref{thm:pullback-equivalence}}
\begin{thm}
For any cocartesian fibration $p \colon \qE \tfib \qB$ and any $b \colon X \to \qB$, the comprehension cocone induces a canonical map over $\qE$ from the oplax colimit of the diagram
\[ \gC{X} \xrightarrow{\gC{b}} \gC{\qB} \xrightarrow{c_p} \qCat\]
to the fiber 
\[
\xymatrix{ E_b \ar[d] \ar[r] \pbexcursion & \qE \ar@{->>}[d]^p \\ X \ar[r]_b & \qB}
\]
and this map is a natural weak equivalence in the Joyal model structure.
\end{thm}
\addtocounter{thm}{-1}
}

The canonical natural transformation of Theorem \ref{thm:pullback-equivalence} defines a natural Joyal equivalence relating the pullback functor $p^*$ to a functor $\tilde{p}^*$ defined by forming oplax colimits of restrictions of the comprehension cocone:
  \[ \xymatrix{ \sSet_{/\qB} \ar@/^3ex/[r]^{\tilde{p}^*} \ar@/_3ex/[r]_{p^*} \ar@{}[r]|{\Downarrow\gamma} & \sSet_{/\qE}}\] 
Both functors $p^*$ and $\tilde{p}^*$ are left Quillen with respect to the sliced Joyal model structures, admitting right Quillen adjoints:

{
\renewcommand{\thethm}{\ref{prop:quillen-exponentiation-adjunctions}}
  \begin{prop}If $p \colon \qE \tfib \qB$ is a cocartesian fibration, the adjunctions
  \[
\begin{tikzcd} \sSet_{/\qE} \arrow[r, bend right, "{\Pi}_p"'] \arrow[r, phantom, "\bot"] & \sSet_{/\qB} \arrow[l, bend right, "{p}^*"'] & \text{and} &  \sSet_{/\qE} \arrow[r, bend right, "\tilde{\Pi}_p"'] \arrow[r, phantom, "\bot"] & \sSet_{/\qB} \arrow[l, bend right, "\tilde{p}^*"']
\end{tikzcd}
\]are Quillen with respect to the sliced Joyal model structures.
  \end{prop}
  \addtocounter{thm}{-1}
  }
By taking mates, there is a canonical natural transformation $\hat{\gamma} \colon \Pi_p \To \tilde\Pi_p$ whose component at any isofibration $q \colon \qF \tfib \qE$ is an equivalence. In this way we obtain an alternate model $\tilde{\Pi}_p$ for the pushforward functor that is more easily  understood: at an isofibration $q \colon \qF \tfib \qE$, $\tilde{\Pi}_pq$ is the pullback along the comprehension cocone of the induced map between lax slices induced by whiskering with $q$:\footnote{The precise meaning of this notation, involving slices of the homotopy coherent nerve of $\qCat$ regarded as a 2-complicial set, is explained in Lemma \ref{lem:lax-slice-whiskering-comp}.}
\[ 
\xymatrix{ \tilde{\Pi}_p(\qF \xtfib{q} \qE) \pbexcursion \ar@{->>}[d] \ar[r] & {\qqCat_2}_{\sslice{\qF}} \ar@{->>}[d]^{q\circ -} \\ \qB \ar[r]_-{\ell^\qE} & {\qqCat_2}_{\sslice{\qE}}}\]
 To prove Proposition \ref{prop:quillen-exponentiation-adjunctions}, we show that the ``whiskering with $q$'' map is an isofibration. This establishes the quasi-categorical analogue of desiderata \ref{itm:desiderata-i} above. We then show further that if $q \colon \qF \tfib \qE$ is a cartesian fibration, then the ``whiskering with $q$'' map has a certain right horn lifting property, thereby proving the quasi-categorical analogue of desiderata \ref{itm:desiderata-ii}:

{
\renewcommand{\thethm}{\ref{cor:cartesian-workhorse}}
  \begin{cor}
   If $p \colon \qE \tfib \qB$ is a cocartesian fibration  and $q \colon \qF \tfib \qE$ is a cartesian fibration between quasi-categories, then the pushforward \[ \Pi_p(q \colon\qF\tfib \qE) \tfib \qB\] is a cartesian fibration between quasi-categories. 
  \end{cor}
    \addtocounter{thm}{-1}
  }
  
We show also that the pullback and pushforward functors along a cocartesian fibration preserve the accompanying class of \emph{cartesian functors} between cartesian fibrations. These results are summarized in the following theorem:

{
\renewcommand{\thethm}{\ref{thm:cosmological-pushforward}}
  \begin{thm}
  For any cocartesian fibration $p\colon \qE \tfib \qB$  between quasi-categories, the pullback-pushforward adjunction restricts to define an adjunction
\[
\begin{tikzcd}[row sep=large] \qCat_{/\qE} \arrow[r, bend right=15, "\Pi_p"'] \arrow[r, phantom, "\bot"] & \qCat_{/\qB} \arrow[l, bend right=15, "p^*"'] \\  \Cart(\qCat)_{/\qE} \arrow[r, bend right=15, "\Pi_p"'] \arrow[r, phantom, "\bot"] \arrow[u, hook] & \Cart(\qCat)_{/\qB} \arrow[l, bend right=15, "p^*"'] \arrow[u, hook]
\end{tikzcd}
\]
    \end{thm}
        \addtocounter{thm}{-1}
  }

As an immediate corollary, we construct ``exponentials'' whose exponents are either cartesian or cocartesian, justifying the appelation ``exponentiable'' for these maps, and prove:

{
\renewcommand{\thethm}{\ref{prop:cartesian-workhorse}}
\begin{prop} 
 If $p \colon \qE \tfib \qB$ is a cocartesian fibration and $q \colon \qF \tfib \qB$ is a cartesian fibration, then  \[ (q \colon \qF \tfib \qB)^{p \colon \qE \tfib \qB}\]  is a cartesian fibration.  \end{prop}
     \addtocounter{thm}{-1}
  }
  
The final two sections of this paper supply some first applications of these results. As the comprehension construction reveals, cartesian fibrations over $\qB$ encode functors $\qB\op\to\qqCat$ valued in the (large) quasi-category of small quasi-categories. In ordinary category theory it is well-known that for any small category $B$ and complete and cocomplete category $C$, the forgetful functor $C^B \to C^{\ob{B}}$ that carries a diagram to its $\ob{B}$-indexed family of objects admits both left and right adjoints, given by left and right Kan extension, and is moreover \emph{monadic} and \emph{comonadic}. Informally, this means that $B$-indexed diagrams can be understood as ``algebras'' or as ``coalgebras'' for a monad or comonad acting on the category of $\ob{B}$-indexed families of objects.

The corresponding result for quasi-categories will be proven in  a sequel to this paper\ifnewcut\else ~\cite{RiehlVerity:2017ts}\fi, but here we demonstrate that the analogous result holds for cartesian fibrations, using a version of Beck's monadicity theorem for quasi-categories proven in \cite{RiehlVerity:2012hc}. Writing $\Cart_{/\qB}$ for the large quasi-category of cartesian fibrations and cartesian functors over $\qB$, we prove:

{
\renewcommand{\thethm}{\ref{thm:comonadic-cart}}
\begin{thm} The forgetful functor 
\[ u\colon \Cart_{/\qB} \longrightarrow \Cart_{/\ob\qB} \cong \prod_{\ob\qB}\qqCat\]
is comonadic and hence also monadic.
\end{thm}
     \addtocounter{thm}{-1}
  }
In a further sequel, we will use this monadicity to establish an equivalence between $\Cart_{/\qB}$ and $\qqCat^{\qB\op}$, both quasi-categories being monadic over $\prod_{\ob\qB}\qqCat$.

Here we include another application of the monadicity of Theorem \ref{thm:comonadic-cart}. Using our analysis of limits and colimits in quasi-categories defined as homotopy coherent nerves in \cite{RiehlVerity:2020oc}, we prove:

{
\renewcommand{\thethm}{\ref{thm:basic-groupoidal-reflection}}
\begin{thm} The inclusion $\qKan\inc\qqCat$ admits both left and right adjoints 
\[ 
\xymatrix@C=3em{\qKan \ar@{^(->}[r]^-\perp_\perp & \qqCat \ar@/_3ex/[l]_{\mathrm{invert}} \ar@/^3ex/[l]^\core}\]
and is monadic and comonadic.
\end{thm}
     \addtocounter{thm}{-1}
  }
  
  The right adjoint here is the familiar functor that takes a quasi-category to its maximal Kan complex core, while the left adjoint is a somewhat more delicate ``groupoidal reflection'' functor. Our final result establishes an analogous groupoidal reflection for cartesian fibrations into the subcategory of \emph{groupoidal cartesian fibrations}, whose fibers are Kan complexes rather than quasi-categories.
  
  {
\renewcommand{\thethm}{\ref{thm:groupoidal-reflection-cartesian}}
  \begin{thm}
There is a left adjoint to the inclusion 
\[
\begin{tikzcd}[sep=large] \Cart\gr_{/\qB} \arrow[r, bend right] \arrow[r, phantom, "\bot"] & \Cart_{/\qB} \arrow[l, bend right, "\mathrm{invert}"']
\end{tikzcd}
\]defining the reflection of a cartesian fibration into a groupoidal cartesian fibration.
\end{thm}
     \addtocounter{thm}{-1}
  }
  
  All of the results mentioned above have duals with cocartesian and cartesian fibrations interchanged. The comprehension functor associated to a cartesian fibration is contravariant and its domain is recovered as the lax colimit of this diagram. It is to avoid this contravariance that we choose to focus the bulk of our presentation on the case of cocartesian fibrations.

This paper is organized as follows. \ifnewcut\else In \S\ref{sec:weighted}, we introduce oplax colimits through the general mechanism of weighted colimits. We prove that oplax weights are \emph{flexible}, which implies that the oplax weighted colimit functor is equivalence-invariant. We also review the collage construction, which allows us to construct flexible weights by instead specifying the shape of their corresponding cocones. In particular, a lax cocone of shape $X$ is indexed by the homotopy coherent realization of the join $X\join\Del^0$.

In \S\ref{sec:cocartesian}, we review some basic aspects of the theory of cocartesian fibrations between quasi-categories. We introduce a quasi-categorical version of the collage construction just discussed, and prove that the quasi-categorical collage of a functor $f \colon \qA \to \qB$ defines a cocartesian fibration over $\Del^1$ that models the oplax colimit of $f$.

In \S\ref{sec:comprehension}, we review the comprehension construction from \cite{RiehlVerity:2017cc}, devoting somewhat more attention to the lax cocones that are the focus of much of the work here.\fi 
\ifnewadd In \S\ref{sec:background}, we provide background material on oplax colimits, cocartesian fibrations, and the comprehension construction from \cite{RiehlVerity:2017cc}. \fi  
Then in \S\ref{sec:pullback}, we prove that pullback along a cocartesian fibration can be modeled as a oplax colimit of a restriction of the comprehension functor. 

The corresponding results for the pushforward functor, including in particular \ref{itm:desiderata-i} and \ref{itm:desiderata-ii}, are then proven in \S\ref{sec:powerful}. The oplax colimits defining the functor $\tilde{p}^*$ in \S\ref{sec:pullback} are properly understood as a variety of $(\infty,2)$-categorical colimits. Consequently, the description of the corresponding right adjoint $\tilde\Pi_p$ involves an $(\infty,2)$-categorical cocone construction, instantiated by forming the slice of 2-complicial set over a vertex. As we explain in \S\ref{sec:2complicial}, 2-complicial sets are simplicial sets in which certain simplices are marked as ``thin.'' This notion is not as unfamiliar as it may seem at first: Kan complexes are precisely the 0-complicial sets while quasi-categories correspond to 1-complicial sets.

In \S\ref{sec:monadicity}, we consider the forgetful functor $\Cart(\qCat)_{/\qB} \to \Cart(\qCat)_{/\ob\qB}$ and construct left and right \emph{biadjoints}: quasi-categorically enriched functors equipped with a natural equivalence of function complexes encoding the adjoint transpose relation. Such data descends to an adjunction between the quasi-categorical cores of these quasi-categorically enriched categories. We then review the monadicity theorem from \cite{RiehlVerity:2012hc} and apply it to prove that this forgetful functor is monadic and comonad as a map between large quasi-categories.

To say that the functor $\Cart_{/\qB} \to \Cart_{/\ob\qB}\cong \prod_{\ob\qB}\qqCat$ is monadic is to say that $\Cart_{/\qB}$ may be recovered as the quasi-category of \emph{homotopy coherent algebras} for a homotopy coherent monad acting on $\prod_{\ob\qB}\qqCat$. In \S\ref{sec:groupoidal-reflection}, we show that $\Cart\gr_{/\qB}$ is similarly the quasi-category of homotopy coherent algebras for the restriction of this homotopy coherent monad along the inclusion $\prod_{\ob\qB}\qKan\inc\prod_{\ob\qB}\qqCat$. We then show that this characterization allows us  to construct the groupoidal reflection functor as a lift of the groupoidal reflection functor $\qqCat\to\qKan$.

This paper is a continuation of a series of papers that redevelop the foundations of $(\infty,1)$-category theory  \cite{RiehlVerity:2012tt, RiehlVerity:2012hc, RiehlVerity:2013cp,RiehlVerity:2015fy,RiehlVerity:2015ke,RiehlVerity:2017cc,RiehlVerity:2018rq,RiehlVerity:2020oc}, the results of which are referenced as I.x.x.x, \ldots, and VIII.x.x.x respectively. However, we deploy relatively few of the tools developed in our previous work to prove the theorems appearing here, and when we do reference prior results, we typically restate them in considerably less generality. Many of the results from previous work recalled here --- for instance Theorem \ref{thm:comprehension} --- are proven in the more abstract setting of any $\infty$-\emph{cosmos}, while in the present manuscript we consider only a single example: the quasi-categorically enriched category $\qCat$ of quasi-categories. As we do not need this notion, we do not recall any specifics here.\footnote{This said, however, we cannot resist appealing to $\infty$-cosmic techniques in a handful of our proofs.} \ifnewadd While this paper was in press, the book \cite{RiehlVerity:2022eo} was published, so in the final version of the present manuscript we have cut a few proofs and instead refer to the corresponding results X.x.x.x that now appear there. \fi




\subsection{Acknowledgments}

The authors are grateful for support from the National Science Foundation (DMS-1551129) and from the Australian Research Council (DP160101519, DP190102432). This work was commenced when the second-named author was visiting the first at Harvard and then at Johns Hopkins and completed while the first-named author was visiting the second at Macquarie. We thank all three institutions for their assistance in procuring the necessary visas as well as for their hospitality. \ifnewadd The final journal version was prepared during a visit supported by the Johns Hopkins Frontier Award, the NSF (DMS-2204304), the ARO MURI (W911NF-20-1-0082), and the US Army DEVCOM Indo-Pacific Fundamental Research Collaboration Opportunities (FA520923C0004). \fi


%% file: background.tex

\section{Background}\label{sec:background}

\ifnewadd
In \S\ref{ssec:weighted}, we introduce oplax colimits through the general mechanism of weighted colimits. We prove that oplax weights are \emph{flexible}, which implies that the oplax weighted colimit functor is equivalence-invariant. We also review the collage construction, which allows us to construct flexible weights by instead specifying the shape of their corresponding cocones. In particular, a lax cocone of shape $X$ is indexed by the homotopy coherent realization of the join $X\join\Del^0$.

In \S\ref{ssec:cocartesian}, we review some basic aspects of the theory of cocartesian fibrations between quasi-categories. We introduce a quasi-categorical version of the collage construction and prove that the quasi-categorical collage of a functor $f \colon \qA \to \qB$ defines a cocartesian fibration over $\Del^1$ that models the oplax colimit of $f$.

In \S\ref{ssec:comprehension}, we review the comprehension construction, devoting somewhat more attention to the lax cocones that are the focus of much of the work here.
\fi 

\subsection{Oplax colimits in simplicial categories}\label{ssec:weighted}

Our aim in this section is to define the \emph{oplax colimit} of a homotopy coherent diagram $\gC{X} \to \SSet$ indexed by the homotopy coherent realization of a simplicial set $X$. Oplax colimits are introduced as particular \emph{weighted colimits}, where the weights in question are simplicial functors that describe the shape of lax cocones. \ifnewcut\else In \S\ref{sec:simp-collage}, we review the general rubric of weighted colimits and explain how these cocone shapes may be presented as simplicial \emph{collages}. In \S\ref{sec:flexible-weight}, we highlight a special class of \emph{flexible weights} that have useful homotopical properties. Finally, in \S\ref{sec:oplax}, we define the weights for oplax colimits as collages and observe that such weights are flexible, a fact that will be exploited in our future work.

\fi Some of this material was previously discussed in \S\refVII{sec:flexible}, where the ``oplax'' weights were called ``pseudo'' weights. See Remark \ref{rmk:oplax-vs-pseudo} for an explanation of this contrast in nomenclature. \ifnewcut\else Collages for weights for limits made an appearance in \S\refVII{ssec:collage} and we refer the proofs of a few of the results appearing below to there, but we reintroduce this construction here to clarify the details in the dual case and because we will require a more extensive analysis of collages than we did in \cite{RiehlVerity:2018rq}.\fi

\ifnewcut\else
\subsection{Weighted colimits and collages}\label{sec:simp-collage}
\fi 

In a simplicially enriched category, the appropriately general notion of colimit allows for the specification of any particular ``shape'' of cone under the diagrams being considered. This specification is given by a simplicial functor referred to as a \emph{weight} for the colimit. 

\begin{defn}[weights for simplicial colimits]\label{defn:simp-weight}
  Suppose $\eD$ is a small simplicial category, which we think of as a
  diagram shape. Then a \emph{weight} on $\eD$ is a simplicial functor $W\colon
  \eD\op\to\SSet$. For any diagram $F\colon \eD\to\eK$ valued in a
  simplicial category $\eK$, a \emph{$W$-cocone} with \emph{nadir} an object
  $e\in \eK$ is a simplicial natural transformation $\iota\colon W\to
  \Fun_{\eK}(F(-),e)$. We say that the $W$-cocone $\iota$ \emph{displays $e$ as
    a $W$-colimit of $F$} if and only if for all objects $e'\in\eK$ the
  simplicial map 
  \begin{equation*}
    \xymatrix@R=0em@C=5em{
      \Fun_{\eK}(e,e') \ar[r]^-{\cong} & \Fun_{\SSet^{\eD\op}}(W,\Fun_{\eK}(F(-),e'))
    }
  \end{equation*}
  given by pre-composition with $\iota$ is an isomorphism. 
 
Many notations are common for the nadir of a weighted colimit cone; here we write $\colim^W\!F$ for the
  colimit of $F$ weighted by $W$. When these exist for all weights and diagrams in
  $\eK$ then $\colim$ defines a simplicial bifunctor that is cocontinuous in both variables:
  \begin{equation*}
    \xymatrix@R=0em@C=7em{
      \SSet^{\eD\op}\times\eK^{\eD}\ar[r]^-{\colim} & \eK
    }
  \end{equation*} 
\end{defn}

A simplicial functor $W\colon\eD\op\to\SSet$ may otherwise be described as
comprising a family of simplicial sets $\{Wd\}_{d\in\obj(\eD)}$ along with right
actions of the hom-spaces of $\eD$
\begin{equation}\label{eq:action-of-weight}
  \xymatrix@R=0em@C=6em{
    Wd'\times\Fun_{\eD}(d,d')\ar[r]^-{*} & Wd
  }
\end{equation}
 which satisfy 
axioms with respect to the identities and composition of $\eD$. This
description leads us to define a simplicially enriched category $\coll(W)$,
called the \emph{collage} of $W$.

\begin{defn}[collages]\label{defn:collage}
  For any weight $W \colon \eD\op\to\SSet$, the \emph{collage} of $W$ is a
  simplicial category $\coll(W)$ that contains $\eD$ as a full simplicial
  subcategory along with precisely one extra object $\top$ whose endomorphism
  space is the point. The function complexes $\Fun_{\coll(W)}(\top,d)$ are all taken
  to be empty and we define:
  \begin{equation*}
    \Fun_{\coll(W)}(d,\top) \defeq Wd\mkern40mu \text{for objects $d\in\eD$.}
  \end{equation*}
  The composition operations between hom-spaces in $\eD$ and those with
  codomain $\top$ are given by the actions depicted in~\eqref{eq:action-of-weight}.
\end{defn}

In the statement of the following result,
$\sSet^{\eD\op}$ denotes the underlying category of the simplicially
  enriched category $\SSet^{\eD\op}$.

\begin{prop}[{collage adjunction, \refVII{prop:collage-adjunction}}]\label{prop:collage-adjunction}$\quad$
\begin{enumerate}[label=(\roman*)]
\item\label{itm:collage-fun} The collage construction defines a fully faithful functor 
  \begin{equation*}
    \xymatrix@R=0em@C=6em{
      \sSet^{\eD\op}\ar[r]^-{\coll} & \prescript{\catone+\eD/}{}{\sCat}
    }
  \end{equation*}
  from the category of $\eD$-indexed weights to the category of simplicial categories under $\catone+\eD$ whose essential image is comprised of those $\langle
  e,F\rangle\colon \catone+\eD\to\eK$ that are bijective on objects, fully faithful
  when restricted to $\eD$ and $\catone$, and have the property that there are no maps in $\eK$ from $e$ to the image of $F$.
\item\label{itm:collage-adj} The collage functor admits a right adjoint, which carries a pair 
$\langle e,F
  \rangle\colon\catone+\eD\to\eK$ to the weight $\Fun_{\eK}(F(-),e)\colon\eD\op\to\SSet$.
\[ \adjdisplay \coll -| \wgt : \prescript{\catone+\eD/}{}{\sCat} ->      \sSet^{\eD\op}. \ifnewadd\qed\fi \] 
\end{enumerate}
\end{prop}

  This adjunction has a useful and important interpretation:
  
  \begin{cor}[{\refVII{cor:collage-bijection}}]\label{cor:collage-bijection} The collage   $\coll(W)$ of a weight
 realises the shape of $W$-cocones, in the sense that simplicial functors 
\[  G \colon \coll(W) \longrightarrow \eK\]
stand in bijection to $W$-cocones under the diagram $G\vert_{\eD}$ with nadir $G(\top)$. \qed
\end{cor}

\ifnewadd 
We record some basic properties about weighted colimits, collages, and left Kan extensions for later use. For proof see \cite[\S 2.1]{RiehlVerity:2021ce-v1}.

\begin{lem}\label{lem:lan-of-weight}
  For any simplicial functor $I \colon \eD \to \eC$ and weight $W \colon
  \eD\op\to\SSet$:
  \begin{enumerate}[label=(\roman*)]
        \item For any diagram $G \colon \eC\to\eK$, we have an isomorphism
  \begin{equation*}
    \colim\nolimits^W\!GI \cong \colim\nolimits^{\lan_{I}W}\!G
  \end{equation*}
  where the colimit on one side exists if and only if the one on the other
  does. 
  \item  Left Kan extension of $W$ along $I$ gives rise to a pushout square in the category of simplicial categories and simplicial functors:
  \begin{equation*}
    \xymatrix@=2em{
      {\catone+\eD}\ar[r]^{\catone+I}\ar@{^(->}[d] &
      {\catone+\eC}\ar@{^(->}[d] \\
      {\coll(W)}\ar[r] & {\coll(\lan_{I}W)}\poexcursion
    } \ifnewadd\qed\fi
  \end{equation*}
  \end{enumerate}
\end{lem}
\fi 
\ifnewcut\else
\begin{lem}\label{lem:lan-of-weight}
  For any simplicial functor $I \colon \eD \to \eC$, weight $W \colon
  \eD\op\to\SSet$, and diagram $G \colon \eC\to\eK$, we have an isomorphism
  \begin{equation*}
    \colim\nolimits^W\!GI \cong \colim\nolimits^{\lan_{I}W}\!G
  \end{equation*}
  where the colimit on one side exists if and only if the one on the other
  does. 
\end{lem}
\begin{proof}
  Simplicial left Kan extension provides an adjunction
  \begin{equation*}
    \adjdisplay \lan_{I} -| -\circ I : \SSet^{\eC\op} -> \SSet^{\eD\op}.
  \end{equation*}
  In particular
  \[
    \Fun_{\SSet^{\eC\op}}(\lan_IW,\Fun_{\eK}(G(-),e)) \cong
    \Fun_{\SSet^{\eD\op}}(W,\Fun_{\eK}(GI(-),e))
  \]
  which shows that $\colim^{\lan_IW}G$ and $\colim^W\!GI$  have
  the same defining universal property.
\end{proof}
\begin{lem}\label{lem:lan-collage-pushout}
  Left Kan extension of $W\colon \eD\op\to\SSet$ along a simplicial functor
  $I\colon \eD\to\eC$ gives rise to a pushout square
  \begin{equation*}
    \xymatrix@=2em{
      {\catone+\eD}\ar[r]^{\catone+I}\ar@{^(->}[d] &
      {\catone+\eC}\ar@{^(->}[d] \\
      {\coll(W)}\ar[r] & {\coll(\lan_{I}W)}\poexcursion
    }
  \end{equation*}
  in the category of simplicial categories and simplicial functors.
\end{lem}
\begin{proof}
  By the defining universal property, a simplicial functor whose domain is the
  pushout of $\catone+\eD\inc\coll(W)$ along $\catone + I$ and whose codomain
  is $\eK$ is given by a pair of functors
  \[
    \gamma\colon \coll(W)\to \eK \qquad \mathrm{and} \qquad
    \langle e,G\rangle \colon\catone+\eC \to \eK.
  \]
  By Corollary \ref{cor:collage-bijection}, the simplicial functor $\gamma$
  represents a $W$-cocone in $\eK$ with nadir $e$ under the diagram
  $GI\colon\eD\to\SSet$. By Lemma \ref{lem:lan-of-weight}, such data
  equivalently describes a $\lan(W)$-shaped cocone under the diagram $G$ with
  nadir $e$. Applying Corollary \ref{cor:collage-bijection} again, we conclude
  that this pushout is given by the simplicial category $\coll(\lan_IW)$, as
  claimed.
\end{proof}
\fi 
In ordinary unenriched category theory, the colimit cone under a $\eD$-shaped diagram may be formed as the left Kan extension along the inclusion $\eD\inc\eD\join\catone$ into the category $\eD\join\catone$ formed by freely adjoining a terminal object ``$\top$'' to $\eD$. The following lemma reveals that the collage plays the roll of the category $\eD\join\catone$ for weighted colimits\ifnewcut\else, a perspective which we will return to in \S\ref{sec:oplax}\fi.

\begin{lem}\label{lem:lan-along-collage} 
  The pointwise left Kan extension of any simplicial functor
  $F\colon\eD\to\eK$ along $I \colon \eD\inc\coll(W)$ exists if and only if
  the colimit $\colim^W\!F$ exists in $\eK$, in which case
  $\lan_IF(\top)\cong\colim^W\!F$. \ifnewadd\qed\fi 
\end{lem}
\ifnewcut\else
\begin{proof}
  Since $\eD\inc\coll(W)$ is fully faithful, when the pointwise left Kan
  extension of any simplicial functor $F\colon\eD\to\eK$ along
  $\eD\in\coll(W)$ exists, it is displayed by an isomorphism:
  \begin{equation*}
    \xymatrix@=1.5em{
      & {\eD}\ar@{^(->}[dl]_I \ar@{}[d]|(.6){\cong}\ar[dr]^{F} & \\
      {\coll(W)}\ar[rr]_-{\lan_IF} && {\eK}
    }
  \end{equation*} 
  By Corollary \ref{cor:collage-bijection}, this data defines a $W$-cocone
  under $F\cong \lan_IF\circ I$ with nadir $\lan_IF(\top) \in \eK$. It is
  easy to verify that the universal property of the left Kan extension
  specializes to describe the universal property of the colimit cocone for
  $\colim^W\!F$, and conversely that the universal property of the weighted
  limit cocone implies the universal property of the left Kan extension.
\end{proof}

\subsection{Flexible weights as projective cell complexes}\label{sec:flexible-weight}
\fi

In order to understand the sense in which certain weighted colimits, including in particular the oplax colimits to be introduced below, are homotopically well
behaved, we recall some facts about weights and simplicial computads from
\S\refII{subsec:collage}:

\begin{defn}[flexible weights and projective cell complexes]\label{defn:proj-cell-cx}
  For a simplicial category $\eD$, the \emph{projective $n$-cell} associated with $[n] \in \Del$ and $d \in \eD$ is the  simplicial natural
  transformation  \[\boundary\Del^n\times \Fun_{\eD}(-,d)\inc
    \Del^n\times \Fun_{\eD}(-,d).\] 
  A natural transformation $\alpha\colon W\to V$ in $\SSet^{\eD\op}$ is a
  \emph{relative projective cell complex} if it factors as a
  countable composite of pushouts of coproducts of projective cells. A weight
  $W$ in $\SSet^{\eD\op}$ is a \emph{flexible weight} if the map
  $!\colon\emptyset\to W$ is a relative projective cell complex\ifnewcut\else, i.e., if $W$ is a projective cell complex\fi.
\end{defn}

\ifnewcut\else
Our interest in colimits weighted by flexible weights is due to the fact that
they are homotopically well behaved.\fi The following result extends without change to pointwise
cofibrant diagrams valued in any model category enriched over the Joyal model
structure on simplicial sets. 

\begin{prop}[{\refII{prop:proj-wlim-homotopical}, \refVII{prop:flexible-weights-are-htpical}}]\label{prop:flexible-weights-are-htpical}$\quad$
  \begin{enumerate}[label=(\roman*)]
  \item For a flexible weight $W \colon \eD\op\to\SSet$ and any diagram $F \colon \eD\to\SSet$, $\colim^W\!F$ may be expressed as a countable composite of
    pushouts of coproducts of maps \[\boundary\Del^n\times Fd\inc\Del^n\times
      Fd.\]
  \item If $\alpha \colon F \to G$ is a simplicial natural transformation
    between two such diagrams whose components are weak equivalences in the
    Joyal model structure,
    then for any flexible weight $W$ the map \[\colim\nolimits^W\!\alpha\colon\colim\nolimits^W\!F\to\colim\nolimits^W\!G\] is a weak
    equivalence in the Joyal model structure. \ifnewadd\qed\fi
  \end{enumerate}
\end{prop}



The collage construction defines a correspondence between flexible weights and \emph{simplicial computads}, a class of ``freely generated'' simplicial categories that define precisely the cofibrant objects \cite[\S 16.2]{Riehl:2014kx} in  the model structure due to  Bergner \cite{Bergner:2007fk}. 

\begin{defn}[simplicial computad]\label{defn:simplicial-computad}
A simplicial category $\eA$, regarded as a simplicial object $[n] \mapsto \eA_n$ in the category of categories with a common set of objects   and identity-on-objects functors, 
 is a \emph{simplicial computad} if and only if:
  \begin{itemize}
    \item each category $\eA_n$ of $n$-\emph{arrows} is freely generated by the reflexive directed graph of \emph{atomic} $n$-arrows, these being those arrows that admit no non-trivial factorizations, and 
       \item if $f$ is an atomic $n$-arrow in $\eA_n$ and $\alpha\colon [m]\to[n]$ is a degeneracy operator in $\Del$ then the degenerated $m$-arrow $f\cdot\alpha$ is atomic in $\eA_m$.
  \end{itemize}
  \ifnewcut\else
A simplicial category  $\eA$ is a simplicial computad if and only if all of its non-identity arrows $f$ can be expressed uniquely as a composite 
\begin{equation}\label{eq:computad-arrow-decomp}
  f = (f_1 \cdot \alpha_1) \circ (f_2 \cdot \alpha_2) \circ \cdots \circ (f_\ell \cdot \alpha_\ell)
\end{equation}
in which each $f_i$ is non-degenerate and atomic and each $\alpha_i\in\Del$ is a degeneracy operator.\fi
\end{defn}

We have the following recognition principle for flexible weights on simplicial
computads,  a mild variant of Proposition~\refII{prop:projcofchar}, proven in \S\refVII{ssec:collage}.

\begin{prop}[{relating flexible weights and simplicial computads, \refVII{thm:flexible-collage}}]
  \label{prop:flexible-weight-as-computad}
  Suppose that $\eD$ is a simplicial computad. Then a weight $W\colon\eD\op\to\SSet$ is flexible if and only if its collage
  $\coll(W)$ is a simplicial computad. \ifnewadd\qed\fi
\end{prop}

\ifnewcut\else
\subsection{Oplax colimits}\label{sec:oplax}
\fi 
\ifnewadd
By the next result, the left adjoint to the homotopy coherent nerve
\[
\adjdisplay \gC -| \hN : \sCat -> \sSet.\]
the \emph{homotopy coherent realization} functor, provides a source of flexible weights. See \S\refVI{sec:coherent-nerve} for a more leisurely presentation with considerably more details.

\begin{prop}[\refVI{prop:gothic-C}]\label{prop:gothic-C}
For any simplicial set $X$, the \emph{homotopy coherent realization} $\gC{X}$ is a simplicial computad. \qed 
\end{prop} 
\fi

\ifnewcut\else Oplax colimits represent particular cones under a homotopy coherent diagram
$\gC{X} \to \eK$ indexed by a simplicial set $X$.  In Definition \ref{defn:oplax-colimit-weight}, we
first present the collage construction that describes the shape of a lax cocone
and then use Proposition \ref{prop:collage-adjunction} to extract the
corresponding weight. To give a concise description of the collage that defines the oplax weight, we make use of the simplicial computad structure on the \emph{homotopy coherent realization} $\gC{X}$ of a simplicial set $X$, our term for the left adjoint to the homotopy coherent nerve
\[
\adjdisplay \gC -| \hN : \sCat -> \sSet.\]
We briefly review this material from  \S\refVI{sec:coherent-nerve}, which gives a more leisurely presentation with considerably more details.

  \begin{ex}[{homotopy coherent simplices as simplicial computads; \S\refVI{sec:htpy-coh-simplex}}]\label{ex:homotopy-coherent-simplex} Recall the simplicial category $\mathfrak{C}\Del^n$ whose objects are integers $0,1,\ldots, n$ and whose function complexes are the cubes \[ \Fun_{\mathfrak{C}\Del^n}(i,j) = \begin{cases} \Cube^{j-i-1} & i < j \\ \Del^0 & i = j \\ \emptyset & i > j \end{cases}\] Here we write $\Cube^k \defeq (\Del^1)^k$. For $i < j$, the vertices of $\Fun_{\mathfrak{C}\Del^n}(i,j)$ are naturally identified with subsets of the closed interval $[i,j] = \{ i \leq t \leq j\}$ containing both endpoints, a set whose cardinality is $j-i-1$; more precisely, $\Fun_{\mathfrak{C}\Del^n}(i,j)$ is the nerve of the poset with these elements, ordered by inclusion. Under this isomorphism, the composition operation  corresponds to the simplicial map 
  \[
  \xymatrix@R=1em{ \Fun_{\gC\Del^n}(i,j)\times\Fun_{\gC\Del^n}(j,k)\ar[r]^-\circ \ar@{}[d]|{\rotatebox{90}{$\cong$}}& \Fun_{\gC\Del^n}(i,j) \ar@{}[d]|{\rotatebox{90}{$\cong$}} \\
\Cube^{\times
    (j-i-1)}\times\Cube^{\times (k-j-1)}\ar[r] & \Cube^{\times (k-i-1)}}\]
     which
  maps the pair of vertices $T \subset [i,j]$ and $S \subset [j,k]$ to $T \cup S \subset [i,k]$.
  
  Again for $i < j$, an $r$-arrow $T^\bullet$ in $\Fun_{\gC\Del^n}(i,j)$ corresponds to a sequence
  \[ T^0 \subset T^1 \subset \cdots \subset T^r\] of subsets of $[i,j] = \{ i \leq t \leq j\}$ and is non-degenerate if and only if each of these inclusions are proper. The composite of a pair of $r$-arrows $T^\bullet \colon i \to j$ and $S^\bullet \colon j \to k$ is the levelwise union $T^\bullet \cup S^\bullet \colon i \to k$ of these sequences.
  
 From this description, it is easy to see that the simplicial category $\gC\Del^n$  is a simplicial computad (Lemma \refVI{lem:simplex-computad}), in which an $r$-arrow $T^\bullet$ from $i$ to $j$ is atomic if and only if the set $T^0 = \{i,j\}$; the only atomic $r$-arrows from $j$ to $j$ are identities. Geometrically,   the atomic arrows in each function complex $\Fun_{\gC\Del^n}(i,j) \cong \Cube^{j-i-1}$ are precisely those simplices that contain the initial vertex in the poset whose nerve defines the simplicial cube.
\end{ex}

If $X$ is a simplicial subset of $\Del^n$, then Lemma \refVI{lem:subcomputad-inclusion} tells us that 
 its homotopy coherent realisation
  $\gC{X}$ is a simplicial subcomputad of $\gC\Del^n$.

\begin{ex}[{homotopy coherent nerves of subsimplices \refVI{ex:subcomputad-of-simplex}}]\label{ex:subcomputad-of-simplex}
In particular: 
 \begin{enumerate}[label=(\roman*)]
  \item\label{itm:boundary} \textbf{boundaries:} The inclusion $\gC\boundary\Del^n\inc\gC\Del^n$ is the identity on objects and full on all function complexes except for the one from $0$ to $n$. The inclusion
         \[
  \xymatrix@R=1em{ \Fun_{\gC\boundary\Del^n}(0,n)\ar@{^(->}[r] \ar@{}[d]|{\rotatebox{90}{$\cong$}}& \Fun_{\gC\Del^n}(0,n)\ar@{}[d]|{\rotatebox{90}{$\cong$}} \\
\boundary\Cube^{n-1}\ar@{^(->}[r] & \Cube^{n-1}}\]
    is isomorphic to the cubical boundary inclusion, where $\boundary\Cube^k$ is the domain of the iterated Leibniz product
  $(\boundary\Del^1\subset\Del^1)^{\leib\times k}$.\footnote{For more details about the Leibniz or ``pushout-product'' construction see \cite[\S 4]{RiehlVerity:2013kx}.}
     
  \item\label{itm:inner-horn} \textbf{inner horns:} The inclusion $\gC{\Horn^{n,k}}\inc\gC\Del^n$ is identity on objects and full on all function complexes except for the one from $0$ to $n$. The inclusion
     \[
  \xymatrix@R=1em{ \Fun_{\gC{\Horn^{n,k}}}(0,n)\ar@{^(->}[r] \ar@{}[d]|{\rotatebox{90}{$\cong$}}& \Fun_{\gC\Del^n}(0,n)\ar@{}[d]|{\rotatebox{90}{$\cong$}} \\
\CHorn^{n-1,k}_1\ar@{^(->}[r] & \Cube^{n-1}}\]
 is isomorphic to the cubical horn inclusion, defined by the the following Leibniz product:
  \begin{equation*}
    (\boundary\Del^1\subset\Del^1)^{\leib\times(j-1)}\leib\times
    (\Del^{\fbv{1}}\subset\Del^1)\leib\times
    (\boundary\Del^1\subset\Del^1)^{\leib\times(k-j)}
  \end{equation*}
  \end{enumerate}
\end{ex}

\begin{defn}[bead shapes]
  We shall call those atomic arrows $T^\bullet\colon 0 \to n$ of $\gC{\Del^{n}}$
  which are not members of $\gC\boundary\Del^n$ \emph{bead shapes}. By Examples \ref{ex:homotopy-coherent-simplex} and \ref{ex:subcomputad-of-simplex}, an $r$-dimensional bead shape $T^\bullet\colon
  0\to n$ is given by a sequence of subsets
  \[ \{0,n\} = T^0 \subset T^1 \subset \cdots \subset T^r = [0,n]\]
  with $T^0 = \{0,n\}$  and $T^r$ equal to the full interval $[0,n] = \{ 0 \leq t \leq n\}$. 
  \end{defn}

More generally, any simplicial category $\gC{X}$ arising as the homotopy coherent realization of a simplicial set $X$ defines a simplicial computad whose atomic arrows $(x,T^\bullet)$, described in Proposition \ref{prop:gothic-C}, are called \emph{beads in $X$}.  As a consequence of this result we find that 
  $r$-simplices of $\gC{X}$ correspond to sequences of abutting beads,
  structures which are called \emph{necklaces} in the work of Dugger and
  Spivak~\cite{DuggerSpivak:2011ms} and Riehl~\cite{Riehl:2011ot}. In this terminology, 
   $\gC X$ is a simplicial computad in which the atomic arrows are those necklaces that consist of a single bead with non-degenerate image.

\begin{prop}[{$\gC X$ as a simplicial computad; \refVI{prop:gothic-C}}]\label{prop:gothic-C}
The homotopy coherent realization $\gC{X}$ of a simplicial set $X$ is a simplicial computad with
  \begin{itemize}
  \item   objects the vertices of $X$ and 
  \item non-degenerate atomic $r$-arrows given by pairs $(x,T^\bullet)$, wherein
  $x$ is a non-degenerate $n$-simplex of $X$ for some $n > r $ and
  $T^\bullet\colon 0\to n$ is an $r$-dimensional bead shape. 
  \end{itemize}
  The domain of $(x,T^\bullet)$ is the initial vertex $x_0$ of $x$ while the codomain is the 
terminal vertex $x_n$.
 \end{prop}
 
 The point of this review is to permit us to define weights for oplax colimits of homotopy coherent diagrams valued in a simplicially (or frequently quasi-categorically) enriched category. In a homotopy coherent diagram, the indexing shape is given by the homotopy coherent realization of a simplicial set $X$. In this context, the join operation $X\join\Del^0$ produces another simplicial set with a freely adjoined cocone vertex. We shall argue that the its homotopy coherent realization defines a collage that presents the weight for oplax colimits. \fi
 
 \begin{rec} For any simplicial set $X$, there is a canonical inclusion $X\inc X\join\Del^0$ into its join
  with the point.  The join $X\join\Del^0$ has a single vertex of $X\join\Del^0$ that is not also a vertex of its subset $X$, which we shall denote by ``$\top$.'' Now for each non-degenerate $n$-simplex $x\in X$ the join $X\join\Del^0$ has
  two corresponding non-degenerate simplices:
  \begin{itemize}
  \item a simplex of dimension $n$ identified with $x$
itself and 
  \item a simplex $(x,\top)$ of dimension $n+1$,
  \end{itemize} and these two cases enumerate all of the
  non-degenerate simplices of $X\join\Del^0$ with the exception of $\top$.
\end{rec}

\ifnewadd Oplax colimits represent particular cones under a homotopy coherent diagram
$\gC{X} \to \eK$ indexed by a simplicial set $X$. In this context, the homotopy coherent realization of of the joint $X\join\Del^0$ defines a collage that presents the weight for oplax colimits.\fi

\begin{defn}[weights for oplax colimits]\label{defn:oplax-colimit-weight}
Applying homotopy coherent
  realisation to the canonical inclusion $X\inc X\join\Del^0$ for any simplicial set $X$, yields a simplicial
  subcomputad $I_X\colon \gC{X} \inc \gC[X\join\Del^0]$\ifnewcut\else. Now from Proposition~\refVI{prop:gothic-C} we know that $\gC[X\join\Del^0]$ may be
  built from $\gC{X}$ by adjoining atomic arrows corresponding to beads
  $((x,\top),T^\bullet)$ and these all have codomain $\top$. 
 It is clear\fi\ifnewadd ~so \fi that the conditions discussed in
Proposition \ref{prop:collage-adjunction}\ref{itm:collage-fun} hold for the inclusion $\langle \top, I_X
  \rangle\colon \catone + \gC{X} \inc \gC[X\join\Del^0]$. Hence, via the counit isomorphism of the collage adjunction, this simplicial category is isomorphic to the collage of the corresponding weight, defining the \emph{weight for oplax colimits} of diagrams of
  shape $\gC{X}$.
  \begin{equation*}
    \xymatrix@R=0em@C=5em{
      \gC{X}\op\ar[r]^{L_X} & {\SSet}
    }\qquad\text{given~by}\qquad L_X(x) \defeq \Fun_{\gC[X\join\Del^0]}(x,\top).
  \end{equation*}
   When $F\colon\gC{X}\to\eK$ is a homotopy coherent diagram of
  shape $X$, then its \emph{oplax colimit} is defined to be the weighted colimit
  \[ \colim\nolimits^{\mathrm{oplax}}F \defeq \colim\nolimits^{L_X}F.\] 
\end{defn}

\begin{rmk}\label{rmk:oplax-vs-pseudo}
The oplax weights being defined here are precisely the ``pseudo'' weights introduced in Definition \refVII{defn:weight-for-pseudo-limits}. The reason for the difference in nomenclature is that in that paper the diagrams considered in \cite{RiehlVerity:2018rq} 
 are valued in Kan complex enriched categories, whereas here the diagrams are valued in quasi-categorically (or simplicially) enriched categories. In a Kan complex, the 1-simplex $\Del^1$ represents an invertible morphism, while in a quasi-category it models a non-invertible morphism.
\end{rmk}

Immediately from Proposition
  \ref{prop:flexible-weight-as-computad}:

\begin{lem}[{\refVII{lem:pseudo-is-flexible}}]\label{lem:oplax-is-flexible}
  For all simplicial sets $X$ the weight $L_X\colon\gC{X}\op\to\SSet$ for oplax
  colimits of diagrams of shape $\gC{X}$ is a flexible weight. \qed
\end{lem}

\subsection{Cocartesian fibrations and quasi-categorical collages}\label{ssec:cocartesian}

In this section, we construct an explicit example of an oplax colimit of diagram of quasi-categories via the \emph{quasi-categorical collage construction}\ifnewcut\else, which we introduce in \S\ref{sec:qcat-collage}\fi. In an important special case, the quasi-categorical collage defines a cocartesian fibration over the 1-simplex, so \ifnewcut\else we begin  in \S\ref{sec:cocart-basics} with a review of this notion, since we shall need it anyways. In \S\ref{sec:cat-of-cocart},\fi we \ifnewadd first \fi introduce the quasi-categorically enriched category of cocartesian fibrations and cartesian functors\ifnewcut\else and observe that pullback defines a functor between such categories\fi. 

\ifnewcut\else
\subsection{Cocartesian fibrations of quasi-categories}\label{sec:cocart-basics}
\fi 

Of the many equivalent definitions of cocartesian fibration (see \S\refIV{sec:cartesian} and \S\refVI{sec:cocartesian}), the following will be the most convenient for this paper:

\begin{defn}[{\cite[2.4.1.8,2.4.2.1]{Lurie:2009fk}, \refIV{cor:lurie-cartesian}}]\label{defn:qcat-cocart} Let $p\colon \qE \tfib \qB$ be an isofibration between quasi-categories.
\begin{enumerate}[label=(\roman*)]
\item\label{itm:qcat-cocart-arrow} A 1-arrow $\chi\colon e\to e'$ of $\qE$ is $p$-\emph{cocartesian} if and only if any
  lifting problem
  \begin{equation*}
    \xymatrix@C=3em{
      \Delta^{\fbv{0,1}}\ar@/^2ex/[rr]!L+/u 4pt/^\chi \ar[r] &
      \Horn^{n,0} \ar[d] \ar[r] &      \qE \ar@{->>}[d]^{p} \\ &
      \Delta^n \ar[r] \ar@{-->}[ur] & \qB}
  \end{equation*}
  has a solution. 
  \item\label{itm:qcat-cocart-fib} An isofibration
  $p\colon\qE\tfib\qB$ is a \emph{cocartesian} fibration of quasi-categories precisely when any arrow
  $\alpha\colon pe \to b$ in $\qB$ admits a lift to an arrow $\chi\colon e\to
  e'$ in $\qE$ which enjoys the lifting property of \ref{itm:qcat-cocart-arrow}.
  \end{enumerate}
\end{defn}

\ifnewcut\else
For efficiency of exposition, we focus largely on the cocartesian fibrations, and leave it to the reader to formulate the dual statements for cartesian fibrations, obtained by replacing each simplicial set by its opposite.

\begin{ex} The product projection $\pi \colon \qA \times \qB \tfib \qB$ defines a \emph{bifibration}, that is, both a cocartesian and a cartesian fibration. A 1-arrow of $\qA \times \qB$ is $\pi$-cocartesian (and also $\pi$-cartesian) just when its component in $\qA$ is an isomorphism.
\end{ex}

\begin{ex}[{\refIV{ex:domain-fibration}, \refIV{ex:comma-fibrations}}]\label{ex:comma-qcat}  For any quasi-category $\qB$, we write $\qB^\cattwo \defeq \qB^{\Del^1}$ for its cotensor with the 1-simplex and $p_0,p_1 \colon \qB^\cattwo \tfib \qB$ for the evaluation maps at the vertices $0,1 \in \Del^1$ respectively. The codomain functor $p_1 \colon \qB^\cattwo \tfib \qB$ is a cocartesian fibration, in which the $p_1$-cocartesian arrows are those whose projections along $p_0$ are invertible. Dually, the domain functor $p_0 \colon \qB^\cattwo \tfib \qB$ is a cartesian fibration, in which the $p_0$-cartesian arrows are those whose codomain components are invertible.

More generally, for any cospan of quasi-categories $f \colon \qB \to \qA$ and $g \colon \qC \to \qA$, the \emph{comma quasi-category} is defined by the pullback
\[ 
\xymatrix{ f\comma g \ar@{->>}[d]_{(p_1,p_0)} \ar[r] \pbexcursion&  \qA^\cattwo \ar@{->>}[d]^{(p_1,p_0)} \\ \qC \times \qB\ar[r]_{g \times f} & \qA \times \qA}\]  and once more $p_1 \colon f \comma g \tfib \qC$ is a cocartesian fibration and $p_0 \colon f \comma g \tfib \qB$ is a cartesian fibration.

In the special case where one of the functors in the cospan is taken to be the identity, we write $f \comma \qA$ and $\qA \comma f$ for what we call the \emph{contravariant} and \emph{covariant} representable comma quasi-categories respectively. In the special case where the functor $a \colon 1 \to \qA$ identifies a vertex of $\qA$, the codomain projection $p_1 \colon a \comma \qA \tfib \qA$ is a cocartesian fibration that encodes the covariant representable functor associated to $a$, while the domain projection $p_0 \colon a \comma \qA \tfib \qA$ is a cartesian fibration encoding the contravariant representable functor. These define the images of $a$ in the co- and contravariant Yoneda embeddings of \S\refVI{sec:yoneda}.
\end{ex}

\begin{lem}[{\refVI{lem:cocart-cyl-are-cocart-arrows}}]\label{lem:ptwise-cocart-cyl} If $p \colon \qE \tfib \qB$ is a cocartesian fibration and $X$ is a simplicial set, then $p^X \colon \qE^X \tfib \qB^X$ is a cocartesian fibration in which a 1-arrow $e \colon \Del^1 \to \qE^X$ is $p^X$-cocartesian just when for each vertex $x \in X$ its component $e(x \cdot \degen^0,\id_{[1]}) \colon \Del^1 \to \qE$ is $p$-cocartesian.
\end{lem}

In \S\refVI{sec:cocartesian}, a $p^X$-cocartesian 1-arrow is called a \emph{pointwise $p$-cocartesian cylinder}.

  \begin{lem}[{\refVI{lem:cocart-cylinder-extensions}}]\label{lem:cocart-cylinder-extensions}
Let $X \inc Y$ be a simplicial subset of a simplicial set $Y$.
\begin{enumerate}[label=(\roman*)]
  \item Any lifting problem
    \begin{equation*}
      \xymatrix@=2em{
        {X\times\Del^1\cup Y\times\Del^{\{0\}}}\ar[r]^-{e}\ar@{^(->}[d] &
        {\qE}\ar@{->>}[d]^{p} \\
        {Y\times\Del^1}\ar[r]_{b}\ar@{.>}[ru]^{\bar{e}} & {\qB}
      }
    \end{equation*}
    in which the cylinder ${X\times\Del^1}\subseteq{X\times\Del^1\cup
      Y\times\Del^{\{0\}}}\stackrel{e}\longrightarrow {E}$ is pointwise
    $p$-cocartesian admits a solution $\bar{e}$ which is also pointwise
    $p$-cocartesian.
  \item Any lifting problem ($n>1$)
    \begin{equation*}
      \xymatrix@=2em{
        {X\times\Del^n\cup Y\times\Horn^{n,n}}\ar[r]^-{e}\ar@{^(->}[d] &
        {\qE}\ar@{->>}[d]^{p} \\
        {Y\times\Del^n}\ar[r]_{b}\ar@{.>}[ru]^{\bar{e}} & {\qB}
      }
    \end{equation*}
    in which the cylinder ${Y\times\Del^{\{n-1,n\}}}\subseteq {X\times\Del^n\cup
      Y\times\Horn^{n,n}}\stackrel{e}\longrightarrow {E}$ is pointwise
    $p$-cocartesian admits a solution $\bar{e}$.
\end{enumerate}
\end{lem}

\begin{defn}\fi If $p$ and $q$ are cocartesian fibrations over $\qB$ then a functor
\[ \xymatrix{ \qE \ar@{->>}[dr]_p \ar[rr]^g & & \qF \ar@{->>}[dl]^{q} \\ & \qB}\] is a \emph{cartesian functor} just when it carries $p$-cocartesian 1-arrows to $q$-cocartesian 1-arrows. \ifnewcut\else
\end{defn}
\fi
As one illustration of the importance of this notion:

\begin{prop}[{\refVIII{prop:equivalence-of-fibrations}}]\label{prop:equivalence-of-fibrations} A cartesian functor between cocartesian fibrations of quasi-categories is an equivalence if and only if it is a fiberwise equivalence:
\[ \xymatrix{ \qE\ar[rr]^g \ar@{->>}[dr]_p & & \qF \ar@{->>}[dl]^q \\ & \qB}\]  i.e., for each $b \in \ob\qB$, the induced functor $g_b \colon \qE_b \to \qF_b$ is an equivalence. \ifnewadd\qed\fi
\end{prop}

\ifnewcut\else
\subsection{The quasi-categorically enriched category of cocartesian fibrations}\label{sec:cat-of-cocart}
\fi 

If $\qB$ is a quasi-category, then we adopt the notation $\qCat_{/\qB}$ for the quasi-categoricaly enriched category of isofibrations over $\qB$ defined as follows.

\begin{defn}\label{defn:fun-over-B} For a quasi-category $\qB$, let $\qCat_{/\qB}$  denote the category whose:
\begin{itemize}
\item objects are isofibrations $p \colon \qE \tfib \qB$ with codomain $\qB$ and
\item whose function complexes $\Fun_{\qB}(p \colon \qE \tfib \qB, q \colon \qF \tfib \qB)$ are defined by the pullbacks
\[ \xymatrix{ \Fun_{\qB}(p \colon \qE \tfib \qB, q \colon \qF \tfib \qB) \ar@{->>}[d] \ar[r] \pbexcursion & \Fun(\qE,\qF) \ar[d]^{q \circ -} \\ \Del^0 \ar[r]^-p & \Fun(\qE,\qB)}\]
\end{itemize}
where $\Fun(\qE,\qF) \cong \qF^\qE$ denotes the usual internal hom in $\qCat$.
\end{defn}

\begin{defn}\label{defn:cart-fun-over-B} For a quasi-category $\qB$, let $\coCart(\qCat)_{/\qB}$  denote the category whose:
\begin{itemize}
\item objects are cocartesian fibrations $p \colon \qE \tfib \qB$ with codomain $\qB$ and
\item whose function complexes $\Fun^c_{\qB}(p \colon \qE \tfib \qB, q \colon \qF \tfib \qB)$ are defined to be the full sub quasi-categories of the function complexes $\Fun_{\qB}(p \colon \qE \tfib \qB, q \colon \qF \tfib \qB)$ of $\qCat_{/\qB}$ defined by restricting the 0-arrows to be cartesian functors over $\qB$.
\end{itemize}
The quasi-categorically enriched category $\Cart(\qCat)_{/\qB}$ of cartesian fibrations and cartesian functors is defined similarly.
\end{defn}

Proposition \refIV{prop:cart-fib-pullback} proves that the pullback of a cocartesian fibration is a cocartesian fibration
\[
\xymatrix{ \qF \ar@{->>}[d]_q \ar[r]^g \pbexcursion & \qE\ar@{->>}[d]^p \\ \qA \ar[r]_f & \qB}\] in which an arrow $\chi$ is $q$-cocartesian if and only $g\chi$ is $p$-cocartesian. It follows that pullback also preserves cartesian functors. Hence:

\begin{prop}\label{prop:pullback-stability}
Pullback along any $f \colon \qA  \to \qB$ defines a quasi-categorically enriched functor
     \[
  \xymatrix@R=1em{ \coCart(\qCat)_{/\qB} \ar[r]^-{f^*} \ar@{}[d]|{\rotatebox{90}{$\supset$}}& \coCart(\qCat)_{/\qA}\ar@{}[d]|{\rotatebox{90}{$\supset$}} \\
\qCat_{/\qB} \ar[r]^-{f^*}& \qCat_{/\qA}} \qed\] 
\end{prop}

We now argue that the pullback functor preserves simplicial tensors. This will be used in \S\ref{sec:powerful} to show that its right adjoint is simplicially enriched, when this functor exists.

\begin{obs}[tensors and pullback]\label{obs:tensor-pullback} Let $X \in \sSet$ be a simplicial set. The tensor of an isofibration $p \colon \qE \tfib \qB$ with $X$ is the right-hand vertical composite, which pulls back to the right-hand vertical composite
\[ \xymatrix{ \qF \times X \ar@{->>}[d]_{\pi_1} \ar[r] \pbexcursion & \qE \times X \ar@{->>}[d]^{\pi_1}\\  \qF \pbexcursion \ar@{->>}[d]_{f^*(p)} \ar[r] & \qE \ar@{->>}[d]^p \\ \qA \ar[r]_f & \qB}\] which defines the tensor of $f^*(p) \colon \qF \tfib \qA$ with $X$.
\end{obs}

The following lemma tells us that this tensor construction respects cartesian functors.

\begin{lem}\label{lem:cart-functor-tensor} For any simplicial set $X$ and cocartesian fibrations $p \colon \qE \tfib \qB$ and $q \colon \qF \tfib \qB$, the isomorphism $\Fun_{\qB}(\qE \times X, \qF) \cong \Fun_{\qB}(\qE,\qF)^X$  restricts to an isomorphism \[\Fun^c_{\qB}(\qE \times X, \qF) \cong \Fun^c_{\qB}(\qE,\qF)^X.\]
\end{lem} 
\begin{proof}
We make use of Theorem \refIV{thm:cart.fun.chars} which provides the following characterization of the sub quasi-category $\Fun^c_{\qB}(\qE,\qF) \subset \Fun_{\qB}(\qE,\qF)$. Any functor $f \colon \qE \to \qF$ over $\qB$ induces a commutative square over $\qB$
\[ \xymatrix{ \qE \ar[r]^f \ar[d] & \qF \ar[d] \\ p \comma \qB \ar[r]_{(f,\id_\qB)} \ar@/^2ex/@{-->}[u]^\ell_\dashv & q \comma \qB \ar@/_2ex/@{-->}[u]_\ell^\vdash}\] whose vertical functors are the canonical ones induced by $p \colon \qE \tfib \qB$ and $q \colon \qF \tfib \qB$. Because $p$ and $q$ are cocartesian, Theorem \refIV{thm:cart.fib.chars} proves the vertical functors admit left adjoints over $\qB$. Theorem \refIV{thm:cart.fun.chars} proves that $f$ is cartesian if and only if the mate of this canonical isomorphism is an isomorphism.

The mate that detects whether $f$ is a cartesian functor lives as a 1-simplex in the simplicial set
\[  \qop{Sq}_{\qB}(p\comma \qB \to \qE,  q \comma \qB \to \qF) \defeq \Fun_{\qB}(\qE , \qF) \times_{\Fun_{\qB}(p \comma \qB, \qF)} \Fun_{\qB}(p \comma \qB, q\comma \qB).\] of commutative squares from $\ell \colon p\comma \qB \to \qE$ to $\ell \colon q \comma \qB \to \qF$. The adjunction over $\qB$ associated to the cocartesian fibration $\qE \times X \xtfib{\pi} \qE\xtfib{p} \qB$ is
\[ \xymatrix{ \qE \times X \ar[rr]^-\perp \ar[dr]_-{p\pi}  & & p\comma\qB \times X \ar[dl]^{p_0\pi} \ar@/_2ex/[ll]_-\ell \\ & \qB }\] the product of the adjunction for $p$ with $X$. In particular, 
\[ \qop{Sq}_{\qB}(p\comma \qB \times X \to \qE \times X, q \comma \qB \to \qF)  \cong \qop{Sq}_{\qB}(p\comma \qB \to \qE, q \comma \qB \to \qF)^X.\] Now a 1-simplex $\qC^X$ is an isomorphism if and only if it is a pointwise isomorphism, which proves that $\Fun^c_{\qB}(\qE \times X, \qF) \cong \Fun^c_{\qB}(\qE,\qF)^X$.
\end{proof}

\ifnewcut\else
\subsection{Collages for quasi-categories}\label{sec:qcat-collage}
\fi

We conclude this section with an example of an oplax colimit. When $X=\Del^1$ a homotopy coherent diagram $\gC{\Del^1} \to \qCat$ is just a functor $f \colon \qA \to \qB$ between quasi-categories. The oplax colimit in simplicial sets is given by the pushout 
\[ 
\xymatrix{ \qA \ar[r]^f \ar[d]_{\id\times\face^0} & \qB \ar[d] \\ \qA \times \Del^1 \ar[r] & \colim\nolimits^{\mathrm{oplax}}f \poexcursion }
\]
Up to equivalence, this oplax colimit is modeled by the \emph{quasi-categorical collage construction} that we now introduce.

\begin{defn}[the quasi-categorical collage construction X.F.5.2]\label{defn:qcat-collage} Consider any cospan $f \colon \qA \to \qC$ and $g \colon \qB \to \qC$, with $\qA$, $\qB$, and $\qC$ all quasi-categories. Define a new simplicial set $\coll(f,g)$ 
by declaring that
\[ \coll(f,g)_n = \Bigl\{ \left(\Del^i \xrightarrow{a} \qA, \Del^j \xrightarrow{b} \qB, \Del^n \xrightarrow{c} \qC\right) \bigg\vert {\begin{array}{ll} c \vert_{\fbv{0,\ldots, i}} &= f(a), \\ c\vert_{\fbv{n-j,\ldots, n}} &= g(b),\end{array}} \begin{array}{c} i,j \geq -1, \\ i+j = n-1.\end{array}\Bigr\}\]
with the convention that conditions indexed by $\Del^{-1}$ are empty (or that each simplicial set is terminally augmented). There are simplicial maps
\[
\begin{tikzcd} B \arrow[r, hook] \arrow[d] \arrow[dr, phantom, "\lrcorner" very near start] & \coll(f,g) \arrow[d, "\rho"] & A \arrow[l, hook'] \arrow[d] \arrow[dl, phantom, "\llcorner" very near start] \\ \{1\} \arrow[r, hook] & \Del^1 & \{0\} \arrow[l, hook']
\end{tikzcd}
\]
where the map $\rho$ sends an $n$-simplex $(a \colon \Del^i \to \qA, b \colon \Del^j \to \qB, c \colon \Del^n \to \qC)$ to the $n$-simplex $[n] \to [1]$ that carries $0,\ldots ,i$ to $0$ and $i+1,\ldots, n$ to $1$. Note that the fiber of $\rho$ over $0$ is isomorphic to $\qA$ while the fiber of $\rho$ over $1$ is isomorphic to $\qB$.
\end{defn}

\ifnewcut\else
\begin{rmk}[on right and left]
As with simplicial collages, we customarily write $\qB+\qA \inc \coll(f,g)$ for the inclusions of the fibers over $1$ and $0$ --- with the fiber over 1 on the left and the fiber over 0 on the right. As with our convention for quasi-categories in Example \ref{ex:comma-qcat}, this positions the covariantly-acting quasi-category on the ``left'' and the contravariantly-acting quasi-category on the ``right.''
\end{rmk}
\fi

\begin{lem}[X.F.5.3]\label{lem:collage-qcat} The map $\rho \colon \coll(f,g) \to \Del^1$ is an inner fibration. In particular, the simplicial set $\coll(f,g)$ is a quasi-category. \ifnewadd\qed\fi
\end{lem}\ifnewcut\else
\begin{proof} Since the fibers of $\rho$ over $0$ and $1$ are the quasi-categories $\qA$ and $\qB$, it suffices to consider inner horns
\[
\xymatrix{ \Horn^{n,k} \ar[r] \ar[d] & \coll(f,g) \ar[d]^\rho \\ \Del^n \ar[r]_{\alpha} \ar@{-->}[ur] & \Del^1}\] for which $\alpha \colon [n] \to [1]$ is a surjection. Suppose $\alpha$ carries $0,\ldots, i$ to $0$ and $i+1,\ldots, n$ to $1$. Note that for any $0 < k < n$, the faces $\fbv{0,\ldots, i}$ and $\fbv{i+1,\ldots,n}$ of $\Del^n$ belong to the horn $\Horn^{n,k}$. In particular, the map $\Horn^{n,k} \to \coll(f,g)$ identifies simplices $a \colon \Del^i \to \qA$ and $\Del^{n-i-1} \to \qB$ together with a horn $\Horn^{n,k} \to \qC$ whose initial and final faces are the images of these simplices under $f \colon \qA \to \qC$ and $g \colon \qB \to \qC$. Since $\qC$ is a quasi-category this horn admits a filler $c \colon \Del^n \to \qC$ and the triple $(a,b,c)$ defines an $n$-simplex in $\coll(f,g)$ that solves the lifting problem.
\end{proof}
\fi

We write $\coll(f,\qB)$ for the collage of $f \colon \qA \to \qB$ with the identity on $\qB$.

\begin{lem}[X.F.5.4]\label{lem:collage-qcat-cocartesian} For any $f \colon \qA \to \qB$, the map $\rho \colon \coll(f,\qB) \to \Del^1$ is a cocartesian fibration. \ifnewadd\qed\fi
\end{lem}
\ifnewcut\else
\begin{proof} To prove the claim, we need only specify cocartesian lifts of the non-degenerate 1-simplex of $\Del^1$ and demonstrate that these edges have the corresponding universal property. To that end, for any vertex $a \in \qA_0$, let $\chi_a \colon \Del^1 \to \coll(f,\qB)$ be the 1-simplex \[ \chi_a := (a \colon \Del^0 \to \qA, fa \colon \Del^0 \to \qB, fa \cdot \degen^0 \colon \Del^1 \to \qB),\] defined by the degenerate edge at $fa \in \qB_0$ lying over the 1-simplex in $\Del^1$. To show that $\chi_a$ is $\rho$-cocartesian, we must construct fillers for any left horn
\[
\xymatrix{ \Del^{\fbv{0,1}} \ar@/^3ex/[rr]^{\chi_a} \ar[r] & \Horn^{n,0} \ar[d] \ar[r] & \coll(f,\qB) \ar[d]^\rho \\ & \Del^n \ar@{-->}[ur] \ar[r]_\beta & \Del^1}\] whose initial edge is $\chi_a$. Note that this condition implies that the bottom map $\beta \colon [n] \to [1]$ carries $0$ to $0$ and the remaining vertices to $1$. The map $\Horn^{n,0} \to \coll(f,\qB)$ defines a horn $\Horn^{n,0} \to \qB$ in the quasi-category $\qB$ whose first edge is degenerate. By Joyal's lemma about filling ``special outer horns,'' such horns admit a filler $b \colon \Del^n \to \qB$ and the triple
\[ (a \colon \Del^0 \to \qA, b \cdot \face^0 \colon \Del^{n-1} \to \qB, b \colon \Del^n \to \qB)\] defines an $n$-simplex in $\coll(f,\qB)$ that solves the lifting problem.
\end{proof}
\fi

\begin{prop}[X.F.5.5]\label{prop:oplax-colimit-of-functor}
For any $f \colon \qA \to \qB$ between quasi-categories, the collage $\coll(f,\qB)$ defines the oplax colimit of $f$ in $\qCat$. That is $\coll(f,\qB)$ defines a cone under the pushout diagram
\[ \xymatrix{ \qA \ar[r]^f \ar@{^(->}[d]_{\id\times\face^0} & \qB \ar@{^(->}[d]  \ar@{^(->}@/^/[ddr] \\ \qA \times \Del^1 \ar[r] \ar@/_/[drr]_h & P \poexcursion \ar@{-->}[dr]^k  \\ & & \coll(f, {\qB})}\] so that the induced map $k$ is inner anodyne, and in particular a weak equivalence in the Joyal model structure.\ifnewadd\qed\fi
\end{prop}
\ifnewcut\else
\begin{proof} The map $k$ is a quotient of the map $h$, which has the following explicit description. For each 
 $n$-simplex $(a,\alpha) \colon \Del^n \to \qA \times \Del^1$ define $i \defeq |\alpha^{-1}(0)| -1$, so that $-1 \leq i \leq n$. Then $h$ carries $(a,\alpha)$ to the $n$-simplex of $\coll(f,\qB)$ corresponding to the triple
\[ (a\vert_{\fbv{0,\ldots, i}} \colon \Del^i \to \qA, fa\vert_{\fbv{i+1,\ldots,n}} \colon \Del^{n-i-1} \to \qB, fa \colon \Del^n \to \qB).\] Note that the composite $\rho h \colon \qA \times \Del^1 \to \Del^1$ is the projection.

It remains to present $k$ as a sequential composite of pushouts of coproducts of inner horn inclusions. To do so, first note that
\[ \coll(f,\qB)_n = \qA_n \coprod \qA_{n-1} \times_{\qB_{n-1}} \qB_n \coprod \cdots \coprod \qA_0 \times_{\qB_0} \qB_n \coprod \qB_n\] 
where each map $\qB_n \to \qB_i$ is the initial face map corresponding to $\Del^{\fbv{0,\ldots, i}} \inc\Del^n$. From the perspective of this decomposition, $P_n$ is the subset containing the sets $\qA_n$ and $\qB_n$ and the subset of $\qA_i \times_{\qB_i} \qB_n$ whose component in $\qB_n$ is in the image of $f$. The $n$-simplices of $\coll(f,\qB)$ that remain to be attached correspond to elements of $\qA_i \times_{\qB_i} \qB_n$, for $0 \leq i < n$, that are not in the image of $f$ in the sense just discussed. Note in particular that $k\colon P_0 \inc \coll(f,\qB)_0$ is an isomorphism and $k \colon P_n \inc \coll(f,\qB)_n$ is an injection for all $n \geq 1$.

To enumerate our attaching maps, we start with the collection of non-degenerate $n$-simplices of $\coll(f,\qB)$ for $n \geq 1$ that are not in the image of $f$ and remove also those elements of $\qA_i \times_{\qB_i} \qB_n$ whose components $b \in \qB_n$ are in the image of the degeneracy map $\degen_i \colon \qB_{n-1} \to \qB_n$. Partially order this set of simplices first in the order of increasing $n$ and the in order of increasing index $i$; that is we lexicographically order the collection of pairs $(n,i)$ with $n \geq 1$ and $0 \leq i < n$. We will filter the inclusion $P\inc\coll(f,\qB)$ as
\[ P \inc P_{(1,0)} \inc P_{(2,0)} \inc P_{(2,1)} \inc P_{(3,0)} \inc \cdots \inc P_{(n,i)} \inc \cdots \inc \colim\cong \coll(f,\qB)\]
where the simplicial set $P_{(n,i)}$ is built from the previous one by a pushout of a coproduct of inner horns indexed by the set of $n$-simplices $(a,b) \in \qA_i \times _{\qB_i} \qB_n$ with $b$ not in the image of $f$ or $\degen_i$. The filler for the horn indexed by $(a,b)$ will attach this $n$ simplex to $\qB_n$ as the missing face of the horn and also the $n+1$ simplex $(a,b\sigma^i) \in \qA_i \times_{\qB_i} \qB_{n+1}$.

Consider a simplex $(a,b) \in \qA_i \times_{\qB_i} \qB_n$ with $b$ not in the image of $f$ or $\degen_i$. Define a horn
\[
\xymatrix{ \Horn^{n+1,i+1} \ar[r] \ar@{^(->}[d] & P_{(n,i)} \ar@{^(->}[d] \\ \Del^{n+1} \ar[r]_-{(a,b\degen^i)} & \coll(f,\qB)}\] For each $0 \leq j < i+1$, the $\face^j$-face of the $n+1$ simplex $(a,b\sigma^i)$ is the $n$-simplex $(a\face^j,b\sigma^i\face^j)$, which lies in $P_{(n,i-1)}$ or in $\qB \inc P$ in the case $i=0$. For each $i+1 < j \leq n+1$, the $\face^j$-face of the $n+1$ simplex $(a,b\sigma^i)$ is the $n$-simplex $(a, b\sigma^i\face^j) = (a,b\face^{j-1}\degen^{i}) \in \qA_i \times_{\qB_i} \qB_n$, which was previously attached to $P_{(n-1,i)}$. So the $\Horn^{n+1,i+1}$ indeed maps to $P_{(n,i)}$, permitting an inductive construction of the next simplicial set in this sequence as the pushout
\[ \xymatrix{\coprod\limits_{\sim} \Horn^{n+1,i+1} \ar[r] \ar@{^(->}[d] & P_{(n,i)} \ar[d] \\ \coprod\limits_{\sim} \Del^{n+1} \ar[r] & P_{(n,i)+1}\poexcursion}
\] where $ P_{(n,i)+1}$ equals  $P_{(n+1,0)}$ in the case $i=n-1$ and $P_{(n,i+1)}$ otherwise.
\end{proof}
\fi  

\begin{cor}[X.F.5.6]\label{cor:lurie-adjunction} Consider a pair of functors between quasi-categories $f \colon \qA \to \qB$ and $u \colon \qB \to \qA$. Then $f$ is left adjoint to $u$ if and only if the collages $\coll(f,\qB)$ and $\coll(\qA, u)$ are equivalent under $\qB + \qA$ and over $\Del^1$. \ifnewadd\qed\fi
\end{cor}
\ifnewcut\else
\begin{proof}
First suppose that $\coll(f,\qB) \simeq \coll(\qA, u)$  under $\qB + \qA$ and over $\Del^1$. By Lemma \ref{lem:collage-qcat-cocartesian} this means that the map $\coll(f,\qB) \to \Del^1$ is both a cocartesian and a cartesian fibration, a \emph{bifibration} in the terminology of \S\refIV{sec:cartesian}. By Proposition \refIV{prop:bifibration-adjunction} it follows that the 1-arrow in $\Del^1$ from $0$ to $1$ induces an adjunction between the fibers $\qA$ and $\qB$. By inspection of that proof, the left adjoint functor so-constructed in the case of the bifibration $\coll(f,\qB) \to \Del^1$ is $f$; substituting the equivalent bifibration $\coll(\qA,u) \to \Del^1$, we see that the right adjoint is equivalent to $u$.

For the converse, we work in the opposite $\infty$-cosmos $\qCat\op$, an $\infty$-cosmos in which ``not all objects are cofibrant,'' as described in Observation \refIV{obs:duals-of-cosmoi}. In that context, Proposition \ref{prop:oplax-colimit-of-functor} proves that $\coll(f,\qB)$ and $\coll(\qA,u)$ construct the contravariant and covariant comma objects associated to the functors $f$ and $u$. If $f \dashv u$ in $\qCat$ then these functors are also adjoint in $\qCat\op$ and Proposition \refI{prop:adjointequiv} then proves that the commas $\coll(f,\qB)$ and $\coll(\qA,u)$ are equivalent under $\qB + \qA$. By construction, this equivalence also lies over $\Del^1$.
\end{proof}
\fi

\subsection{The comprehension construction}\label{ssec:comprehension}

In this section we review the \emph{comprehension construction} from \cite{RiehlVerity:2017cc}. It constructs, for any cocartesian fibration $p \colon \qE \tfib\qB$ of quasi-categories, a ``straightening,'' which has the form of a simplicial functor $c_p \colon \gC\qB \to\qCat$ that sends each vertex $b \in \qB$ to the fiber $\qE_b$. It also constructs a canonical lax cocone $\ell^\qE \colon \gC[\qB\join\Del^0] \to \qCat$ of shape $\qB$ under this diagram with nadir $\qE$.

\ifnewcut\else
The underlying mechanics of the comprehension construction are reviewed in \S\ref{sec:cocart-cocones} and the comprehension construction itself is given in \S\ref{sec:comprehension-review}.

\subsection{Cocartesian transformations between lax cocones}\label{sec:cocart-cocones}
\fi 

Corollary \ref{cor:collage-bijection} tells us that the collage of a weight $W$ realizes the shape of $W$-cocones. Applying this result to the weights for oplax colimits introduced in Definition \ref{defn:oplax-colimit-weight}, we obtain the following definition of a \emph{lax cocone}.

\begin{defn}[{lax cocones \refVI{defn:lax-cocone}}]
  Suppose that $X$ is a simplicial set. Then a \emph{lax cocone of shape $X$\/} in $\SSet$ is defined
  to be a simplicial functor $\ell^B\colon\gC[X\join\Del^0]\to\SSet$
    \begin{equation*}
    \xymatrix@=1.5em{
      & {\catone+\gC{X}}\ar@{^(->}[dl]\ar[dr]^{\langle B,B_\bullet \rangle} & \\
      {\gC[X\join\Del^0]}\ar[rr]_-{\ell^B} && {\SSet}
    }
  \end{equation*}
  The restriction of a lax cocone
  $\ell^B\colon\gC[X\join\Del^0]\to\eK$ to a functor $B_{\bullet}\colon\gC{X}\to \SSet$ is
  called its \emph{base\/}. We say
  that $\ell^B$ is a lax cocone \emph{under the diagram\/} $B_\bullet$; the object
  $B\in\sSet$ obtained by evaluating $\ell^B$ at the object $\top$ is called the
  \emph{nadir} of that lax cocone.
  \end{defn}

  \ifnewcut\else
  \begin{rmk}
  In the original Definition \refVI{defn:lax-cocone}, the target was required to be a quasi-categorically enriched category $\eK$ and the base of a lax cocone was required to factor through
  through the inclusion $g_*\eK\subseteq \eK$ of the maximal Kan complex enriched subcategory. The point of this requirement was so that the transpose of the base diagram defined a diagram $X \to \nrvhc g_*\eK$ valued in the large quasi-category of objects and morphisms in $\eK$. But in this paper we will frequently  consider lax cocones whose codomain is $\SSet$, in which case the requirement that the base diagram factor through the local groupoidal core won't necessarily make sense and should be dropped. We shall also be more interested in the totality of the data of a lax cocone than in studying the base diagram in isolation.
  \end{rmk}
  
  \begin{defn}[{canonical lax cocones \refVI{lem:constant-lax-cones}}]\label{defn:canonical-cocone}\fi
  \ifnewadd
  \begin{ex}[{canonical lax cocones \refVI{lem:constant-lax-cones}}]\label{ex:canonical-cocone} \fi
    For any simplicial set $X$, there exists a lax cocone 
      \begin{equation*}
    \xymatrix@=1.5em{
      & {\catone+\gC{X}}\ar@{^(->}[dl]\ar[dr]^{\langle X,1 \rangle} & \\
      {\gC[X\join\Del^0]}\ar[rr]_-{k^X} && {\SSet}
    }
  \end{equation*}
  whose base is constant at the terminal quasi-category $1$ and whose nadir is $X$ that we refer to as the \emph{canonical $X$-shaped lax cocone}. 
  \ifnewadd 
\end{ex}\fi
  
\ifnewcut\else
To define $k^X \colon \gC[X \join\Del^0] \to \SSet$, it remains to define the images in $\Fun(1,X)\cong X$ of the atomic $r$-arrows with domain a vertex $x_0$ in $X$ and with codomain $\top$.   Applying Proposition \ref{prop:gothic-C}, each atomic $r$-arrow from $x_0$ to $\top$ corresponds to a bead $((x, \top), T^\bullet)$ represented by some non-degenerate $n$-simplex $x \in X$ whose initial vertex is $x_0$ and atomic $r$-arrow $T^\bullet \colon 0 \to n+1$ in $\gC\Del^{n+1}$ defined by
  \[ \{0,n+1\} = T^0 \subset T^1 \subset \cdots \subset T^r = [0,n+1].\]
  Define $\mu \colon \Del^r \to \Del^n$ by $\mu(i) = \max (T^i \backslash \{n+1\})$. Then the image of the bead $((x, \top), T^\bullet)$ is the $r$-simplex $x \cdot \mu \in X$ in the hom-quasi-category $X \cong \Fun(1,X)$ from the image of the domain vertex $x_0$ of $x$.
\end{defn}
\fi

\begin{obs}[{whiskering lax cocones \refVI{obs:whiskering-cocone}}]\label{obs:whiskering-cocone}
Let $\ell^A \colon \gC[X \join\Del^0] \to \SSet$ be a lax cocone with base diagram $A_\bullet \colon \gC{X} \to \SSet$ and nadir $\ell^A_\top = A$, and let $f \colon A \to B$ be any map of simplicial sets. Then there is a \emph{whiskered lax cocone} $f \cdot \ell^A \colon \gC[X \join\Del^0] \to \SSet$ with the same base diagram $A_\bullet \colon \gC{X} \to \SSet$ and with nadir $B$, whose components from a vertex $x \in X$ to $\top$ are defined by whiskering with $f$:
\[ \Fun_{\gC[X\join\Del^0]}(x,\top) \xrightarrow{ \ell^X_{x,\top}} \Fun(A_x,A) \xrightarrow{f \circ -} \Fun(A_x,B)\]
\end{obs}

\begin{lem}\label{lem:canonical-whiskering} 
For any map of simplicial sets $f \colon Y \to X$, the canonical lax cocone of shape $X$ restricts along $\gC[f \join\id] \colon \gC[Y\join\Del^0] \to \gC[X\join\Del^0]$ to the whiskered composite 
 \[
\begin{tikzcd}[row sep=small, column sep=tiny] & \catone+ \gC{Y} \arrow[dl, hook'] \arrow[rr, "\catone+\gC{f}"] & & \catone+ \gC{X}  \arrow[dl, hook'] \arrow[dr, "{\langle X,1 \rangle}"] & \arrow[dr, phantom, "="] & & \catone+\gC{Y} \arrow[dl, hook'] \arrow[dr, "{\langle X,1\rangle}"] \\ \gC[Y\join\Del^0] \arrow[rr, "{\gC[f \join\id]}"'] & & \gC[X\join\Del^0]\arrow[rr, "k^X"'] && \SSet &  \gC[Y\join\Del^0] \arrow[rr, "f \cdot k^Y"']  && \SSet
\end{tikzcd}
\]
of the canonical lax cocone of shape $Y$ with $f \colon Y \to X$.
\end{lem}
\begin{proof}
By direct verification from   \ifnewcut\else Definition \ref{defn:canonical-cocone}\fi\ifnewadd  Lemma \refVI{lem:constant-lax-cones} \fi and Observation \ref{obs:whiskering-cocone}.
\end{proof}

 \S\refVI{sec:cocones} introduces a mechanism for producing new lax cocones from given ones: namely as domain components of cocartesian cocones over a given codomain lax cocone. 

\begin{defn}[cocartesian cocones {\refVI{defn:cocart-cocone}}]
    Suppose we are given a simplicial set $X$ and lax cocones $\ell^E,\ell^B
  \colon\gC[X\join\Del^0]\to \SSet$ of shape $X$
  with bases $E_\bullet$ and $B_\bullet$ respectively. Suppose also that we are given a
  simplicial natural transformation
  \begin{equation*}
    \xymatrix@R=0em@C=8em{
      {\gC[X\join\Del^0]}\ar@/_2ex/[]!R(0.5);[r]_{\ell^B}\ar@/^2ex/[]!R(0.5);[r]^{\ell^E}
      \ar@{}[r]|-{\textstyle \Downarrow p} & {\SSet.}
    }
  \end{equation*}
  Then we say that the triple $(\ell^E,\ell^B,p)$ is a \emph{cocartesian cocone} if
  and only if
  \begin{enumerate}[label=(\roman*)]
  \item\label{itm:cocart-cocone-i} the nadir of the natural transformation $p$, that being its component
    $p\colon \qE\tfib \qB$ at the object $\top$, is a cocartesian fibration between quasi-categories
  \item\label{itm:cocart-cocone-ii} for all $0$-simplices $x\in X$ the naturality square is a pullback, and
    \begin{equation*}
      \xymatrix{
        {E_x}\ar@{->>}[d]_{p_x}\ar[r]^{\ell^E_{{x}}} \pbexcursion & {\qE}\ar@{->>}[d]^{p} \\
        {B_x}\ar[r]_{\ell^B_{{x}}} & {\qB}
      }
    \end{equation*}
      \item\label{itm:cocart-cocone-iii} for all non-degenerate $1$-simplices $f\colon x\to y\in X$ the $1$-arrow   is $p$-cocartesian.
    \begin{equation*}
      \xymatrix@C=5em@R=1em{
        {E_x}\ar[dd]_{E_{f}}\ar[dr]^{\ell^E_{{x}}} & \\
        {}\ar@{}[r]|(0.35){\textstyle\Downarrow \ell^E_{{f}}} & {\qE} \\
        {E_y}\ar[ur]_{\ell^E_{{y}}} & 
      }
    \end{equation*}

  \end{enumerate}
\ifnewcut\else  In this situation we also say that the pair $(\ell^E,p)$ defines a \emph{cocartesian
  cocone over} $\ell^B$.\fi
  \end{defn}

\begin{lem}[{pullbacks of cocartesian cocones \refVI{lem:cocart-cocone-pb}}]\label{lem:cocart-cocone-pb} Suppose given:
\begin{itemize}
\item a pullback diagram of quasi-categories in which $p$ and $q$ are cocartesian fibrations;
\begin{equation}\label{eq:cocart-cocone-pb} \xymatrix{ \qF \ar@{->>}[d]_q \ar[r]^g \pbexcursion & \qE \ar@{->>}[d]^p \\ \qA \ar[r]_f & \qB}\end{equation} 
\item a lax cocone $\ell^A \colon \gC[X\join\Del^0] \to \eK$  with nadir $\qA$; and
\item a cocartesian cocone $(\ell^E,\ell^B,p)$ whose nadir is $p \colon \qE \tfib \qB$ and whose codomain cocone $\ell^B = f \cdot \ell^A$ is obtained from the lax cocone $\ell^A$ by whiskering with $f\colon \qA \to \qB$.
\end{itemize}
Then there is a cocartesian cone $(\ell^F,\ell^A,q)$ whose codomain is $\ell^A$, whose nadir is $q \colon \qF \tfib \qA$, and whose domain component is a lax cocone $\ell^F$ that whiskers with $g$ to the lax cone $\ell^E = g \cdot \ell^F$.  \qed
\end{lem}

Conversely, a cocartesian cocone $(\ell^F,\ell^A,q)$ with nadir $q \colon \qF \tfib \qA$ can be whiskered with a pullback square \eqref{eq:cocart-cocone-pb} to define a cocartesian cocone $(g\cdot \ell^F, f\cdot \ell^A, p)$ with nadir $p \colon \qE \tfib \qB$ and whose domain and codomain are whiskered lax cocones as defined in Observation \ref{obs:whiskering-cocone}.

\begin{rmk}\label{rmk:arbitrary-pullbacks}
If the map $f$ of Lemma \ref{lem:cocart-cocone-pb} is replaced by any map  of simplicial sets $f \colon X \to \qB$, whose domain is not necessarily a quasi-category, it is still possible to pull back the data of a cocartesian cocone $(\ell^E,\ell^B,p)$ whose codomain lax cocone $\ell^B = f \cdot \ell^X$ is obtained by whiskering a lax cocone with nadir $X$. This constructs a simplicial natural transformation  $(\ell^F,\ell^X,q)$ whose nadir is the pullback $q \colon F \to X$ of $p$ along $f$. Since this is not a map between quasi-categories, it does not really make sense to call it a cocartesian fibration. Nonetheless, this construction produces a lax cocone $\ell^F$ of shape $X$, which will have some utility. See Remark \ref{rmk:generalized-unstraightening}.
\end{rmk}

\ifnewcut\else
\subsection{The comprehension construction}\label{sec:comprehension-review}
\fi 

A cocartesian fibration $p \colon \qE \tfib \qB$ between quasi-categories has a ``straightening'' called the \emph{comprehension functor} $c_p \colon \gC\qB \to \qCat$, a homotopy coherent diagram of shape $\qB$ that sends each vertex $b$ to the fiber $\qE_b$ of $p$ over $b$. This arises as the base diagram of the domain of a cartesian cocone over the canonical $\qB$-shaped lax cocone. 

\begin{thm}[{\refVI{defn:basic-comprehension}}]\label{thm:comprehension}
For any cocartesian fibration $p\colon \qE\tfib \qB$ of quasi-categories, there is a cocartesian cocone 
  \begin{equation*}
    \xymatrix@R=0em@C=8em{
      {\gC[\qB\join\Del^0]}\ar@/_2ex/[]!R(0.5);[r]_{k^\qB}\ar@/^2ex/[]!R(0.5);[r]^{\ell^\qE}
      \ar@{}[r]|-{\textstyle \Downarrow p} & {\qCat.}
    }
  \end{equation*}
  of shape $\qB$ in $\qCat$ with nadir $p \colon \qE\to\qB$ over the canonical lax cocone $k^\qB$. The base of the domain component defines the comprehension functor $c_p$, which acts on an object $b \colon 1 \to \qB$ of $\gC\qB$ by forming the pullback 
  \[
      \xymatrix@=1.7em{
      {\qE_b}\pbexcursion\ar[r]^{\ell^\qB_{{b}}}\ar@{->>}[d]_-{p_b} &
      {\qE}\ar@{->>}[d]^-p \\
      {1}\ar[r]_b & {\qB}
    }
    \] and acts on 1-arrows $f \colon a \to b$ of $\qB$ 
 by factoring the codomain of a $p$-cocartesian lift $\ell^\qE_{{f}}$ of $f$ through the pullback at the front of the diagram:
  \begin{equation}\label{eq:comprehension-on-1-arrows}
        \begin{xy}
      0;<1.4cm,0cm>:<0cm,0.75cm>::
      *{\xybox{
          \POS(1,0)*+{1}="one"
          \POS(0,1)*+{1}="two"
          \POS(3,0.5)*+{\qB}="three"
          \ar@{=} "one";"two"
          \ar@/_5pt/ "one";"three"_{b}^(0.1){}="otm"
          \ar@/^10pt/ "two";"three"^{a}_(0.5){}="ttm"|(0.325){\hole}
          \ar@{=>} "ttm"-<0pt,7pt> ; "otm"+<0pt,10pt> ^(0.3){f}
          \POS(1,2.5)*+{\qE_{b}}="one'"
          \POS(1,2.5)*{\pbcorner}
          \POS(0,3.5)*+{\qE_{a}}="two'"
          \POS(0,3.6)*{\pbcorner}
          \POS(3,3)*+{\qE}="three'"
          \ar@/_5pt/ "one'";"three'"_{\ell^\qE_{{b}}}
          \ar@/^10pt/ "two'";"three'"^{\ell^\qE_{{a}}}_(0.55){}="ttm'"
          \ar@{->>} "one'";"one"_(0.325){p_{b}}
          \ar@{->>} "two'";"two"_{p_{a}}
          \ar@{->>} "three'";"three"^{p}
          \ar@{..>} "two'";"one'"_*!/^2pt/{\scriptstyle \qE_f}
          \ar@{..>}@/_10pt/ "two'";"three'"^(0.44){}="otm'"
          \ar@{=>} "ttm'"-<0pt,4pt> ; "otm'"+<0pt,4pt> ^(0.3){\ell^\qE_{{f}}}
        }}
    \end{xy}
  \end{equation}
These cocartesian lifts define 
  components of the lax cocone
        \begin{equation*}
    \xymatrix@=1.5em{
      & {\catone+\gC{\qB}}\ar@{^(->}[dl]\ar[dr]^{\langle \qE,c_p \rangle} & \\
      {\gC[\qB\join\Del^0]}\ar[rr]_-{\ell^\qE} && {\qCat}
    }
  \end{equation*}
  with nadir $\qE$ under the comprehension functor. \qed 
  \end{thm}

By Observation \refVI{obs:unique-comprehension}, any pair of lax cocones that arise as the domain of a cocartesian cocone over a common codomain define vertices in a contractible Kan complex and are in particular equivalent as diagrams. This is used to prove the following result relating comprehension functors with pullbacks.

\begin{prop}[{comprehension and change of base \refVI{prop:comprehension-cob}}]\label{prop:comprehension-cob}
  Suppose that we are given a pullback
  \begin{equation*}
    \xymatrix@=2em{
      {\qF}\pbexcursion\ar@{->>}[d]_{q}\ar[r]^g &
      {\qE}\ar@{->>}[d]^{p} \\
      {\qA}\ar[r]_{f} & {\qB}
    }
  \end{equation*}
  of quasi-categories in which $p$ and thus $q$ are cocartesian
  fibrations. Then the diagrams\[
  \gC\qA \xrightarrow{c_q} \qCat\qquad \mathrm{and} \qquad \gC\qA \xrightarrow{\gC{f}} \gC\qB \xrightarrow{c_p} \qCat\] are connected by a homotopy coherent natural isomorphism. \qed
\end{prop}

\begin{rmk}\label{rmk:generalized-unstraightening}
If $\qA$ is not a quasi-category, it is not possible to directly construct the comprehension functor for the pullback of $p$ along $f$. However, by  Lemmas \ref{lem:canonical-whiskering} and Remark \ref{rmk:arbitrary-pullbacks},  the cocartesian cocone over the canonical $\qB$-shaped lax cocone can be pulled back along any map of simplicial sets $f \colon X \to \qB$ to define a cocartesian cocone over the canonical $X$-shaped lax cocone. Thus, a posteriori, we can think of the base of the lax cocone
\[ \gC[X \join\Del^0] \xrightarrow{\gC[f\join\id]} \gC[\qB\join\Del^0] \xrightarrow{\ell^\qE} \qCat\] as defining a comprehension functor for the pullback of $p \colon \qE \tfib \qB$ along $f \colon X \to \qB$. 
\end{rmk}


%% file: weighted.tex

\section{Weighted colimits in simplicial categories}\label{sec:weighted}

Our aim in this section is to define the \emph{oplax colimit} of a homotopy coherent diagram $\gC{X} \to \SSet$ indexed by the homotopy coherent realization of a simplicial set $X$, our name for the left adjoint of the homotopy coherent nerve functor $\hN \colon \sCat\to\sSet$. Oplax colimits are introduced as particular \emph{weighted colimits}, where the weights in question are simplicial functors that describe the shape of lax cocones. In \S\ref{sec:simp-collage}, we review the general rubric of weighted colimits and explain how these cocone shapes may be presented as simplicial \emph{collages}. In \S\ref{sec:flexible-weight}, we highlight a special class of \emph{flexible weights} that have useful homotopical properties. Finally, in \S\ref{sec:oplax}, we define the weights for oplax colimits as collages and observe that such weights are flexible, a fact that will be exploited in our future work.

Some of this material was previously discussed in \S\refVII{sec:flexible}, where the ``oplax'' weights were called ``pseudo'' weights. See Remark \ref{rmk:oplax-vs-pseudo} for an explanation of this contrast in nomenclature. Collages for weights for limits made an appearance in \S\refVII{ssec:collage} and we refer the proofs of a few of the results appearing below to there, but we reintroduce this construction here to clarify the details in the dual case and because we will require a more extensive analysis of collages than we did in \cite{RiehlVerity:2018rq}.

\subsection{Weighted colimits and collages}\label{sec:simp-collage}

In a simplicially enriched category, the appropriately general notion of colimit allows for the specification of any particular ``shape'' of cone under the diagrams being considered. This specification is given by a simplicial functor referred to as a \emph{weight} for the colimit. 

\begin{defn}[weights for simplicial colimits]\label{defn:simp-weight}
  Suppose $\eD$ is a small simplicial category, which we think of as a
  diagram shape. Then a \emph{weight} on $\eD$ is a simplicial functor $W\colon
  \eD\op\to\SSet$. For any diagram $F\colon \eD\to\eK$ of shape $\eD$ valued in a
  simplicial category $\eK$, a \emph{$W$-cocone} with \emph{nadir} an object
  $e\in \eK$ is a simplicial natural transformation $\iota\colon W\to
  \Fun_{\eK}(F(-),e)$. We say that the $W$-cocone $\iota$ \emph{displays $e$ as
    a $W$-colimit of $F$} if and only if for all objects $e'\in\eK$ the
  simplicial map 
  \begin{equation*}
    \xymatrix@R=0em@C=5em{
      \Fun_{\eK}(e,e') \ar[r]^-{\cong} & \Fun_{\SSet^{\eD\op}}(W,\Fun_{\eK}(F(-),e'))
    }
  \end{equation*}
  given by post-composition with $\iota$ is an isomorphism. 
 
Many notations are common for the nadir of a weighted colimit cone; here we write $\colim^W\!F$ for the
  colimit of $F$ weighted by $W$. When these objects exist for all weights and diagrams in
  $\eK$ then we may extend $\colim$ to a simplicial bifunctor:
  \begin{equation*}
    \xymatrix@R=0em@C=7em{
      \SSet^{\eD\op}\times\eK^{\eD}\ar[r]^-{\colim} & \eK
    }
  \end{equation*} that is cocontinuous in both variables.
\end{defn}

A simplicial functor $W\colon\eD\op\to\SSet$ may otherwise be described as
comprising a family of simplicial sets $\{Wd\}_{d\in\obj(\eD)}$ along with right
actions
\begin{equation}\label{eq:action-of-weight}
  \xymatrix@R=0em@C=6em{
    Wd'\times\Fun_{\eD}(d,d')\ar[r]^-{*} & Wd
  }
\end{equation}
of the hom-spaces of $\eD$ which must collectively satisfy the customary
axioms with respect to the identities and composition of $\eD$. This
description leads us to define a simplicially enriched category $\coll(W)$,
called the \emph{collage} of $W$.

\begin{defn}[collages]\label{defn:collage}
  For any weight $W \colon \eD\op\to\SSet$, the \emph{collage} of $W$ is a
  simplicial category $\coll(W)$ that contains $\eD$ as a full simplicial
  subcategory along with precisely one extra object $\top$ whose endomorphism
  space is the point. The function complexes $\Fun_{\coll(W)}(\top,d)$ are all taken
  to be empty and we define:
  \begin{equation*}
    \Fun_{\coll(W)}(d,\top) \defeq Wd\mkern40mu \text{for objects $d\in\eD$.}
  \end{equation*}
  The composition operations between hom-spaces in $\eD$ and those with
  codomain $\top$ are given by the actions depicted in~\eqref{eq:action-of-weight}.
\end{defn}

In the statement of the following result,
$\sSet^{\eD\op}$ denotes the underlying category of the simplicially
  enriched category $\SSet^{\eD\op}$.

\begin{prop}[{collage adjunction, \refVII{prop:collage-adjunction}}]\label{prop:collage-adjunction}$\quad$
\begin{enumerate}[label=(\roman*)]
\item\label{itm:collage-fun} The collage construction defines a fully faithful functor 
  \begin{equation*}
    \xymatrix@R=0em@C=6em{
      \sSet^{\eD\op}\ar[r]^-{\coll} & \prescript{\catone+\eD/}{}{\sCat}
    }
  \end{equation*}
  from the category of $\eD$-indexed weights to the category of simplicial categories under $\catone+\eD$ whose essential image is comprised of those $\langle
  e,F\rangle\colon \catone+\eD\to\eK$ that are bijective on objects, fully faithful
  when restricted to $\eD$ and $\catone$, and have the property that there are no maps in $\eK$ from $e$ to the image of $F$.
\item\label{itm:collage-adj} The collage functor admits a right adjoint
\[ \adjdisplay \coll -| \wgt : \prescript{\catone+\eD/}{}{\sCat} ->      \sSet^{\eD\op}.\] which carries a pair 
$\langle e,F
  \rangle\colon\catone+\eD\to\eK$ to the weight $\Fun_{\eK}(F(-),e)\colon\eD\op\to\SSet$.
\end{enumerate}
\end{prop}

  This adjunction has a useful and important interpretation:
  
  \begin{cor}[{\refVII{cor:collage-bijection}}]\label{cor:collage-bijection} The collage   $\coll(W)$ of a weight
 realises the shape of $W$-cocones, in the sense that simplicial functors 
\[  G \colon \coll(W) \longrightarrow \eK\]
stand in bijection to $W$-cocones under the diagram $G\vert_{\eD}$ with nadir $G(\top)$. \qed
\end{cor}

\begin{lem}\label{lem:lan-of-weight}
  For any simplicial functor $I \colon \eD \to \eC$, weight $W \colon
  \eD\op\to\SSet$, and diagram $G \colon \eC\to\eK$, we have an isomorphism
  \begin{equation*}
    \colim\nolimits^W\!GI \cong \colim\nolimits^{\lan_{I}W}\!G
  \end{equation*}
  where the colimit on one side exists if and only if the one on the other
  does.
\end{lem}

\begin{proof}
  Simplicial left Kan extension provides an adjunction
  \begin{equation*}
    \adjdisplay \lan_{I} -| -\circ I : \SSet^{\eC\op} -> \SSet^{\eD\op}.
  \end{equation*}
  In particular
  \[
    \Fun_{\SSet^{\eC\op}}(\lan_IW,\Fun_{\eK}(G(-),e)) \cong
    \Fun_{\SSet^{\eD\op}}(W,\Fun_{\eK}(GI(-),e))
  \]
  which shows that $\colim^{\lan_IW}G$ and $\colim^W\!GI$  have
  the same defining universal property.
\end{proof}

\begin{lem}\label{lem:lan-collage-pushout}
  Left Kan extension of $W\colon \eD\op\to\SSet$ along a simplicial functor
  $I\colon \eD\to\eC$ gives rise to a pushout square
  \begin{equation*}
    \xymatrix@=2em{
      {\catone+\eD}\ar[r]^{\catone+I}\ar@{^(->}[d] &
      {\catone+\eC}\ar@{^(->}[d] \\
      {\coll(W)}\ar[r] & {\coll(\lan_{I}W)}\poexcursion
    }
  \end{equation*}
  in the category of simplicial categories and simplicial functors.
\end{lem}

\begin{proof}
  By the defining universal property, a simplicial functor whose domain is the
  pushout of $\catone+\eD\inc\coll(W)$ along $\catone + I$ and whose codomain
  is $\eK$ is given by a pair of functors
  \[
    \gamma\colon \coll(W)\to \eK \qquad \mathrm{and} \qquad
    \langle e,G\rangle \colon\catone+\eC \to \eK.
  \]
  By Corollary \ref{cor:collage-bijection}, the simplicial functor $\gamma$
  represents a $W$-cocone in $\eK$ with nadir $e$ under the diagram
  $GI\colon\eD\to\SSet$. By Lemma \ref{lem:lan-of-weight}, such data
  equivalently describes a $\lan(W)$-shaped cocone under the diagram $G$ with
  nadir $e$. Applying Corollary \ref{cor:collage-bijection} again, we conclude
  that this pushout is given by the simplicial category $\coll(\lan_IW)$, as
  claimed.
\end{proof}

In ordinary unenriched category theory, the colimit cone under a $\eD$-shaped diagram may be formed as the left Kan extension along the inclusion $\eD\inc\eD\join\catone$ into the category $\eD\join\catone$ formed by freely adjoining a terminal object ``$\top$'' to $\eD$. The following lemma reveals that the collage plays the roll of the category $\eD\join\catone$ for weighted colimits, a perspective which we will return to in \S\ref{sec:oplax}.

\begin{lem}\label{lem:lan-along-collage} 
  The pointwise left Kan extension of any simplicial functor
  $F\colon\eD\to\eK$ along $I \colon \eD\inc\coll(W)$ exists if and only if
  the colimit $\colim^W\!F$ exists in $\eK$, and then
  $\lan_IF(\top)\cong\colim^W\!F$.
\end{lem}

\begin{proof}
  Since $\eD\inc\coll(W)$ is fully faithful, when the pointwise left Kan
  extension of any simplicial functor $F\colon\eD\to\eK$ along
  $\eD\in\coll(W)$ exists, it is displayed by an isomorphism:
  \begin{equation*}
    \xymatrix@=1.5em{
      & {\eD}\ar@{^(->}[dl]_I \ar@{}[d]|(.6){\cong}\ar[dr]^{F} & \\
      {\coll(W)}\ar[rr]_-{\lan_IF} && {\eK}
    }
  \end{equation*} 
  By Corollary \ref{cor:collage-bijection}, this data defines a $W$-cocone
  under $F\cong \lan_IF\circ I$ with nadir $\lan_IF(\top) \in \eK$. It is
  easy to verify that the universal property of the left Kan extension
  specializes to describe the universal property of the colimit cocone for
  $\colim^W\!F$, and conversely that the universal property of the weighted
  limit cocone implies the universal property of the left Kan extension.
\end{proof}

\subsection{Flexible weights as projective cell complexes}\label{sec:flexible-weight}

In order to understand the sense in which certain weighted colimits, including in particular the oplax colimits to be introduced below, are homotopically well
behaved, we recall some facts about weights and simplicial computads from
\S\refII{subsec:collage}:

\begin{defn}[flexible weights and projective cell complexes]\label{defn:proj-cell-cx}
  Suppose that $\eD$ is a simplicial category. Then a simplicial natural
  transformation of the form \[\boundary\Del^n\times \Fun_{\eD}(-,d)\inc
    \Del^n\times \Fun_{\eD}(-,d),\] for some $[n]\in\Del$ and object
  $d\in\eD$, is said to be the \emph{projective $n$-cell} associated with $d$.
  A natural transformation $\alpha\colon W\to V$ in $\SSet^{\eD\op}$ is a
  \emph{relative projective cell complex} if it may be constructed as a
  countable composite of pushouts of coproducts of projective cells. A weight
  $W$ in $\SSet^{\eD\op}$ is a \emph{flexible weight} if the map
  $!\colon\emptyset\to W$ is a relative projective cell complex, i.e., if $W$
  is a projective cell complex.
\end{defn}

Our interest in colimits weighted by flexible weights is due to the fact that
they are homotopically well behaved. We state the following result for diagrams
valued in simplicial sets, but its proof extends without change to pointwise
cofibrant diagrams valued in any model category enriched over the Joyal model
structure on simplicial sets. 

\begin{prop}[{\refII{prop:proj-wlim-homotopical}, \refVII{prop:flexible-weights-are-htpical}}]\label{prop:flexible-weights-are-htpical}$\quad$
  \begin{enumerate}[label=(\roman*)]
  \item For a flexible weight $W \colon \eD\op\to\SSet$ and any diagram $F \colon \eD\to\SSet$, $\colim^W\!F$ may be expressed as a countable composite of
    pushouts of coproducts of maps \[\boundary\Del^n\times Fd\inc\Del^n\times
      Fd.\]
  \item If $\alpha \colon F \to G$ is a simplicial natural transformation
    between two such diagrams whose components are weak equivalences in the
    Joyal model structure,
    then for any flexible weight $W$ the map \[\colim\nolimits^W\!\alpha\colon\colim\nolimits^W\!F\to\colim\nolimits^W\!G\] is a weak
    equivalence in the Joyal model structure.
  \end{enumerate}
\end{prop}



The collage construction defines a correspondence between flexible weights and \emph{simplicial computads}, a class of ``freely generated'' simplicial categories that define precisely the cofibrant objects \cite[\S 16.2]{Riehl:2014kx} in  the model structure due to  Bergner \cite{Bergner:2007fk}. 

\begin{defn}[simplicial computad]\label{defn:simplicial-computad}
A simplicial category $\eA$, regarded as a simplicial object $[n] \mapsto \eA_n$ in the category of categories with a common set of objects $\mathrm{ob}\eA$ and identity-on-objects functors, 
 is a \emph{simplicial computad} if and only if:
  \begin{itemize}
    \item each category $\eA_n$ of $n$-\emph{arrows} is freely generated by the reflexive directed graph of \emph{atomic} $n$-arrows, these being those arrows that admit no non-trivial factorizations, and if 
       \item if $f$ is an atomic $n$-arrow in $\eA_n$ and $\alpha\colon [m]\to[n]$ is a degeneracy operator in $\Del$ then the degenerated $m$-arrow $f\cdot\alpha$ is atomic in $\eA_m$.
  \end{itemize}
A simplicial category  $\eA$ is a simplicial computad if and only if all of its non-identity arrows $f$ can be expressed uniquely as a composite 
\begin{equation}\label{eq:computad-arrow-decomp}
  f = (f_1 \cdot \alpha_1) \circ (f_2 \cdot \alpha_2) \circ \cdots \circ (f_\ell \cdot \alpha_\ell)
\end{equation}
in which each $f_i$ is non-degenerate and atomic and each $\alpha_i\in\Del$ is a degeneracy operator.
\end{defn}

We have the following recognition principle for flexible weights on simplicial
computads,  a mild variant of Proposition~\refII{prop:projcofchar}, proven in \S\refVII{ssec:collage}.

\begin{prop}[{relating flexible weights and simplicial computads, \refVII{thm:flexible-collage}}]
  \label{prop:flexible-weight-as-computad}
  Suppose that $\eD$ is a simplicial computad and that $W\colon\eD\op\to\SSet$
  is a weight. Then $W$ is a flexible weight if and only if its collage
  $\coll(W)$ is a simplicial computad.
\end{prop}

\subsection{Oplax colimits}\label{sec:oplax}

Oplax colimits represent particular cones under a homotopy coherent diagram
$\gC{X} \to \eK$ indexed by a simplicial set $X$. In Definition \ref{defn:oplax-colimit-weight}, we
first present the collage construction that describes the shape of a lax cocone
and then use Proposition \ref{prop:collage-adjunction} to extract the
corresponding weight. To give a concise description of the collage that defines the oplax weight, we make use of the simplicial computad structure on the \emph{homotopy coherent realization} $\gC{X}$ of a simplicial set $X$, our term for the left adjoint to the homotopy coherent nerve
\[
\adjdisplay \gC -| \hN : \sCat -> \sSet.\]
We briefly review this material from  \S\refVI{sec:coherent-nerve}, which gives a more leisurely presentation with considerably more details.

  \begin{ex}[{homotopy coherent simplices as simplicial computads; \S\refVI{sec:htpy-coh-simplex}}]\label{ex:homotopy-coherent-simplex} Recall the simplicial category $\mathfrak{C}\Del^n$ whose objects are integers $0,1,\ldots, n$ and whose function complexes are the cubes \[ \Fun_{\mathfrak{C}\Del^n}(i,j) = \begin{cases} \Cube^{j-i-1} & i < j \\ \Del^0 & i = j \\ \emptyset & i > j \end{cases}\] Here we write $\Cube^k \defeq (\Del^1)^k$. For $i < j$, the vertices of $\Fun_{\mathfrak{C}\Del^n}(i,j)$ are naturally identified with subsets of the closed interval $[i,j] = \{ i \leq t \leq j\}$ containing both endpoints, a set whose cardinality is $j-i-1$; more precisely, $\Fun_{\mathfrak{C}\Del^n}(i,j)$ is the nerve of the poset with these elements, ordered by inclusion. Under this isomorphism, the composition operation  corresponds to the simplicial map 
  \[
  \xymatrix@R=1em{ \Fun_{\gC\Del^n}(i,j)\times\Fun_{\gC\Del^n}(j,k)\ar[r]^-\circ \ar@{}[d]|{\rotatebox{90}{$\cong$}}& \Fun_{\gC\Del^n}(i,j) \ar@{}[d]|{\rotatebox{90}{$\cong$}} \\
\Cube^{\times
    (j-i-1)}\times\Cube^{\times (k-j-1)}\ar[r] & \Cube^{\times (k-i-1)}}\]
     which
  maps the pair of vertices $T \subset [i,j]$ and $S \subset [j,k]$ to $T \cup S \subset [i,k]$.
  
  Again for $i < j$, an $r$-arrow $T^\bullet$ in $\Fun_{\gC\Del^n}(i,j)$ corresponds to a sequence
  \[ T^0 \subset T^1 \subset \cdots \subset T^r\] of subsets of $[i,j] = \{ i \leq t \leq j\}$ and is non-degenerate if and only if each of these inclusions are proper. The composite of a pair of $r$-arrows $T^\bullet \colon i \to j$ and $S^\bullet \colon j \to k$ is the levelwise union $T^\bullet \cup S^\bullet \colon i \to k$ of these sequences.
  
 From this description, it is easy to see that the simplicial category $\gC\Del^n$  is a simplicial computad (Lemma \refVI{lem:simplex-computad}), in which an $r$-arrow $T^\bullet$ from $i$ to $j$ is atomic if and only if the set $T^0 = \{i,j\}$; the only atomic $r$-arrows from $j$ to $j$ are identities. Geometrically,   the atomic arrows in each function complex $\Fun_{\gC\Del^n}(i,j) \cong \Cube^{j-i-1}$ are precisely those simplices that contain the initial vertex in the poset whose nerve defines the simplicial cube.
\end{ex}

If $X$ is a simplicial subset of $\Del^n$, then Lemma \refVI{lem:subcomputad-inclusion} tells us that 
 its homotopy coherent realisation
  $\gC{X}$ is a simplicial subcomputad of $\gC\Del^n$.

\begin{ex}[{homotopy coherent nerves of subsimplices \refVI{ex:subcomputad-of-simplex}}]\label{ex:subcomputad-of-simplex}
In particular: 
 \begin{enumerate}[label=(\roman*)]
  \item\label{itm:boundary} \textbf{boundaries:} The inclusion $\gC\boundary\Del^n\inc\gC\Del^n$ is the identity on objects and full on all function complexes except for the one from $0$ to $n$. The inclusion
         \[
  \xymatrix@R=1em{ \Fun_{\gC\boundary\Del^n}(0,n)\ar@{^(->}[r] \ar@{}[d]|{\rotatebox{90}{$\cong$}}& \Fun_{\gC\Del^n}(0,n)\ar@{}[d]|{\rotatebox{90}{$\cong$}} \\
\boundary\Cube^{n-1}\ar@{^(->}[r] & \Cube^{n-1}}\]
    is isomorphic to the cubical boundary inclusion, where $\boundary\Cube^k$ is the domain of the iterated Leibniz product
  $(\boundary\Del^1\subset\Del^1)^{\leib\times k}$.\footnote{For more details about the Leibniz or ``pushout-product'' construction see \cite[\S 4]{RiehlVerity:2013kx}.}
     
  \item\label{itm:inner-horn} \textbf{inner horns:} The inclusion $\gC{\Horn^{n,k}}\inc\gC\Del^n$ is identity on objects and full on all function complexes except for the one from $0$ to $n$. The inclusion
     \[
  \xymatrix@R=1em{ \Fun_{\gC{\Horn^{n,k}}}(0,n)\ar@{^(->}[r] \ar@{}[d]|{\rotatebox{90}{$\cong$}}& \Fun_{\gC\Del^n}(0,n)\ar@{}[d]|{\rotatebox{90}{$\cong$}} \\
\CHorn^{n-1,k}_1\ar@{^(->}[r] & \Cube^{n-1}}\]
 is isomorphic to the cubical horn inclusion, defined by the the following Leibniz product:
  \begin{equation*}
    (\boundary\Del^1\subset\Del^1)^{\leib\times(j-1)}\leib\times
    (\Del^{\fbv{1}}\subset\Del^1)\leib\times
    (\boundary\Del^1\subset\Del^1)^{\leib\times(k-j)}
  \end{equation*}
  \end{enumerate}
\end{ex}

\begin{defn}[bead shapes]
  We shall call those atomic arrows $T^\bullet\colon 0 \to n$ of $\gC{\Del^{n}}$
  which are not members of $\gC\boundary\Del^n$ \emph{bead shapes}. By Examples \ref{ex:homotopy-coherent-simplex} and \ref{ex:subcomputad-of-simplex}, an $r$-dimensional bead shape $T^\bullet\colon
  0\to n$ is given by a sequence of subsets
  \[ \{0,n\} = T^0 \subset T^1 \subset \cdots \subset T^r = [0,n]\]
  with $T^0 = \{0,n\}$  and $T^r$ equal to the full interval $[0,n] = \{ 0 \leq t \leq n\}$. 
  \end{defn}

More generally, any simplicial category $\gC{X}$ arising as the homotopy coherent realization of a simplicial set $X$ defines a simplicial computad whose atomic arrows $(x,T^\bullet)$, described in Proposition \ref{prop:gothic-C}, are called \emph{beads in $X$}.  As a consequence of this result we find that 
  $r$-simplices of $\gC{X}$ correspond to sequences of abutting beads,
  structures which are called \emph{necklaces} in the work of Dugger and
  Spivak~\cite{DuggerSpivak:2011ms} and Riehl~\cite{Riehl:2011ot}. In this terminology, 
   $\gC X$ is a simplicial computad in which the atomic arrows are those necklaces that consist of a single bead with non-degenerate image.

\begin{prop}[{$\gC X$ as a simplicial computad; \refVI{prop:gothic-C}}]\label{prop:gothic-C}
The homotopy coherent realization $\gC{X}$ of a simplicial set $X$ is a simplicial computad with
  \begin{itemize}
  \item   objects the vertices of $X$ and 
  \item non-degenerate atomic $r$-arrows given by pairs $(x,T^\bullet)$, wherein
  $x$ is a non-degenerate $n$-simplex of $X$ for some $n > r $ and
  $T^\bullet\colon 0\to n$ is an $r$-dimensional bead shape. 
  \end{itemize}
  The domain of $(x,T^\bullet)$ is the initial vertex $x_0$ of $x$ while the codomain is the 
terminal vertex $x_n$.
 \end{prop}
 
 The point of this review is to permit us to define weights for oplax colimits of homotopy coherent diagrams valued in a simplicially (or frequently quasi-categorically) enriched category. In a homotopy coherent diagram, the indexing shape is given by the homotopy coherent realization of a simplicial set $X$. In this context, the join operation $X\join\Del^0$ produces another simplicial set with a freely adjoined cocone vertex. We shall argue that the its homotopy coherent realization defines a collage that presents the weight for oplax colimits.
 
 \begin{rec} For any simplicial set $X$, there is a canonical inclusion $X\inc X\join\Del^0$ into its join
  with the point.  The join $X\join\Del^0$ has a single vertex of $X\join\Del^0$ that is not also a vertex of its subset $X$, which we shall denote by ``$\top$.'' Now for each non-degenerate $n$-simplex $x\in X$ the join $X\join\Del^0$ has
  two corresponding non-degenerate simplices:
  \begin{itemize}
  \item a simplex of dimension $n$ identified with $x$
itself and 
  \item a simplex $(x,\top)$ of dimension $n+1$,
  \end{itemize} and these two cases enumerate all of the
  non-degenerate simplices of $X\join\Del^0$ with the exception of $\top$.
\end{rec}

\begin{defn}[weights for oplax colimits]\label{defn:oplax-colimit-weight}
Applying homotopy coherent
  realisation to the canonical inclusion $X\inc X\join\Del^0$, we obtain a simplicial
  subcomputad $I_X\colon \gC{X} \inc \gC[X\join\Del^0]$ for any simplicial set $X$.

  Now from Proposition~\refVI{prop:gothic-C} we know that $\gC[X\join\Del^0]$ may be
  built from $\gC{X}$ by adjoining atomic arrows corresponding to beads
  $((x,\top),T^\bullet)$ and these all have codomain $\top$.
It is clear that the conditions discussed in
Proposition \ref{prop:collage-adjunction}\ref{itm:collage-fun} hold for the inclusion $\langle \top, I_X
  \rangle\colon \catone + \gC{X} \inc \gC[X\join\Del^0]$. Hence, via the counit isomorphism of the collage adjunction, this simplicial category is isomorphic to the collage of the corresponding weight
  \begin{equation*}
    \xymatrix@R=0em@C=5em{
      \gC{X}\op\ar[r]^{L_X} & {\SSet}
    }\qquad\text{given~by}\qquad L_X(x) \defeq \Fun_{\gC[X\join\Del^0]}(x,\top).
  \end{equation*}
  
  We refer to $L_X$ as the \emph{weight for oplax colimits} of diagrams of
  shape $\gC{X}$. When $F\colon\gC{X}\to\eK$ is a homotopy coherent diagram of
  shape $X$, then its \emph{oplax colimit} is defined to be the weighted colimit
  \[ \colim\nolimits^{\mathrm{oplax}}F \defeq \colim\nolimits^{L_X}F,\] should this exist in $\eK$. 
\end{defn}

\begin{rmk}\label{rmk:oplax-vs-pseudo}
The oplax weights being defined here are precisely the ``pseudo'' weights introduced in Definition \refVII{defn:weight-for-pseudo-limits}. The reason for the difference in nomenclature is that in that paper the diagrams considered in \cite{RiehlVerity:2018rq} 
 are valued in Kan complex enriched categories, whereas here the diagrams are valued in quasi-categorically (or simplicially) enriched categories. In a Kan complex, the 1-simplex $\Del^1$ represents an invertible morphism, while in a quasi-category it models a non-invertible morphism.
\end{rmk}

Immediately from Proposition
  \ref{prop:flexible-weight-as-computad}:

\begin{lem}[{\refVII{lem:pseudo-is-flexible}}]\label{lem:oplax-is-flexible}
  For all simplicial sets $X$ the weight $L_X\colon\gC{X}\op\to\SSet$ for oplax
  colimits of diagrams of shape $\gC{X}$ is a flexible weight.
\end{lem}


%% file: cartesian.tex

\section{Cocartesian fibrations and quasi-categorical collages}\label{sec:cocartesian}

In this section, we construct an explicit example of an oplax colimit of diagram of quasi-categories via the \emph{quasi-categorical collage construction}, which we introduce in \S\ref{sec:qcat-collage}. In an important special case, the quasi-categorical collage defines a cocartesian fibration over the 1-simplex, so we begin in \S\ref{sec:cocart-basics} with a review of this notion, since we shall need it anyways. In \S\ref{sec:cat-of-cocart}, we introduce the quasi-categorically enriched category of cocartesian fibrations and cartesian functors and observe that pullback defines a functor between such categories. 

\subsection{Cocartesian fibrations of quasi-categories}\label{sec:cocart-basics}

Of the many equivalent definitions of cocartesian fibration (see \S\refIV{sec:cartesian} and \S\refVI{sec:cocartesian}), the following will be the most convenient for this paper:

\begin{defn}[{\refIV{cor:lurie-cartesian}}]\label{defn:qcat-cocart} Let $p\colon \qE \tfib \qB$ be an isofibration between quasi-categories.
\begin{enumerate}[label=(\roman*)]
\item\label{itm:qcat-cocart-arrow} A 1-arrow $\chi\colon e\to e'$ of $\qE$ is $p$-\emph{cocartesian} if and only if any
  lifting problem
  \begin{equation*}
    \xymatrix@C=3em{
      \Delta^{\fbv{0,1}}\ar@/^2ex/[rr]!L+/u 4pt/^\chi \ar[r] &
      \Horn^{n,0} \ar[d] \ar[r] &      \qE \ar@{->>}[d]^{p} \\ &
      \Delta^n \ar[r] \ar@{-->}[ur] & \qB}
  \end{equation*}
  has a solution. 
  \item\label{itm:qcat-cocart-fib} An isofibration
  $p\colon\qE\tfib\qB$ is a \emph{cocartesian} fibration of quasi-categories precisely when any arrow
  $\alpha\colon pe \to b$ in $\qB$ admits a lift to an arrow $\chi\colon e\to
  e'$ in $\qE$ which enjoys the lifting property of \ref{itm:qcat-cocart-arrow}.
  \end{enumerate}
\end{defn}

For efficiency of exposition, we focus largely on the cocartesian fibrations, and leave it to the reader to formulate the dual statements for cartesian fibrations, obtained by replacing each simplicial set by its opposite.

\begin{ex} The product projection $\pi \colon \qA \times \qB \tfib \qB$ defines a \emph{bifibration}, that is, both a cocartesian and a cartesian fibration. A 1-arrow of $\qA \times \qB$ is $\pi$-cocartesian (and also $\pi$-cartesian) just when its component in $\qA$ is an isomorphism.
\end{ex}

\begin{ex}[{\refIV{ex:domain-fibration}, \refIV{ex:comma-fibrations}}]\label{ex:comma-qcat}  For any quasi-category $\qB$, we write $\qB^\cattwo \defeq \qB^{\Del^1}$ for its cotensor with the 1-simplex and $p_0,p_1 \colon \qB^\cattwo \tfib \qB$ for the evaluation maps at the vertices $0,1 \in \Del^1$ respectively. The codomain functor $p_1 \colon \qB^\cattwo \tfib \qB$ is a cocartesian fibration, in which the $p_1$-cocartesian arrows are those whose projections along $p_0$ are invertible. Dually, the domain functor $p_0 \colon \qB^\cattwo \tfib \qB$ is a cartesian fibration, in which the $p_0$-cartesian arrows are those whose codomain components are invertible.

More generally, for any cospan of quasi-categories $f \colon \qB \to \qA$ and $g \colon \qC \to \qA$, the \emph{comma quasi-category} is defined by the pullback
\[ 
\xymatrix{ f\comma g \ar@{->>}[d]_{(p_1,p_0)} \ar[r] \pbexcursion&  \qA^\cattwo \ar@{->>}[d]^{(p_1,p_0)} \\ \qC \times \qB\ar[r]_{g \times f} & \qA \times \qA}\]  and once more $p_1 \colon f \comma g \tfib \qC$ is a cocartesian fibration and $p_0 \colon f \comma g \tfib \qB$ is a cartesian fibration.

In the special case where one of the functors in the cospan is taken to be the identity, we write $f \comma \qA$ and $\qA \comma f$ for what we call the \emph{contravariant} and \emph{covariant} representable comma quasi-categories respectively. In the special case where the functor $a \colon 1 \to \qA$ identifies a vertex of $\qA$, the codomain projection $p_1 \colon a \comma \qA \tfib \qA$ is a cocartesian fibration that encodes the covariant representable functor associated to $a$, while the domain projection $p_0 \colon a \comma \qA \tfib \qA$ is a cartesian fibration encoding the contravariant representable functor. These define the images of $a$ in the co- and contravariant Yoneda embeddings of \S\refVI{sec:yoneda}.
\end{ex}

\begin{lem}[{\refVI{lem:cocart-cyl-are-cocart-arrows}}]\label{lem:ptwise-cocart-cyl} If $p \colon \qE \tfib \qB$ is a cocartesian fibration and $X$ is a simplicial set, then $p^X \colon \qE^X \tfib \qB^X$ is a cocartesian fibration in which a 1-arrow $e \colon \Del^1 \to \qE^X$ is $p^X$-cocartesian just when for each vertex $x \in X$ its component $e(x \cdot \degen^0,\id_{[1]}) \colon \Del^1 \to \qE$ is $p$-cocartesian.
\end{lem}

In \S\refVI{sec:cocartesian}, a $p^X$-cocartesian 1-arrow is called a \emph{pointwise $p$-cocartesian cylinder}.

  \begin{lem}[{\refVI{lem:cocart-cylinder-extensions}}]\label{lem:cocart-cylinder-extensions}
Let $X \inc Y$ be a simplicial subset of a simplicial set $Y$.
\begin{enumerate}[label=(\roman*)]
  \item Any lifting problem
    \begin{equation*}
      \xymatrix@=2em{
        {X\times\Del^1\cup Y\times\Del^{\{0\}}}\ar[r]^-{e}\ar@{^(->}[d] &
        {\qE}\ar@{->>}[d]^{p} \\
        {Y\times\Del^1}\ar[r]_{b}\ar@{.>}[ru]^{\bar{e}} & {\qB}
      }
    \end{equation*}
    in which the cylinder ${X\times\Del^1}\subseteq{X\times\Del^1\cup
      Y\times\Del^{\{0\}}}\stackrel{e}\longrightarrow {E}$ is pointwise
    $p$-cocartesian admits a solution $\bar{e}$ which is also pointwise
    $p$-cocartesian.
  \item Any lifting problem ($n>1$)
    \begin{equation*}
      \xymatrix@=2em{
        {X\times\Del^n\cup Y\times\Horn^{n,n}}\ar[r]^-{e}\ar@{^(->}[d] &
        {\qE}\ar@{->>}[d]^{p} \\
        {Y\times\Del^n}\ar[r]_{b}\ar@{.>}[ru]^{\bar{e}} & {\qB}
      }
    \end{equation*}
    in which the cylinder ${Y\times\Del^{\{n-1,n\}}}\subseteq {X\times\Del^n\cup
      Y\times\Horn^{n,n}}\stackrel{e}\longrightarrow {E}$ is pointwise
    $p$-cocartesian admits a solution $\bar{e}$.
\end{enumerate}
\end{lem}

\begin{defn} If $p$ and $q$ are cocartesian fibrations over $\qB$ then a functor
\[ \xymatrix{ \qE \ar@{->>}[dr]_p \ar[rr]^g & & \qF \ar@{->>}[dl]^{q} \\ & \qB}\] is a \emph{cartesian functor} just when it carries $p$-cocartesian 1-arrows to $q$-cocartesian 1-arrows.
\end{defn}

As one illustration of the importance of this notion, we have the following important proposition:

\begin{prop}[{\refVIII{prop:equivalence-of-fibrations}}]\label{prop:equivalence-of-fibrations} A cartesian functor
\[ \xymatrix{ \qE\ar[rr]^g \ar@{->>}[dr]_p & & \qF \ar@{->>}[dl]^q \\ & \qB}\] between cocartesian fibrations of quasi-categories is an equivalence if and only if it is a fiberwise equivalence: for each $b \in \ob\qB$, the induced functor $g_b \colon \qE_b \to \qF_b$ is an equivalence. 
\end{prop}

\subsection{The quasi-categorically enriched category of cocartesian fibrations}\label{sec:cat-of-cocart}

If $\qB$ is a quasi-category, then we adopt the notation $\qCat_{/\qB}$ for the quasi-categoricaly enriched category of isofibrations over $\qB$ defined as follows.

\begin{defn}\label{defn:fun-over-B} For a quasi-category $\qB$, let $\qCat_{/\qB}$  denote the category whose:
\begin{itemize}
\item objects are isofibrations $p \colon \qE \tfib \qB$ with codomain $\qB$ and
\item whose function complexes $\Fun_{\qB}(p \colon \qE \tfib \qB, q \colon \qF \tfib \qB)$ are defined by the pullbacks
\[ \xymatrix{ \Fun_{\qB}(p \colon \qE \tfib \qB, q \colon \qF \tfib \qB) \ar@{->>}[d] \ar[r] \pbexcursion & \Fun(\qE,\qF) \ar[d]^{q \circ -} \\ \Del^0 \ar[r]^-p & \Fun(\qE,\qB)}\]
\end{itemize}
where $\Fun(\qE,\qF) \cong \qF^\qE$ denotes the usual internal hom in $\qCat$.
\end{defn}

\begin{defn}\label{defn:cart-fun-over-B} For a quasi-category $\qB$, let $\coCart(\qCat)_{/\qB}$  denote the category whose:
\begin{itemize}
\item objects are cocartesian fibrations $p \colon \qE \tfib \qB$ with codomain $\qB$ and
\item whose function complexes $\Fun^c_{\qB}(p \colon \qE \tfib \qB, q \colon \qF \tfib \qB)$ are defined to be the full sub quasi-categories of the function complexes $\Fun_{\qB}(p \colon \qE \tfib \qB, q \colon \qF \tfib \qB)$ of $\qCat_{/\qB}$ defined by restricting the 0-arrows to be cartesian functors over $\qB$.
\end{itemize}
The quasi-categorically enriched category $\Cart(\qCat)_{/\qB}$ of cartesian fibrations and cartesian functors is defined similarly.
\end{defn}

Proposition \refIV{prop:cart-fib-pullback} proves that the pullback of a cocartesian fibration is a cocartesian fibration
\[
\xymatrix{ \qF \ar@{->>}[d]_q \ar[r]^g \pbexcursion & \qE\ar@{->>}[d]^p \\ \qA \ar[r]_f & \qB}\] in which an arrow $\chi$ is $q$-cocartesian if and only $g\chi$ is $p$-cocartesian. It follows that pullback also preserves cartesian functors. Hence:

\begin{prop}\label{prop:pullback-stability}
Pullback along any $f \colon \qA  \to \qB$ defines a quasi-categorically enriched functor
     \[
  \xymatrix@R=1em{ \coCart(\qCat)_{/\qB} \ar[r]^-{f^*} \ar@{}[d]|{\rotatebox{90}{$\supset$}}& \coCart(\qCat)_{/\qA}\ar@{}[d]|{\rotatebox{90}{$\supset$}} \\
\qCat_{/\qB} \ar[r]^-{f^*}& \qCat_{/\qA}}\]
\end{prop}

We now argue that the pullback functor preserves simplicial tensors. This will be used in \S\ref{sec:powerful} to show that its right adjoint is simplicially enriched, when this functor exists.

\begin{obs}[tensors and pullback]\label{obs:tensor-pullback} Let $X \in \sSet$ be a simplicial set. The tensor of an isofibration $p \colon \qE \tfib \qB$ with $X$ is the right-hand vertical composite, which pulls back to the right-hand vertical composite
\[ \xymatrix{ \qF \times X \ar@{->>}[d]_{\pi_1} \ar[r] \pbexcursion & \qE \times X \ar@{->>}[d]^{\pi_1}\\  \qF \pbexcursion \ar@{->>}[d]_{f^*(p)} \ar[r] & \qE \ar@{->>}[d]^p \\ \qA \ar[r]_f & \qB}\] which defines the tensor of $f^*(p) \colon \qF \tfib \qA$ with $X$.
\end{obs}

The following lemma tells us that this tensor construction respects cartesian functors.

\begin{lem}\label{lem:cart-functor-tensor} For any simplicial set $X$ and cocartesian fibrations $p \colon \qE \tfib \qB$ and $q \colon \qF \tfib \qB$, the isomorphism $\Fun_{\qB}(\qE \times X, \qF) \cong \Fun_{\qB}(\qE,\qF)^X$  restricts to an isomorphism \[\Fun^c_{\qB}(\qE \times X, \qF) \cong \Fun^c_{\qB}(\qE,\qF)^X.\]
\end{lem} 
\begin{proof}
We make use of Theorem \refIV{thm:cart.fun.chars} which provides the following characterization of the sub quasi-category $\Fun^c_{\qB}(\qE,\qF) \subset \Fun_{\qB}(\qE,\qF)$. Any functor $f \colon \qE \to \qF$ over $\qB$ induces a commutative square over $\qB$
\[ \xymatrix{ \qE \ar[r]^f \ar[d] & \qF \ar[d] \\ p \comma \qB \ar[r]_{(f,\id_\qB)} \ar@/^2ex/@{-->}[u]^\ell_\dashv & q \comma \qB \ar@/_2ex/@{-->}[u]_\ell^\vdash}\] whose vertical functors are the canonical ones induced by $p \colon \qE \tfib \qB$ and $q \colon \qF \tfib \qB$. Because $p$ and $q$ are cocartesian, Theorem \refIV{thm:cart.fib.chars} proves the vertical functors admit left adjoints over $\qB$. Theorem \refIV{thm:cart.fun.chars} proves that $f$ is cartesian if and only if the mate of this canonical isomorphism is an isomorphism.

The mate that detects whether $f$ is a cartesian functor lives as a 1-simplex in the simplicial set
\[  \qop{Sq}_{\qB}(p\comma \qB \to \qE,  q \comma \qB \to \qF) \defeq \Fun_{\qB}(\qE , \qF) \times_{\Fun_{\qB}(p \comma \qB, \qF)} \Fun_{\qB}(p \comma \qB, q\comma \qB).\] of commutative squares from $\ell \colon p\comma \qB \to \qE$ to $\ell \colon q \comma \qB \to \qF$. The adjunction over $\qB$ associated to the cocartesian fibration $\qE \times X \xtfib{\pi} \qE\xtfib{p} \qB$ is
\[ \xymatrix{ \qE \times X \ar[rr]^-\perp \ar[dr]_-{p\pi}  & & p\comma\qB \times X \ar[dl]^{p_0\pi} \ar@/_2ex/[ll]_-\ell \\ & \qB }\] the product of the adjunction for $p$ with $X$. In particular, 
\[ \qop{Sq}_{\qB}(p\comma \qB \times X \to \qE \times X, q \comma \qB \to \qF)  \cong \qop{Sq}_{\qB}(p\comma \qB \to \qE, q \comma \qB \to \qF)^X.\] Now a 1-simplex $\qC^X$ is an isomorphism if and only if it is a pointwise isomorphism, which proves that $\Fun^c_{\qB}(\qE \times X, \qF) \cong \Fun^c_{\qB}(\qE,\qF)^X$.
\end{proof}

\subsection{Collages for quasi-categories}\label{sec:qcat-collage}

We conclude this section with an example of an oplax colimit. When $X=\Del^1$ a homotopy coherent diagram $\gC{\Del^1} \to \qCat$ is just a functor $f \colon \qA \to \qB$ between quasi-categories. The oplax colimit in simplicial sets is given by the pushout 
\[ 
\xymatrix{ \qA \ar[r]^f \ar[d]_{\id\times\face^0} & \qB \ar[d] \\ \qA \times \Del^1 \ar[r] & \colim\nolimits^{\mathrm{oplax}}f \poexcursion }
\]
Our aim is to prove Proposition \ref{prop:oplax-colimit-of-functor} which demonstrates that  this oplax colimit is modeled up to equivalence by the \emph{quasi-categorical collage construction} that we now introduce.

\begin{defn}[the quasi-categorical collage construction]\label{defn:qcat-collage} Consider any cospan $f \colon \qA \to \qC$ and $g \colon \qB \to \qC$, with $\qA$, $\qB$, and $\qC$ all quasi-categories. Define a new simplicial set $\coll(f,g)$ 
by declaring that
\[ \coll(f,g)_n = \Bigl\{ \left(\Del^i \xrightarrow{a} \qA, \Del^j \xrightarrow{b} \qB, \Del^n \xrightarrow{c} \qC\right) \bigg\vert {\begin{array}{ll} c \vert_{\fbv{0,\ldots, i}} &= f(a), \\ c\vert_{\fbv{n-j,\ldots, n}} &= g(b),\end{array}} \begin{array}{c} i,j \geq -1, \\ i+j = n-1.\end{array}\Bigr\}\]
with the convention that conditions indexed by $\Del^{-1}$ are empty (or that each simplicial set is terminally augmented). There are simplicial maps
\[
\begin{tikzcd} B \arrow[r, hook] \arrow[d] \arrow[dr, phantom, "\lrcorner" very near start] & \coll(f,g) \arrow[d, "\rho"] & A \arrow[l, hook'] \arrow[d] \arrow[dl, phantom, "\llcorner" very near start] \\ \{1\} \arrow[r, hook] & \Del^1 & \{0\} \arrow[l, hook']
\end{tikzcd}
\]
 the top ones being the evident inclusions. The map $\rho$ sends an $n$-simplex $(a \colon \Del^i \to \qA, b \colon \Del^j \to \qB, c \colon \Del^n \to \qC)$ to the $n$-simplex $[n] \to [1]$ that carries $0,\ldots ,i$ to $0$ and $i+1,\ldots, n$ to $1$. Note that the fiber of $\rho$ over $0$ is isomorphic to $\qA$ while the fiber of $\rho$ over $1$ is isomorphic to $\qB$.
\end{defn}

\begin{rmk}[on right and left]
As with simplicial collages, we customarily write $\qB+\qA \inc \coll(f,g)$ for the inclusions of the fibers over $1$ and $0$ --- with the fiber over 1 on the left and the fiber over 0 on the right. As with our convention for quasi-categories in Example \ref{ex:comma-qcat}, this positions the covariantly-acting quasi-category on the ``left'' and the contravariantly-acting quasi-category on the ``right.''
\end{rmk}

\begin{lem}\label{lem:collage-qcat} The map $\rho \colon \coll(f,g) \to \Del^1$ is an inner fibration. In particular, the simplicial set $\coll(f,g)$ is a quasi-category.
\end{lem}
\begin{proof} Since the fibers of $\rho$ over $0$ and $1$ are the quasi-categories $\qA$ and $\qB$, it suffices to consider inner horns
\[
\xymatrix{ \Horn^{n,k} \ar[r] \ar[d] & \coll(f,g) \ar[d]^\rho \\ \Del^n \ar[r]_{\alpha} \ar@{-->}[ur] & \Del^1}\] for which $\alpha \colon [n] \to [1]$ is a surjection. Suppose $\alpha$ carries $0,\ldots, i$ to $0$ and $i+1,\ldots, n$ to $1$. Note that for any $0 < k < n$, the faces $\fbv{0,\ldots, i}$ and $\fbv{i+1,\ldots,n}$ of $\Del^n$ belong to the horn $\Horn^{n,k}$. In particular, the map $\Horn^{n,k} \to \coll(f,g)$ identifies simplices $a \colon \Del^i \to \qA$ and $\Del^{n-i-1} \to \qB$ together with a horn $\Horn^{n,k} \to \qC$ whose initial and final faces are the images of these simplices under $f \colon \qA \to \qC$ and $g \colon \qB \to \qC$. Since $\qC$ is a quasi-category this horn admits a filler $c \colon \Del^n \to \qC$ and the triple $(a,b,c)$ defines an $n$-simplex in $\coll(f,g)$ that solves the lifting problem.
\end{proof}

We write $\coll(f,\qB)$ for the collage of $f \colon \qA \to \qB$ with the identity on $\qB$.

\begin{lem}\label{lem:collage-qcat-cocartesian} For any $f \colon \qA \to \qB$, the map $\rho \colon \coll(f,\qB) \to \Del^1$ is a cocartesian fibration. 
\end{lem}
\begin{proof} To prove the claim, we need only specify cocartesian lifts of the non-degenerate 1-simplex of $\Del^1$ and demonstrate that these edges have the corresponding universal property. To that end, for any vertex $a \in \qA_0$, let $\chi_a \colon \Del^1 \to \coll(f,\qB)$ be the 1-simplex \[ \chi_a := (a \colon \Del^0 \to \qA, fa \colon \Del^0 \to \qB, fa \cdot \degen^0 \colon \Del^1 \to \qB),\] defined by the degenerate edge at $fa \in \qB_0$ lying over the 1-simplex in $\Del^1$. To show that $\chi_a$ is $\rho$-cocartesian, we must construct fillers for any left horn
\[
\xymatrix{ \Del^{\fbv{0,1}} \ar@/^3ex/[rr]^{\chi_a} \ar[r] & \Horn^{n,0} \ar[d] \ar[r] & \coll(f,\qB) \ar[d]^\rho \\ & \Del^n \ar@{-->}[ur] \ar[r]_\beta & \Del^1}\] whose initial edge is $\chi_a$. Note that this condition implies that the bottom map $\beta \colon [n] \to [1]$ carries $0$ to $0$ and the remaining vertices to $1$. The map $\Horn^{n,0} \to \coll(f,\qB)$ defines a horn $\Horn^{n,0} \to \qB$ in the quasi-category $\qB$ whose first edge is degenerate. By Joyal's lemma about filling ``special outer horns,'' such horns admit a filler $b \colon \Del^n \to \qB$ and the triple
\[ (a \colon \Del^0 \to \qA, b \cdot \face^0 \colon \Del^{n-1} \to \qB, b \colon \Del^n \to \qB)\] defines an $n$-simplex in $\coll(f,\qB)$ that solves the lifting problem.
\end{proof}

\begin{prop}\label{prop:oplax-colimit-of-functor}
For any $f \colon \qA \to \qB$ between quasi-categories, the collage $\coll(f,\qB)$ defines the oplax colimit of $f$ in $\qCat$. That is $\coll(f,\qB)$ defines a cone under the pushout diagram
\[ \xymatrix{ \qA \ar[r]^f \ar@{^(->}[d]_{\id\times\face^0} & \qB \ar@{^(->}[d]  \ar@{^(->}@/^/[ddr] \\ \qA \times \Del^1 \ar[r] \ar@/_/[drr]_h & P \poexcursion \ar@{-->}[dr]^k  \\ & & \coll(f, {\qB})}\] so that the induced map $k$ is inner anodyne, and in particular a weak equivalence in the Joyal model structure.
\end{prop}
\begin{proof} The map $k$ is a quotient of the map $h$, which has the following explicit description. For each 
 $n$-simplex $(a,\alpha) \colon \Del^n \to \qA \times \Del^1$ define $i \defeq |\alpha^{-1}(0)| -1$, so that $-1 \leq i \leq n$. Then $h$ carries $(a,\alpha)$ to the $n$-simplex of $\coll(f,\qB)$ corresponding to the triple
\[ (a\vert_{\fbv{0,\ldots, i}} \colon \Del^i \to \qA, fa\vert_{\fbv{i+1,\ldots,n}} \colon \Del^{n-i-1} \to \qB, fa \colon \Del^n \to \qB).\] Note that the composite $\rho h \colon \qA \times \Del^1 \to \Del^1$ is the projection.

It remains to present $k$ as a sequential composite of pushouts of coproducts of inner horn inclusions. To do so, first note that
\[ \coll(f,\qB)_n = \qA_n \coprod \qA_{n-1} \times_{\qB_{n-1}} \qB_n \coprod \cdots \coprod \qA_0 \times_{\qB_0} \qB_n \coprod \qB_n\] 
where each map $\qB_n \to \qB_i$ is the initial face map corresponding to $\Del^{\fbv{0,\ldots, i}} \inc\Del^n$. From the perspective of this decomposition, $P_n$ is the subset containing the sets $\qA_n$ and $\qB_n$ and the subset of $\qA_i \times_{\qB_i} \qB_n$ whose component in $\qB_n$ is in the image of $f$. The $n$-simplices of $\coll(f,\qB)$ that remain to be attached correspond to elements of $\qA_i \times_{\qB_i} \qB_n$, for $0 \leq i < n$, that are not in the image of $f$ in the sense just discussed. Note in particular that $k\colon P_0 \inc \coll(f,\qB)_0$ is an isomorphism and $k \colon P_n \inc \coll(f,\qB)_n$ is an injection for all $n \geq 1$.

To enumerate our attaching maps, we start with the collection of non-degenerate $n$-simplices of $\coll(f,\qB)$ for $n \geq 1$ that are not in the image of $f$ and remove also those elements of $\qA_i \times_{\qB_i} \qB_n$ whose components $b \in \qB_n$ are in the image of the degeneracy map $\degen_i \colon \qB_{n-1} \to \qB_n$. Partially order this set of simplices first in the order of increasing $n$ and the in order of increasing index $i$; that is we lexicographically order the collection of pairs $(n,i)$ with $n \geq 1$ and $0 \leq i < n$. We will filter the inclusion $P\inc\coll(f,\qB)$ as
\[ P \inc P_{(1,0)} \inc P_{(2,0)} \inc P_{(2,1)} \inc P_{(3,0)} \inc \cdots \inc P_{(n,i)} \inc \cdots \inc \colim\cong \coll(f,\qB)\]
where the simplicial set $P_{(n,i)}$ is built from the previous one by a pushout of a coproduct of inner horns indexed by the set of $n$-simplices $(a,b) \in \qA_i \times _{\qB_i} \qB_n$ with $b$ not in the image of $f$ or $\degen_i$. The filler for the horn indexed by $(a,b)$ will attach this $n$ simplex to $\qB_n$ as the missing face of the horn and also the $n+1$ simplex $(a,b\sigma^i) \in \qA_i \times_{\qB_i} \qB_{n+1}$.

Consider a simplex $(a,b) \in \qA_i \times_{\qB_i} \qB_n$ with $b$ not in the image of $f$ or $\degen_i$. Define a horn
\[
\xymatrix{ \Horn^{n+1,i+1} \ar[r] \ar@{^(->}[d] & P_{(n,i)} \ar@{^(->}[d] \\ \Del^{n+1} \ar[r]_-{(a,b\degen^i)} & \coll(f,\qB)}\] For each $0 \leq j < i+1$, the $\face^j$-face of the $n+1$ simplex $(a,b\sigma^i)$ is the $n$-simplex $(a\face^j,b\sigma^i\face^j)$, which lies in $P_{(n,i-1)}$ or in $\qB \inc P$ in the case $i=0$. For each $i+1 < j \leq n+1$, the $\face^j$-face of the $n+1$ simplex $(a,b\sigma^i)$ is the $n$-simplex $(a, b\sigma^i\face^j) = (a,b\face^{j-1}\degen^{i}) \in \qA_i \times_{\qB_i} \qB_n$, which was previously attached to $P_{(n-1,i)}$. So the $\Horn^{n+1,i+1}$ indeed maps to $P_{(n,i)}$, permitting an inductive construction of the next simplicial set in this sequence as the pushout
\[ \xymatrix{\coprod\limits_{\sim} \Horn^{n+1,i+1} \ar[r] \ar@{^(->}[d] & P_{(n,i)} \ar[d] \\ \coprod\limits_{\sim} \Del^{n+1} \ar[r] & P_{(n,i)+1}\poexcursion}
\] where $ P_{(n,i)+1}$ equals  $P_{(n+1,0)}$ in the case $i=n-1$ and $P_{(n,i+1)}$ otherwise.
\end{proof}

\begin{cor}\label{cor:lurie-adjunction} Consider a pair of functors between quasi-categories $f \colon \qA \to \qB$ and $u \colon \qB \to \qA$. Then $f$ is left adjoint to $u$ if and only if the collages $\coll(f,\qB)$ and $\coll(\qA, u)$ are equivalent under $\qB + \qA$ and over $\Del^1$. 
\end{cor}
\begin{proof}
First suppose that $\coll(f,\qB) \simeq \coll(\qA, u)$  under $\qB + \qA$ and over $\Del^1$. By Lemma \ref{lem:collage-qcat-cocartesian} this means that the map $\coll(f,\qB) \to \Del^1$ is both a cocartesian and a cartesian fibration, a \emph{bifibration} in the terminology of \S\refIV{sec:cartesian}. By Proposition \refIV{prop:bifibration-adjunction} it follows that the 1-arrow in $\Del^1$ from $0$ to $1$ induces an adjunction between the fibers $\qA$ and $\qB$. By inspection of that proof, the left adjoint functor so-constructed in the case of the bifibration $\coll(f,\qB) \to \Del^1$ is $f$; substituting the equivalent bifibration $\coll(\qA,u) \to \Del^1$, we see that the right adjoint is equivalent to $u$.

For the converse, we work in the opposite $\infty$-cosmos $\qCat\op$, an $\infty$-cosmos in which ``not all objects are cofibrant,'' as described in Observation \refIV{obs:duals-of-cosmoi}. In that context, Proposition \ref{prop:oplax-colimit-of-functor} proves that $\coll(f,\qB)$ and $\coll(\qA,u)$ construct the contravariant and covariant comma objects associated to the functors $f$ and $u$. If $f \dashv u$ in $\qCat$ then these functors are also adjoint in $\qCat\op$ and Proposition \refI{prop:adjointequiv} then proves that the commas $\coll(f,\qB)$ and $\coll(\qA,u)$ are equivalent under $\qB + \qA$. By construction, this equivalence also lies over $\Del^1$.
\end{proof}


%% file: comprehension-review.tex

\section{The comprehension construction}\label{sec:comprehension}

In this section we review the \emph{comprehension construction} from \cite{RiehlVerity:2017cc} in considerably less generality than given in that source. It constructs, for any cocartesian fibration $p \colon \qE \tfib\qB$ of quasi-categories, a ``straightening,'' which has the form of a simplicial functor $c_p \colon \gC\qB \to\qCat$ that sends each vertex $b \in \qB$ to the fiber $\qE_b$. It also constructs a canonical lax cocone $\ell^\qE \colon \gC[\qB\join\Del^0] \to \qCat$ of shape $\qB$ under this diagram with nadir $\qE$, the lax colimits of restrictions of which will be used in \S\ref{sec:pullback} to model pullbacks along the functor $p$.

The underlying mechanics of the comprehension construction are reviewed in \S\ref{sec:cocart-cocones} and the comprehension construction itself is given in \S\ref{sec:comprehension-review}.

\subsection{Cocartesian transformations between lax cocones}\label{sec:cocart-cocones}

Corollary \ref{cor:collage-bijection} tells us that the collage of a weight $W$ realizes the shape of $W$-cocones. Applying this result to the weights for oplax colimits introduced in Definition \ref{defn:oplax-colimit-weight}, we obtain the following definition of a \emph{lax cocone}.

\begin{defn}[{lax cocones \refVI{defn:lax-cocone}}]
  Suppose that $X$ is a simplicial set. Then a \emph{lax cocone of shape $X$\/} in $\SSet$ is defined
  to be a simplicial functor $\ell^B\colon\gC[X\join\Del^0]\to\SSet$
    \begin{equation*}
    \xymatrix@=1.5em{
      & {\catone+\gC{X}}\ar@{^(->}[dl]\ar[dr]^{\langle B,B_\bullet \rangle} & \\
      {\gC[X\join\Del^0]}\ar[rr]_-{\ell^B} && {\SSet}
    }
  \end{equation*}
  The restriction of a lax cocone
  $\ell^B\colon\gC[X\join\Del^0]\to\eK$ to a functor $B_{\bullet}\colon\gC{X}\to \SSet$ is
  called its \emph{base\/}. We say
  that $\ell^B$ is a lax cocone \emph{under the diagram\/} $B_\bullet$; the object
  $B\in\sSet$ obtained by evaluating $\ell^B$ at the object $\top$ is called the
  \emph{nadir} of that lax cocone.
  \end{defn}
  
  \begin{rmk}
  In the original Definition \refVI{defn:lax-cocone}, the target was required to be a quasi-categorically enriched category $\eK$ and the base of a lax cocone was required to factor through
  through the inclusion $g_*\eK\subseteq \eK$ of the maximal Kan complex enriched subcategory. The point of this requirement was so that the transpose of the base diagram defined a diagram $X \to \nrvhc g_*\eK$ valued in the large quasi-category of objects and morphisms in $\eK$. But in this paper we will frequently  consider lax cocones whose codomain is $\SSet$, in which case the requirement that the base diagram factor through the local groupoidal core won't necessarily make sense and should be dropped. We shall also be more interested in the totality of the data of a lax cocone than in studying the base diagram in isolation.
  \end{rmk}
  
  \begin{defn}[{canonical lax cocones \refVI{lem:constant-lax-cones}}]\label{defn:canonical-cocone} For any simplicial set $X$, there exists a lax cocone 
      \begin{equation*}
    \xymatrix@=1.5em{
      & {\catone+\gC{X}}\ar@{^(->}[dl]\ar[dr]^{\langle X,1 \rangle} & \\
      {\gC[X\join\Del^0]}\ar[rr]_-{k^X} && {\SSet}
    }
  \end{equation*}
  whose base is constant at the terminal quasi-category $1$ and whose nadir is $X$ that we refer to as the \emph{canonical $X$-shaped lax cocone}. 
  
To define $k^X \colon \gC[X \join\Del^0] \to \SSet$, it remains to define the images in $\Fun(1,X)\cong X$ of the atomic $r$-arrows with domain a vertex $x_0$ in $X$ and with codomain $\top$.   Applying Proposition \ref{prop:gothic-C}, each atomic $r$-arrow from $x_0$ to $\top$ corresponds to a bead $((x, \top), T^\bullet)$ represented by some non-degenerate $n$-simplex $x \in X$ whose initial vertex is $x_0$ and atomic $r$-arrow $T^\bullet \colon 0 \to n+1$ in $\gC\Del^{n+1}$ defined by
  \[ \{0,n+1\} = T^0 \subset T^1 \subset \cdots \subset T^r = [0,n+1].\]
  Define $\mu \colon \Del^r \to \Del^n$ by $\mu(i) = \max (T^i \backslash \{n+1\})$. Then the image of the bead $((x, \top), T^\bullet)$ is the $r$-simplex $x \cdot \mu \in X$ in the hom-quasi-category $X \cong \Fun(1,X)$ from the image of the domain vertex $x_0$ of $x$.
\end{defn}

\begin{obs}[{whiskering lax cocones \refVI{obs:whiskering-cocone}}]\label{obs:whiskering-cocone}
Let $\ell^A \colon \gC[X \join\Del^0] \to \SSet$ be a lax cocone with base diagram $A_\bullet \colon \gC{X} \to \SSet$ and nadir $\ell^A_\top = A$, and let $f \colon A \to B$ be any map of simplicial sets. Then there is a \emph{whiskered lax cocone} $f \cdot \ell^A \colon \gC[X \join\Del^0] \to \SSet$ with the same base diagram $A_\bullet \colon \gC{X} \to \SSet$ and with nadir $B$, whose components from a vertex $x \in X$ to $\top$ are defined by whiskering with $f$:
\[ \Fun_{\gC[X\join\Del^0]}(x,\top) \xrightarrow{ \ell^X_{x,\top}} \Fun(A_x,A) \xrightarrow{f \circ -} \Fun(A_x,B)\]
\end{obs}

\begin{lem}\label{lem:canonical-whiskering} 
For any map of simplicial sets $f \colon Y \to X$, the canonical lax cocone of shape $X$ restricts along $\gC[f \join\id] \colon \gC[Y\join\Del^0] \to \gC[X\join\Del^0]$ to the whiskered composite 
 \[
\begin{tikzcd}[row sep=small, column sep=tiny] & \catone+ \gC{Y} \arrow[dl, hook'] \arrow[rr, "\catone+\gC{f}"] & & \catone+ \gC{X}  \arrow[dl, hook'] \arrow[dr, "{\langle X,1 \rangle}"] & \arrow[dr, phantom, "="] & & \catone+\gC{Y} \arrow[dl, hook'] \arrow[dr, "{\langle X,1\rangle}"] \\ \gC[Y\join\Del^0] \arrow[rr, "{\gC[f \join\id]}"'] & & \gC[X\join\Del^0]\arrow[rr, "k^X"'] && \SSet &  \gC[Y\join\Del^0] \arrow[rr, "f \cdot k^Y"']  && \SSet
\end{tikzcd}
\]
of the canonical lax cocone of shape $Y$ with $f \colon Y \to X$.
\end{lem}
\begin{proof}
By direct verification from Definition \ref{defn:canonical-cocone} and Observation \ref{obs:whiskering-cocone}.
\end{proof}

 \S\refVI{sec:cocones} introduces a mechanism for producing new lax cocones from given ones: namely as domain components of cocartesian cocones over a given codomain lax cocone. 

\begin{defn}[cocartesian cocones {\refVI{defn:cocart-cocone}}]
    Suppose we are given a simplicial set $X$ and lax cocones $\ell^E,\ell^B
  \colon\gC[X\join\Del^0]\to \SSet$ of shape $X$
  with bases $E_\bullet$ and $B_\bullet$ respectively. Suppose also that we are given a
  simplicial natural transformation
  \begin{equation*}
    \xymatrix@R=0em@C=8em{
      {\gC[X\join\Del^0]}\ar@/_2ex/[]!R(0.5);[r]_{\ell^B}\ar@/^2ex/[]!R(0.5);[r]^{\ell^E}
      \ar@{}[r]|-{\textstyle \Downarrow p} & {\SSet.}
    }
  \end{equation*}
  Then we say that the triple $(\ell^E,\ell^B,p)$ is a \emph{cocartesian cocone} if
  and only if
  \begin{enumerate}[label=(\roman*)]
  \item\label{itm:cocart-cocone-i} the nadir of the natural transformation $p$, that being its component
    $p\colon \qE\tfib \qB$ at the object $\top$, is a cocartesian fibration between quasi-categories
  \item\label{itm:cocart-cocone-ii} for all $0$-simplices $x\in X$ the naturality square
    \begin{equation*}
      \xymatrix@=2em{
        {E_x}\ar@{->>}[d]_{p_x}\ar[r]^{\ell^E_{{x}}} \pbexcursion & {\qE}\ar@{->>}[d]^{p} \\
        {B_x}\ar[r]_{\ell^B_{{x}}} & {\qB}
      }
    \end{equation*}
    is a pullback, and
  \item\label{itm:cocart-cocone-iii} for all non-degenerate $1$-simplices $f\colon x\to y\in X$ the $1$-arrow
    \begin{equation*}
      \xymatrix@C=5em@R=1.5em{
        {E_x}\ar[dd]_{E_{f}}\ar[dr]^{\ell^E_{{x}}} & \\
        {}\ar@{}[r]|(0.35){\textstyle\Downarrow \ell^E_{{f}}} & {\qE} \\
        {E_y}\ar[ur]_{\ell^E_{{y}}} & 
      }
    \end{equation*}
    is $p$-cocartesian.
  \end{enumerate}
  In this situation we also say that the pair $(\ell^E,p)$ defines a \emph{cocartesian
  cocone over} $\ell^B$.
\end{defn}

\begin{lem}[{pullbacks of cocartesian cocones \refVI{lem:cocart-cocone-pb}}]\label{lem:cocart-cocone-pb} Suppose given:
\begin{itemize}
\item a pullback diagram of quasi-categories
\begin{equation}\label{eq:cocart-cocone-pb} \xymatrix{ \qF \ar@{->>}[d]_q \ar[r]^g \pbexcursion & \qE \ar@{->>}[d]^p \\ \qA \ar[r]_f & \qB}\end{equation} in which $p$ and $q$ are cocartesian fibrations;
\item a lax cocone $\ell^A \colon \gC[X\join\Del^0] \to \eK$  with nadir $\qA$; and
\item a cocartesian cocone $(\ell^E,\ell^B,p)$ whose nadir component is $p \colon \qE \tfib \qB$ and whose codomain cocone $\ell^B = f \cdot \ell^A$ is obtained from the lax cocone $\ell^A$ by whiskering with $f\colon \qA \to \qB$ as in Observation \ref{obs:whiskering-cocone}.
\end{itemize}
Then there is a cocartesian cone $(\ell^F,\ell^A,q)$ whose codomain is $\ell^A$, whose nadir component is $q \colon \qF \tfib \qA$, and whose domain component is a lax cocone $\ell^F$ that whiskers with $g$ to the lax cone $\ell^E = g \cdot \ell^F$.
\end{lem}

Conversely, a cocartesian cocone $(\ell^F,\ell^A,q)$ with nadir component $q \colon \qF \tfib \qA$ can be whiskered with a pullback square \eqref{eq:cocart-cocone-pb} to define a cocartesian cocone $(g\cdot \ell^F, f\cdot \ell^A, p)$ with nadir $p \colon \qE \tfib \qB$ and whose domain and codomain are whiskered lax cocones as defined in Observation \ref{obs:whiskering-cocone}.

\begin{rmk}\label{rmk:arbitrary-pullbacks}
If the map $f$ of Lemma \ref{lem:cocart-cocone-pb} is replaced by any map  of simplicial sets $f \colon X \to \qB$, whose domain is not necessarily a quasi-category, it is still possible to pull back the data of a cocartesian cocone $(\ell^E,\ell^B,p)$ whose codomain lax cocone $\ell^B = f \cdot \ell^X$ is obtained by whiskering a lax cocone with nadir $X$. This constructs a simplicial natural transformation  $(\ell^F,\ell^X,q)$ whose nadir component is the pullback $q \colon F \to X$ of $p$ along $f$. Since this is not a map between quasi-categories, it doesn't really make sense to call it a cocartesian fibration. Nonetheless, this construction produces a lax cocone $\ell^F$ of shape $X$, which will have some utility. See Remark \ref{rmk:generalized-unstraightening}.
\end{rmk}

\subsection{The comprehension construction}\label{sec:comprehension-review}

The comprehension functor associated to a cocartesian fibration $p \colon \qE \tfib \qB$ between quasi-categories is the base diagram of the domain of a cartesian cocone over the canonical $\qB$-shaped lax cocone. 

\begin{thm}[{\refVI{defn:basic-comprehension}}]\label{thm:comprehension}
For any cocartesian fibration $p\colon \qE\tfib \qB$ of quasi-categories, there is a cocartesian cocone 
  \begin{equation*}
    \xymatrix@R=0em@C=8em{
      {\gC[\qB\join\Del^0]}\ar@/_2ex/[]!R(0.5);[r]_{k^\qB}\ar@/^2ex/[]!R(0.5);[r]^{\ell^\qE}
      \ar@{}[r]|-{\textstyle \Downarrow p} & {\qCat.}
    }
  \end{equation*}
  of shape $\qB$ in $\qCat$ with nadir component $p \colon \qE\to\qB$ over the canonical lax cocone $k^\qB$ of Definition \ref{defn:canonical-cocone} whose base is constant at $1$. The base of the domain component defines the comprehension functor $c_p$, which acts on an object $b \colon 1 \to \qB$ of $\gC\qB$ by forming the pullback 
  \[
      \xymatrix@=1.7em{
      {\qE_b}\pbexcursion\ar[r]^{\ell^\qB_{{b}}}\ar@{->>}[d]_-{p_b} &
      {\qE}\ar@{->>}[d]^-p \\
      {1}\ar[r]_b & {\qB}
    }
    \] and acts on 1-arrows $f \colon a \to b$ of $\qB$ 
 by factoring the codomain of a $p$-cocartesian lift $\ell^\qE_{{f}}$ of $f$ through the pullback at the front of the diagram:
  \begin{equation}\label{eq:comprehension-on-1-arrows}
        \begin{xy}
      0;<1.4cm,0cm>:<0cm,0.75cm>::
      *{\xybox{
          \POS(1,0)*+{1}="one"
          \POS(0,1)*+{1}="two"
          \POS(3,0.5)*+{\qB}="three"
          \ar@{=} "one";"two"
          \ar@/_5pt/ "one";"three"_{b}^(0.1){}="otm"
          \ar@/^10pt/ "two";"three"^{a}_(0.5){}="ttm"|(0.325){\hole}
          \ar@{=>} "ttm"-<0pt,7pt> ; "otm"+<0pt,10pt> ^(0.3){f}
          \POS(1,2.5)*+{\qE_{b}}="one'"
          \POS(1,2.5)*{\pbcorner}
          \POS(0,3.5)*+{\qE_{a}}="two'"
          \POS(0,3.6)*{\pbcorner}
          \POS(3,3)*+{\qE}="three'"
          \ar@/_5pt/ "one'";"three'"_{\ell^\qE_{{b}}}
          \ar@/^10pt/ "two'";"three'"^{\ell^\qE_{{a}}}_(0.55){}="ttm'"
          \ar@{->>} "one'";"one"_(0.325){p_{b}}
          \ar@{->>} "two'";"two"_{p_{a}}
          \ar@{->>} "three'";"three"^{p}
          \ar@{..>} "two'";"one'"_*!/^2pt/{\scriptstyle \qE_f}
          \ar@{..>}@/_10pt/ "two'";"three'"^(0.44){}="otm'"
          \ar@{=>} "ttm'"-<0pt,4pt> ; "otm'"+<0pt,4pt> ^(0.3){\ell^\qE_{{f}}}
        }}
    \end{xy}
  \end{equation}
These cocartesian lifts define 
  components of the lax cocone
        \begin{equation*}
    \xymatrix@=1.5em{
      & {\catone+\gC{\qB}}\ar@{^(->}[dl]\ar[dr]^{\langle \qE,c_p \rangle} & \\
      {\gC[\qB\join\Del^0]}\ar[rr]_-{\ell^\qE} && {\qCat}
    }
  \end{equation*}
  with nadir $\qE$ under the comprehension functor.
  \end{thm}
  
  \begin{rmk}\label{rmk:straightening} The comprehension functor $c_p \colon \gC\qB \to \qCat$ is a ``straightening'' of the cocartesian fibration $p \colon \qE \tfib \qB$, a homotopy coherent diagram of shape $\qB$ that sends each vertex $b$ to the fiber $\qE_b$ of $p$ over $b$.
  \end{rmk}

By Observation \refVI{obs:unique-comprehension}, any pair of lax cocones that arise as the domain of a cocartesian cocone over a common codomain define vertices in a contractible Kan complex and are in particular equivalent as diagrams. This is used to prove the following result relating comprehension functors with pullbacks.

\begin{prop}[{comprehension and change of base \refVI{prop:comprehension-cob}}]\label{prop:comprehension-cob}
  Suppose that we are given a pullback
  \begin{equation*}
    \xymatrix@=2em{
      {\qF}\pbexcursion\ar@{->>}[d]_{q}\ar[r]^g &
      {\qE}\ar@{->>}[d]^{p} \\
      {\qA}\ar[r]_{f} & {\qB}
    }
  \end{equation*}
  of quasi-categories in which $p$ and thus $q$ are cocartesian
  fibrations. Then the diagrams\[
  \gC\qA \xrightarrow{c_q} \qCat\qquad \mathrm{and} \qquad \gC\qA \xrightarrow{\gC{f}} \gC\qB \xrightarrow{c_p} \qCat\] are connected by a homotopy coherent natural isomorphism.
\end{prop}

\begin{rmk}\label{rmk:generalized-unstraightening}
If $\qA$ is not a quasi-category, it is not possible to directly construct the comprehension functor for the pullback of $p$ along $f$. However, by  Lemmas \ref{lem:canonical-whiskering} and Remark \ref{rmk:arbitrary-pullbacks},  the cocartesian cocone over the canonical $\qB$-shaped lax cocone can be pulled back along any map of simplicial sets $f \colon X \to \qB$ to define a cocartesian cocone over the canonical $X$-shaped lax cocone. Thus, a posteriori, we can think of the base of the lax cocone
\[ \gC[X \join\Del^0] \xrightarrow{\gC[f\join\id]} \gC[\qB\join\Del^0] \xrightarrow{\ell^\qE} \qCat\] as defining a comprehension functor for the pullback of $p \colon \qE \tfib \qB$ along $f \colon X \to \qB$. 
\end{rmk}


%% file: pullback.tex

\section{Pullback along a cocartesian fibration as an oplax colimit}\label{sec:pullback}

Our aim in this section is to provide an equivalent model of the pullback functor 
\[ p^* \colon \sSet_{/\qB} \longrightarrow \sSet_{/\qE}\] along a cocartesian fibration $p \colon \qE \tfib \qB$ between quasi-categories. In the next section, we will use this to construct an equivalent model of its right adjoint, the pushforward $\Pi_p$, whose homotopical properties are more easily established. Before commencing with our work, we briefly sketch the connection between the pullback and pushforward.

The category of simplicial sets, as a presheaf topos, is locally
  cartesian closed, so pullback along any map $p \colon E \to B$ admits a right
  adjoint:
  \[
\begin{tikzcd} \sSet_{/E} \arrow[r, bend right, "{\Pi}_p"'] \arrow[r, phantom, "\bot"] & \sSet_{/B} \arrow[l, bend right, "{p}^*"'] 
\end{tikzcd}
\]
By Observation \ref{obs:tensor-pullback}, pullback along $p$ is a simplicially enriched functor that preserves tensors with
  simplicial sets, so by \cite[4.85]{kelly:ect} it follows that the adjunction $p^* \dashv \Pi_p$ is simplicially enriched.

The right adjoint may be described explicitly:

\begin{lem}\label{lem:strict-test-objects} The $n$-simplices of the pushforward $\Pi_pq \colon \Pi_pF \to B$ of $q \colon F \to E$ correspond to pair comprised of an $n$-simplex $b \colon \Del^n \to B$ together with a map $E_b \to F$ in $\sSet_{/E}$, whose domain is defined by the pullback
\[
\xymatrix@=2em{
        {E_b}\pbexcursion\ar[r]^{e_b}\ar[d]_{p_b} & {E}\ar[d]^{p} \\
        {\Del^n}\ar[r]_{b} & {B}}
\] Moreover, \ifnewcut\else under the representation
  \begin{equation*}
    \Pi_pF_n \cong \coprod_{b\colon\Del^n\to B} \Hom_E\left( 
      \vcenter{\xymatrix@1@R=2.5ex{E_b\ar[d]^{e_b} \\ E}},
      \vcenter{\xymatrix@1@R=2.5ex{F\ar[d]^{q} \\ E }}
    \right)
  \end{equation*}\fi 
a simplicial operator $\alpha\colon [m]\to
  [n]$ acts on an $n$-simplex by pre-composition with
  \[ \xymatrix{ E_{b \cdot \alpha} \pbexcursion \ar[r]^-{E_\alpha} \ar[d]_{p_{b\cdot\alpha}} & E_b \ar[d]^{p_b} \\ \Del^m \ar[r]_\alpha & \Del^n}\ifnewadd \qed\fi\] 
\end{lem}
\ifnewcut\else
\begin{proof}
  A generalised element $g\colon X\to\Pi_pF$  gives a map
  \begin{equation*}
    \xymatrix@=2em{
      X\ar[dr]_{b}\ar[rr]^{g} & & {\Pi_pF}\ar[dl]^{\Pi_pq} \\ & B & }
  \end{equation*}
  in $\sSet_{/B}$ in which we define $b\defeq g\bar{q}$.  Transposing this under
  the adjunction $p^*\dashv\Pi_p$, we get a corresponding map
  \begin{equation*}
    \vcenter{\xymatrix@=2em{
        {E_b}\ar[rr]^{f}\ar[dr]_{e_b} & & F\ar[dl]^{q} \\
        & {E} &
      }}
    \mkern20mu\text{where}\mkern20mu
    \vcenter{\xymatrix@=2em{
        {E_b}\pbexcursion\ar[r]^{e_b}\ar[d]_{p_b} & {E}\ar[d]^{p} \\
        {X}\ar[r]_{b} & {B}
      }}
  \end{equation*}
  in $\sSet_{/E}$. The claim follows by specializing to the case $X=\Del^n$.
  \end{proof}
  \fi

To study the pushforward construction along a
  cocartesian fibration $p\colon\qE\tfib\qB$, we will
  replace the test objects $e_b\colon \qE_b\to \qE$ involved the description of
  $\Pi_p$ given in Lemma \ref{lem:strict-test-objects} 
   by weakly equivalent test objects
  in the Joyal model structure. These new test objects $\tilde{e}_b\colon\tilde{E}_b\to \qE$
  will arise as certain weighted colimits of a fixed diagram
  $c_p\colon \gC{\qB}\to \qCat\subseteq\SSet$, namely the straightening of the cocartesian fibration $p$ defined using the  comprehension construction.

The construction of the replacement to the pullback functor is given in \S\ref{sec:alt-pullback}, and the proof that the pullback replacement is equivalent to the strict pullback is given in \S\ref{sec:pullback-comparison}. 

\subsection{A replacement for pullback along a cocartesian fibration}\label{sec:alt-pullback}

\begin{ntn}\label{ntn:straightened-functor}
  For the remainder of this section shall fix a cocartesian fibration of
  quasi-categories $p\colon \qE\tfib \qB$ as well as a corresponding comprehension functor
  \begin{equation*}
    \xymatrix@R=0em@C=5em{ \gC{\qB}\ar[r]^-{c_p} &
      {\qCat}, }
  \end{equation*}
  the ``straightening'' of the cocartesian fibration $p$.
  We also fix the associated lifted lax cocone with nadir $\qE$ described in Theorem \ref{thm:comprehension}:
      \begin{equation}\label{eq:lax-cocone-under-comprehension}
    \xymatrix@=1.5em{
      & {\catone+\gC{\qB}}\ar@{^(->}[dl]\ar[dr]^{\langle \qE,c_p \rangle} & \\
      {\gC[\qB\join\Del^0]}\ar[rr]_-{\ell^\qE} && {\qCat}
    }
  \end{equation}
\end{ntn}

Our aim is to prove that this lax cocone is a colimit cocone. We will achieve this as the $b=\id_{\qB}$ special case of our first main theorem:

\begin{thm}\label{thm:pullback-equivalence}
For any cocartesian fibration $p \colon \qE \tfib \qB$ and any $b \colon X \to \qB$, the comprehension cocone induces a canonical map over $\qE$ from the oplax colimit of the diagram
\[ \gC{X} \xrightarrow{\gC{b}} \gC{\qB} \xrightarrow{c_p} \qCat\]
to the fiber 
\[
\xymatrix{ E_b \ar[d] \ar[r] \pbexcursion & \qE \ar@{->>}[d]^p \\ X \ar[r]_b & \qB}
\]
and this map is a natural weak equivalence in the Joyal model structure.
\end{thm}

Before proving this result, we tighten up its statement. As we explain presently, there is a functor $\tilde{p}^* \colon \sSet_{/\qB} \to \sSet_{/\qE}$ that acts on objects by carrying a generalized element $b \colon X \to \qB$ to a canonical map $\colim^{\mathrm{oplax}}{(c_p \circ \gC{b})} \to \qE$. After defining this more formally, we construct a comparison natural transformation
\begin{equation}\label{eq:comp-nat-trans}
\begin{tikzcd} \sSet_{/\qB} \arrow[r, bend left, "\tilde{p}^*"] \arrow[r, bend right, "p^*"'] \arrow[r, phantom, "\scriptstyle\Downarrow\gamma"] & \sSet_{/\qE}
\end{tikzcd}
\end{equation}
Theorem \ref{thm:pullback-equivalence} asserts that this map is a componentwise Joyal equivalence. We first describe the action of the functor $\tilde{p}^* \colon \sSet_{/\qB} \to \sSet_{/\qE}$ on objects before establishing the functoriality of this construction. Recall from Remark \ref{rmk:generalized-unstraightening} that the comprehension functor $c_p \colon \gC{\qB} \to \qCat$ can be used to define a ``straightening'' of its pullbacks:
\[ 
    \xymatrix@R=0em@C=4em{
      {\gC{X}}\ar[r]^{\gC{b}} &
      {\gC{\qB}}\ar[r]^-{c_p} & \SSet
    }\]
    even in the case where $X$ is not a quasi-category.

\begin{defn}\label{defn:cocart-pullback-replacement} Given an generalized element $b \colon X \to \qB$ in $\sSet_{/\qB}$, define a simplicial set
\[ \tilde{E}_b \defeq \colim\nolimits^{\mathrm{oplax}}{\left(\begin{tikzcd}[ampersand replacement=\&] \gC{X} \arrow[r, "\gC{b}"] \& \gC{\qB} \arrow[r, "c_p"] \& \SSet\end{tikzcd}\right)}.\]
The oplax colimit $\tilde{\qE}_b$ is the nadir of the universal lax cocone under the diagram $c_p \circ \gC{X}$. This is the weighted colimit weighted by the weight for oplax colimits $L_X$ introduced in Definition \ref{defn:oplax-colimit-weight}.

The simplicial functor 
  \begin{equation}\label{eq:comprehension-cone}
\begin{tikzcd}[sep=small] & \catone+\gC{X}  \arrow[dl, hook'] \arrow[rr, "\catone+\gC{b}"] & &  \catone+\gC{\qB}  \arrow[dl, hook'] \arrow[dr, "{\langle E, \qE \rangle}"]  \\  {\gC[X\join\Del^0]} \arrow[rr, "{\gC[b\join\id]}"'] & & {\gC[\qB\join\Del^0]}\arrow[rr, "\ell^\qE"'] && \SSet\end{tikzcd}
    \end{equation}
defines a lax cocone $\ell^\qE\vert_b \colon \gC{[X\join\Del^0]} \to \SSet$ under the diagram $c_p \circ \gC{b}$ with nadir $\qE$, inducing a unique simplicial map $\tilde{e}_b \colon \tilde{E}_b \to \qE$ from the oplax colimit. This constructs an object of $\SSet_{/\qE}$.
\end{defn}

To establish the functoriality of this construction, it will be convenient to re-express the oplax colimits of Definition \ref{defn:cocart-pullback-replacement}.

\begin{lem}\label{lem:oplax-comp-colimit-fun} For any simplicial map $b \colon X \to B$ define a weight $L_b \colon \gC{\qB}\op \to \SSet$ by taking the left Kan extension along $\gC{b} \colon \gC{X}\op \to \gC{\qB}\op$ of the weight for oplax colimits. 
\begin{enumerate}[label=(\roman*)]
\item\label{itm:oplax-comp-colimit-fun-weight}  Then for any diagram $F \colon \gC{\qB} \to \SSet$, there is an isomorphism
\[ \colim\nolimits^{\mathrm{oplax}}{(F \circ \gC{b})} \coloneq \colim\nolimits^{L_{X}}{(F \circ \gC{b})} \cong \colim\nolimits^{L_b}{F}.\]
\item\label{itm:oplax-comp-colimit-fun-collage} The weight $L_b \colon \gC{\qB}\op \to \SSet$  is flexible and its collage is given by the pushout
\[
\begin{tikzcd} \catone+\gC{X} \arrow[d, hook] \arrow[dr, phantom, "\ulcorner" very near end] \arrow[r, "\catone+\gC{b}"] & \catone+ \gC{\qB}  \arrow[d, hook] \\ \gC{[X\join\Del^0]} \arrow[r] & \gC((X\join\Del^0)\! \underset{X+\Del^0}{\cup} (B+\Del^0)) \cong \coll L_b
\end{tikzcd}
\]
\end{enumerate}
\end{lem}
\ifnewadd
\begin{proof} Statement \ref{itm:oplax-comp-colimit-fun-weight} and the second part of \ref{itm:oplax-comp-colimit-fun-collage} follow from Lemma \ref{lem:lan-of-weight}.  Since $\coll(L_b)$ is the homotopy coherent realization of the pushout  of simplicial sets, Proposition \ref{prop:gothic-C} tells us that it is a simplicial computad and thus, by Proposition \ref{prop:flexible-weight-as-computad}, $L_b$ is a flexible weight.
  \end{proof}\fi 
\ifnewcut\else
\begin{proof}
By Lemma \ref{lem:lan-of-weight}, the weighted colimit of a restricted diagram is isomorphic to the colimit of the original diagram 
weighted by the left Kan extension of the weight. This specializes to prove \ref{itm:oplax-comp-colimit-fun-weight}.

By Lemma \ref{lem:lan-collage-pushout} the collage of the left Kan extended weight $L_b$ is computed by the pushout of \ref{itm:oplax-comp-colimit-fun-collage}. Since $\coll(L_b)$ is the homotopy coherent realization of the pushout  of simplicial sets, Proposition \ref{prop:gothic-C} tells us that it is a simplicial computad and thus, by Proposition \ref{prop:flexible-weight-as-computad}, $L_b$ is a flexible weight.
\end{proof}\fi

\begin{obs}\label{obs:oplax-comp-colimit-functoriality}
The utility of Lemma \ref{lem:oplax-comp-colimit-fun} is as follows. Suppose now that we have a map
\[\begin{tikzcd}[sep=small] X \arrow[dr, "b"'] \arrow[rr, "u"] & & Y \arrow[dl, "c"] \\ & \qB \end{tikzcd}\]
in the slice category $\sSet_{/\qB}$. This gives rise to a commutative diagram of simplicial computads
\begin{equation}\label{eq:collage-of-map-pushout}
\begin{tikzcd} \gC{[X\join\Del^0]} \arrow[d, "\gC{u^\triangleright}"'] & \catone+ \gC{X} \arrow[l, hook'] \arrow[d, "\catone+\gC{u}"'] \arrow[r, "\catone+\gC{b}"] &\catone+ \gC{\qB} \arrow[d, equals] \\ \gC[Y\join\Del^0] & \catone+\gC{Y} \arrow[l, hook'] \arrow[r, "\catone+\gC{c}"'] & \catone+\gC{\qB}
\end{tikzcd}
\end{equation}
inducing a simplicial computad morphism $\coll(L_b) \to \coll(L_c)$ in the category $\slicel{\catone+\gC{\qB}}{\sCptd}$.

This construction is functorial, defining the horizontal functor in the following square
\[
\begin{tikzcd} \sSet_{/\qB} \arrow[r, "\coll L_\bullet"] \arrow[d, "L_\bullet"', dashed] & \slicel{\catone+\gC{\qB}}{\sCptd} \arrow[d, hook] \\ \sSet^{\gC{\qB}\op} \arrow[r, "\coll", hook] & \slicel{\catone+\gC{\qB}}{\sCat}
\end{tikzcd}
\]
By the description of the essential image of the collage functor given in Proposition \ref{prop:collage-adjunction}, we see that $\coll L_\bullet$ factors as indicated defining a functor $L_\bullet \colon \sSet_{/\qB} \to \slicel{\catone+\gC{\qB}}{\sCat}$.

Finally note that the collage $\coll L_{\qB\op} \cong \gC{\qB\join\Del^0}$ for the weight for oplax colimits of shape $\qB$ defines a cone under the pushout diagram of Lemma \ref{lem:oplax-comp-colimit-fun}\ref{itm:oplax-comp-colimit-fun-collage}. Thus the codomain of the functor $\coll(L_\bullet)$ lifts to the slice category
\[ 
\begin{tikzcd}  \sSet_{/\qB} \arrow[r, "\coll L_\bullet"] & \left(\slicel{\catone+\gC{\qB}}{\sCptd}\right)_{/\gC{\qB\join\Del^0}}
\end{tikzcd}
\]
Correspondingly, by the fully faithfulness of the collage construction, we can equally regard $L_\bullet$ as a functor
\[ 
\begin{tikzcd}  \sSet_{/\qB} \arrow[r, "L_\bullet"] & \left(\sSet^{\gC{\qB}\op}\right)_{/L_{\qB}}
\end{tikzcd}
\]
landing in the full subcategory spanned by the flexible weights.
\end{obs}

Observation \ref{obs:oplax-comp-colimit-functoriality} allows us to extend Definition \ref{defn:cocart-pullback-replacement} to a functor.

\begin{defn}\label{defn:cocart-pullback-replacement-fun} Define $\tilde{p}^* \colon \sSet_{/\qB} \to \sSet_{/\qE}$ to be the composite functor
\[
\tilde{p}^* \coloneq \begin{tikzcd} \sSet_{/\qB} \arrow[r, "L_\bullet"] & \left(\sSet^{\gC{\qB}\op}\right)_{/L_{\qB\op}} \arrow[r, "\colim^{-}{c_p}"] & \sSet_{/\tilde{\qE}_\qB} \arrow[r, "\tilde{e}_\qB"] & \sSet_{/\qE}
\end{tikzcd}
\]
where $\tilde{\qE}_\qB \coloneq \colim^{\mathrm{oplax}}{c_p}$ and $\tilde{e}_\qB \colon \tilde{\qE}_\qB \to \qE$ is the map induced by the lax cocone \eqref{eq:lax-cocone-under-comprehension}.
\end{defn}

For later use, we record a few properties of the functor just constructed.

\begin{lem}\label{lem:pullback-replacement-colimits} The functor $\tilde{p}^* \colon \sSet_{/\qB} \to \sSet_{/\qE}$ preserves colimits.
\end{lem}
\begin{proof}
In Definition \ref{defn:cocart-pullback-replacement-fun} the functor under consideration is defined as a composite of three functors, the latter two of which manifestly preserve colimits. Since colimits in a slice category over an object are created by the forgetful functor, it remains only to prove that the functor $L_\bullet \colon \sSet_{/\qB} \to\sSet^{\gC{\qB}\op}$ preserves colimits. Since Proposition \ref{prop:collage-adjunction} demonstrates that the inclusion $\sSet^{\gC{\qB}\op}  \inc  \slicel{\catone+\gC{\qB}}{\sCat}$ is full and coreflective, to show that this functor preserves colimits, it suffices to show that 
\[ 
\begin{tikzcd}  \sSet_{/\qB} \arrow[r, "\coll L_\bullet"] & \slicel{\catone+\gC{\qB}}{\sCat}
\end{tikzcd}
\]
preserves them.

By Observation \ref{obs:oplax-comp-colimit-functoriality}, the action of this functor on objects and morphisms is defined by the pushout of Lemma \ref{lem:oplax-comp-colimit-fun}\ref{itm:oplax-comp-colimit-fun-collage}, which we regard as a diagram in $\slicel{\catone}{\sCat}$. The functors
\[
\begin{tikzcd} \sSet_{/\qB} \arrow[r, "\catone+\gC(-)"] & \slicel{\catone}{\sCat}
\end{tikzcd} \qquad \text{and} \qquad
\begin{tikzcd} \sSet_{/\qB} \arrow[r, "\gC(-)\join\Del^0"] & \slicel{\catone}{\sCat}
\end{tikzcd} 
\]
both preserve colimits. Thus, the functor from $\sSet_{/\qB}$ to the category of pushout diagrams in $\slicel{\catone}{\sCat}$ with one vertex fixed at $\catone+\gC{\qB}$ preserves colimits. The pushout preserves colimits as well so we conclude that $\coll(L_\bullet)$ and hence $\tilde{p}^*\colon \sSet_{/\qB} \to \sSet_{/\qE}$ preserves colimits, as desired.
\end{proof}

\begin{lem}\label{lem:pullback-replacement-monos} The functor $\tilde{p}^* \colon \sSet_{/\qB} \to \sSet_{/\qE}$ preserves monomorphisms.
\end{lem}
\begin{proof}
From the Definition \ref{defn:cocart-pullback-replacement-fun}, to see that $\tilde{p}^*$ preserves monomorphisms \begin{equation*}
\begin{tikzcd}[sep=small] X \arrow[rr, "u", hook] \arrow[dr, "b"'] & & Y \arrow[dl, "c"] \\ & \qB
\end{tikzcd}
  \end{equation*}
 over $B$ it suffices to show that the composite of the functors $L_\bullet$ and $\colim_{\gC{\qB}}(-,c_p)$ preserve monomorphisms. To do so, we'll prove that the comparison functor $L_u \colon L_b \to L_c$ between weights in $\sSet^{\gC{\qB}\op}$ is a relative projective cell complex, as defined in \ref{defn:proj-cell-cx}. A theorem of Gambino \cite{Gambino:2010wl}  implies that $\colim_{\gC{\qB}}(-,c_p)$ carries relative projective cell complexes to monomorphisms in simplicial sets; see also \cite[11.5.1]{Riehl:2014kx}
 
 Recall that the weight $L_b$ is constructed as a collage defined by a pushout, which is the homotopy coherent realization of a pushout of simplicial sets.  The natural transformation $L_u$ is  encoded by the map between collages constructed as the pushout \eqref{eq:collage-of-map-pushout}; again this map is the homotopy coherent realization of a map of simplicial sets.  Since the left-hand horizontal inclusions are also simplicial subcomputad inclusions, it follows from the standard argument that the induced 
 map $\coll(L_u) \colon \coll(L_b) \inc \coll(L_c)$ between the pushouts is a simplicial subcomputad inclusion and by the relative analogue Proposition~\refII{prop:projcofchar} of Proposition \ref{prop:flexible-weight-as-computad}, $L_u \colon L_b\to L_c$ is a relative projective cell complex, as desired.
\end{proof}

\subsection{Comparison with the strict pullback}\label{sec:pullback-comparison}

Now that we've precisely defined a functor $\tilde{p}^* \colon \sSet_{/\qB} \to \sSet_{/\qE}$ that carries a generalized element to the oplax colimit of the restricted comprehension functor, our next task is to define the natural transformation \eqref{eq:comp-nat-trans} alluded to in the statement of Theorem \ref{thm:pullback-equivalence}.  To explain the existence of the natural map $\gamma_b \colon \tilde{E}_b \to E_b$ for $b \colon X \to \qB$, recall that the $p$-cocartesian lifts with
  codomain $\qE$ used to define the action of $c_p\colon \gC{\qB}\to\qCat$ on
  arrows in the image of $\gC{b}\colon \gC{X}\to\gC{\qB}$ lie over arrows with codomain $\qB$ which
  have a given factorisation through $b\colon X\to \qB$. This is depicted in the following diagram by the arrow $b\gamma$ and its $p$-cocartesian lift $\ell^\qE_{b\gamma}$:
  \begin{equation*}
    \begin{xy}
      0;<1.4cm,0cm>:<0cm,0.75cm>::
      *{\xybox{
          \POS(1,-0.5)*+{\Del^0}="one"
          \POS(0,1.5)*+{\Del^0}="two"
          \POS(2.25,0.5)*+{X}="three"
          \POS(4,0.5)*+{\qB}="four"
          \POS(1,4)*+{\qE_{bx'}}="one'"
          \POS(1,3.75)*{\pbcorner}
          \POS(0,6)*+{\qE_{bx}}="two'"
          \POS(0,6)*{\pbcorner}
          \ar@{=} "one";"two"
          \ar@/_5pt/ "one";"three"_{x'}^(0.15){}="otm"
          \ar@/^10pt/ "two";"three"^{x}_(0.6){}="ttm"|(0.42){\hole}
          \ar "three";"four"_(0.5){b}
          \ar@{=>} "ttm"-<0pt,7pt> ; "otm"+<0pt,10pt> ^(0.3){\gamma}
          \POS(2.25,2.5)*+{E_b}="three'"
          \POS(4,2.5)*+{\qE}="four'"
          \POS(2.25,2.5)*{\pbcorner}
          \ar "three'";"four'"_{\ell^\qE_b}
          \ar@{->>} "four'";"four"^{p}
          \ar@/_5pt/ "one'";"four'"^(0.6){\ell^\qE_{bx'}}|(0.42){}="otm'"
          \ar@/^15pt/ "two'";"four'"^{\ell^\qE_{bx}}_(0.58){}="ttm'"
          \ar@{->>} "one'";"one"_(0.325){p_{bx'}}
          \ar@{->>} "two'";"two"_{p_{bx}}
          \ar@{->>} "three'";"three"^{p_b}
          \ar "two'";"one'"_*!/^2pt/{\scriptstyle e_{b\gamma}}
          \ar@{=>} "ttm'"-<0pt,4pt> ; "otm'"+<0pt,4pt> ^(0.3){\ell^\qE_{b\gamma}}
          \ar@{.>}@/_6pt/ "one'";"three'"^(0.4){}="otm''"
          \ar@{.>}@/^12pt/ "two'";"three'"_(0.65){}="ttm''"
          \ar@{:>} "ttm''";"otm''"
        }}
    \end{xy}
  \end{equation*}
From the diagram, it is clear that such $p$-cocartesian arrows factor through $e_b \colon
  E_b\to \qE$ to give the
  dotted arrows with codomain $E_b$ as drawn; this is the main component of the proof of Lemma \ref{lem:cocart-cocone-pb}. This idea is formalized as follows: 

\begin{lem} For any $b \colon X \to \qB$, the diagram $c_p \circ \gC{b}$ is the base of a lax cocone with nadir $E_b$
\[
\begin{tikzcd}[sep=small] & \catone+\gC{X} \arrow[dl, hook'] \arrow[dr, "{\langle E_b,c_p\circ \gC{b}\rangle}"] \\ \gC[X\join\Del^0] \arrow[rr, "\ell^E"'] & & \SSet
\end{tikzcd}
\]
Hence, the universal property of the oplax colimit defines a natural map $\tilde{E}_b \to E_b$ over $E$.
\end{lem}
\begin{proof}
Apply Lemma \ref{lem:cocart-cocone-pb} to the cocartesian cone of Theorem \ref{thm:comprehension} as described in Remark \ref{rmk:arbitrary-pullbacks}.
 \end{proof}

To prove Theorem \ref{thm:pullback-equivalence} we must verify that   \eqref{eq:comp-nat-trans} is a componentwise Joyal weak equivalence. We first demonstrate this for generalized elements $b \colon \Del^n \to B$ whose domains are simplices and then use the results of Lemmas \ref{lem:pullback-replacement-colimits} and \ref{lem:pullback-replacement-monos} to extend these results to the general case.

\begin{ex}\label{ex:pullback-equivalence-0} By definition $\tilde{p}^*(b \colon \Del^0 \to \qB)$ is the oplax colimit of the diagram
\[
\begin{tikzcd} \gC\Del^0 \arrow[r, "\gC{b}"] & \gC{\qB} \arrow[r, "c_p"] & \qCat \end{tikzcd}
\]
that sends the unique object to the fiber $\qE_b$ of $p \colon \qE \tfib \qB$ over $b \colon \Del^0 \to \qB$. The weight  for lax cocones of shape $\Del^0$ is the terminal weight so the weighted colimit is just the ordinary colimit of this one object diagram. Thus $\tilde{p}^* \colon \sSet_{/\qB} \to \sSet_{/\qE}$ sends $b \colon \Del^0 \to \qB$ to $\qE_b \to \qE$, which is isomorphic to the strict pullback $p^*(b \colon \Del^0 \to \qB)$.
\end{ex}

For $b \colon \Del^1 \to \qB$, $\tilde{p}^*(b \colon \Del^1 \to \qB)$ is the oplax colimit of the diagram 
\[
\begin{tikzcd} \gC\Del^1 \arrow[r, "\gC{b}"] & \gC{\qB} \arrow[r, "c_p"] & \qCat \end{tikzcd}
\]
 whose image is diagram $e_b \colon \qE_{b_0} \to \qE_{b_1}$ of quasi-categories constructed in \eqref{eq:comprehension-over-arrow}. In this case, the oplax colimit has a simple description: it is given by an ``mapping cylinder'' formed by attaching $\qE_{b_1}$ along the codomain edge of the cylinder $\qE_{b_0} \times \Del^1$ via the map $e_b$. We now show this simplicial set is Joyal weak equivalent to the strict fiber $\qE_b$.
  
\begin{prop}\label{prop:pullback-equivalence-1} The data formed by applying the comprehension construction 
\begin{equation}\label{eq:comprehension-over-arrow}
\begin{tikzcd}[row sep=tiny] \qE_{b_0} \arrow[rrrd, bend left=10, ""{name=1,right}, "\ell^\qE_0"] \arrow[ddr, dotted, "e_b"'] \arrow[ddd, "p_0"'] \\  & & & \qE_b \arrow[ddd, two heads, "p_b"] \arrow[rr] \arrow[dddrr, phantom, "\lrcorner" very near start] & & \qE \arrow[ddd, two heads, "p"]  \\ 
& \qE_{b_1} \arrow[urr, bend right=10, "\ell^\qE_1"'] \arrow[from=1,Rightarrow, dotted, shorten >= 0.5em, shorten <= 0.5em, "\chi"]\arrow[ddrr, phantom, "\lrcorner" pos=.05]  \\
1 \arrow[rrrd, bend left=10, ""{name=2,right}, "0"] \arrow[ddr, equals]  \\ & & & \Del^1 \arrow[rr, "b"] & & \qB\\ 
& 1 \arrow[urr, bend right=10, "1"'] \arrow[from=uuu, crossing over, two heads, "p_1"' pos=.6] \arrow[from=2,Rightarrow,  shorten >= 0.5em, shorten <= 0.5em, "\kappa"] 
\end{tikzcd}
\end{equation}
to the pullback $p_b \colon \qE_b \to \Del^1$ of $p \colon \qE \tfib \qB$ along $b \colon \Del^1 \to \qB$
 induces a Joyal weak equivalence
\[ 
\begin{tikzcd} {\qE_{b_0}} \arrow[r, "e_b"] \arrow[d, hook, "\id \times \face^0"']  \arrow[dr, phantom, "\ulcorner" very near end] & {\qE_{b_1}} \arrow[d, hook] \arrow[ddr, bend left, "\ell^\qE_1"]  \\ {\qE_{b_0}} \times \Del^1 \arrow[r] \arrow[drr, bend right, "\chi"']& \tilde{\qE}_b \arrow[dr, dashed, "\gamma_b", "\rotatebox{90}{$\sim$}"']  \\ & & \qE_b
\end{tikzcd}
\] 
\end{prop}
\begin{proof}
By Proposition \ref{prop:oplax-colimit-of-functor}, the oplax colimit $\tilde{\qE}_b$ is weakly equivalent to the quasi-categorical collage $\coll(e_b,\qE_{b_1})$, introduced in Definition \ref{defn:qcat-collage}. Moreover, Proposition \ref{prop:oplax-colimit-of-functor}   demonstrates that the equivalence $k \colon \tilde{\qE}_b \we \coll(e_b,\qE_{b_1})$  is inner anodyne. In particular, there exists a lift
\[
\begin{tikzcd} \tilde{\qE}_b \arrow[d, tail, "\rotatebox{90}{$\sim$}", "k"'] \arrow[r, "\gamma_b"] & \qE_b \arrow[d, "p_b", two heads] \\ \coll(e_b,\qE_{b_1}) \arrow[r, "\rho"'] \arrow[ur, dashed, "\ell"'] & \Del^1
\end{tikzcd}
\]
 defining a direct comparison map $\ell \colon \coll(e_b,\qE_{b_1}) \to \qE_b$ over $\Del^1$.
 
 To prove that $\ell$ is an equivalence, observe by Lemma \ref{lem:collage-qcat-cocartesian} and Proposition \ref{prop:pullback-stability} that $\rho$ and $p_b$ are both cocartesian fibrations. Hence, Proposition \ref{prop:equivalence-of-fibrations} tells us that if $\ell$ is a cartesian functor, then to demonstrate that $\ell$ is an equivalence, we need only show that it restricts to an equivalence on the fibers over 0 and 1. Indeed, $\ell$ is an isomorphism on both fibers, so now our only remaining task is to demonstrate that it is a cartesian functor.
 
 The proof of Lemma \ref{lem:collage-qcat-cocartesian} reveals that the non-degenerate cocartesian edges of $\coll(e_b,\qE_{b_1})$ are those represented by the degenerate edge of some vertex lying over the non-degenerate 1-simplex in $\Del^1$. Such edges lie in the image of the functor $\qE_{b_0}\times \Del^1 \to \qE_b$ used to define the map $\tilde{E}_b \to \qE_b$, and this functor in turn is defined to be a representative for the cocartesian lift of the 1-arrow between the two objects of $\Del^1$ to a map with domain $\qE_{b_0} \to \qE_b$. In particular, it defines a cocartesian cylinder in the sense of Lemma \ifnewcut\else \ref{lem:ptwise-cocart-cyl}\fi\ifnewadd\refVI{lem:cocart-cyl-are-cocart-arrows}\fi, which tells us that its components indexed by vertices of $\qE_{b_0}$ are $p$-cocartesian 1-arrows.
 This proves that $\ell$ carries $\rho$-cocartesian arrows to $p$-cocartesian arrows, and thus $\ell$ is an equivalence.
 \end{proof}

We argue inductively that $\gamma_b \colon \tilde{\qE}_b \to \qE_b$ is an equivalence for any $n$-simplex $b \colon \Del^n \to \qB$ under the assumption that this is true for simplices of lower dimension. Our strategy mirrors that adopted for the 1-simplex: we construct a quasi-categorical model for the oplax colimit of a homotopy coherent diagram $c_p \circ \gC{b} \colon \gC\Del^n \to \qCat$, i.e., a quasi-category equivalent to the simplicial set $\tilde{\qE}_b$ defined as the oplax colimit of $c_p \circ \gC{b}$, and then show that this is equivalent to the strict pullback $\qE_b$. The inductive step makes use of the following weights.

\begin{ntn}[weights for the inductive comparison]\label{ntn:oplax-simplex-weight} To compare the weights $L_{\Del^{n-1}}$ and $L_{\Del^{n}}$ for oplax colimits of a homotopy coherent $n-1$-simplex and $n$-simplex, we left Kan extend the former along the inclusion $\face^{n} \colon (\gC{\Del^{n-1}})\op\inc (\gC{\Del^{n}})\op$, writing $L_{\Del^{n-1}} \colon (\gC\Del^{n})\op \to \sSet$ for the left Kan extension of $L_{\Del^{n-1}}$. Explicitly, this weight is defined by 
\[
\xymatrix@R=5pt{ L_{\Del^{n-1}} \colon (\gC{\Del^{n}})\op \ar[r] & \SSet \\ i \ar@{}[r]|-{\mapsto} & { \begin{cases} \Fun_{\gC{\Del^{n}}}(i,n) & i < n \\ \emptyset & i = n \end{cases}}}
\]
Let $Y^{n}$ denote the representable weight
\[
\xymatrix@R=5pt{ Y^{n} \colon  (\gC{\Del^{n}})\op\ar[r] & \SSet \\ i \ar@{}[r]|-{\mapsto} & \Fun_{\gC{\Del^n}}(i,n)}
\]
Note there is a natural inclusion $L_{\Del^{n-1}}\inc Y^{n}$ that is the identity in all components except the one indexed by the object $n \in \gC{\Del^{n}}$.
\end{ntn}

\begin{lem}\label{lem:induction-all-in-the-weights} $\quad$ 
\begin{enumerate}[label=(\roman*)]
\item\label{itm:weight-induction-i} The following diagram defines a pushout of weights in $\SSet^{(\gC\Del^n)\op}$:
\[
\begin{tikzcd} L_{\Del^{n-1}} \arrow[r, hook] \arrow[d, "\id \times \face^0"'] \arrow[dr, phantom, "\ulcorner" very near end] & Y^n \arrow[d] \\ L_{\Del^{n-1}} \times \Del^1 \arrow[r] & L_{\Del^n}
\end{tikzcd}
\]
\item\label{itm:weight-induction-ii} Let $F \colon \gC\Del^n \to \SSet$ be a homotopy coherent diagram whose $L_{\Del^{n-1}}$-weighted colimit is Joyal weakly equivalent to the simplicial set $\qE_{n-1}$. Then the oplax colimit of $F$ is Joyal weakly equivalent to the pushout along a canonical map $\iota_n$ induced by the diagram $F$.
\[ 
\begin{tikzcd} \qE_{n-1} \arrow[r, "\iota_n"] \arrow[d, "\id \times \face^0"'] \arrow[dr, phantom, "\ulcorner" very near end] & F_n \arrow[d] \\ \qE_{n-1} \times \Del^1 \arrow[r] & \qE_n
\end{tikzcd}
\]
\end{enumerate}
\end{lem}
\begin{proof}
The pushout in \ref{itm:weight-induction-i} can be verified componentwise at each $i \in \gC\Del^n$ at which point this relationship is evident from the definitions. 

The pushout of \ref{itm:weight-induction-ii} follows.  If $\qE_{n-1}$ is isomorphic to the $L_{\Del^{n-1}}$-weighted colimit of $F$, then the pushout diagram of  \ref{itm:weight-induction-ii} is obtained by applying the cocontinuous functor $\colim^-\!F$ to the pushout diagram of  \ref{itm:weight-induction-i}. In this case, the map $\iota_n$ has a natural explicit description. By Lemma \ref{lem:lan-of-weight}, the $L_{\Del^{n-1}}$-weighted colimit of $F$ coincides with the oplax colimit of the restricted diagram $F \circ \gC(\face^n)\op$. The functor $F$ itself defines a canonical lax cocone under this restricted diagram with nadir $F_n$. Hence there is a natural comparison $\iota_n$ from the $L_{\Del^{n-1}}$-weighted colimit to $F_n$. 

Observe from Proposition \ref{prop:flexible-weight-as-computad} and Lemma \ref{lem:oplax-is-flexible} that all of the weights appearing in \ref{itm:weight-induction-i} are flexible.  Proposition \ref{prop:flexible-weights-are-htpical} then demonstrates that the pushout being constructed is equivalence-invariant.
\end{proof}

This lemma provides the inductive step in the following computation:

\begin{prop}\label{prop:pullback-equivalence-on-simplices} For any simplex $b \colon \Del^n \to \qB$, the component $\gamma_b \colon \tilde{\qE}_b \to \qE_b$ from the oplax colimit to the strict pullback is a Joyal weak equivalence.
\end{prop}
\begin{proof}
The base cases for $n=0$ and $n=1$ appear as Example \ref{ex:pullback-equivalence-0} and Proposition \ref{prop:pullback-equivalence-1}. For the induction step, suppose we have shown this is a componentwise weak equivalence for all $n-1$-simplices in $\qB$. By Lemma \ref{lem:lan-of-weight} and Notation \ref{ntn:oplax-simplex-weight}, the $L_{\Del^{n-1}}$-weighted colimit of the diagram $\gC\Del^n \xrightarrow{\gC{b}} \gC{\qB} \xrightarrow{c_p} \sSet$ is isomorphic to the oplax weighted colimit of the restricted diagram
\[ 
\begin{tikzcd} \gC\Del^{n-1} \arrow[r, "\gC{\face^n}"] & \gC\Del^n \arrow[r, "\gC{b}"] & \gC{\qB} \arrow[r, "c_p"] & \sSet.
\end{tikzcd}
\]
By the inductive hypothesis, this weighted colimit $\tilde{\qE}_{b\cdot\face^n}$ is weakly equivalent to the pullback $\qE_{b\cdot\face^n}$. By Lemma \ref{lem:induction-all-in-the-weights}, the diagram
\[ 
\begin{tikzcd} \qE_{b\cdot\face^n} \arrow[r, "\iota_n"] \arrow[d, "\id \times \face^0"']  & \qE_{b_n} \arrow[d] \\ \qE_{b\cdot\face^n} \times \Del^1 \arrow[r] & \tilde{\qE}_b
\end{tikzcd}
\]
is then a pushout up to Joyal weak equivalence. So it follows from Proposition \ref{prop:oplax-colimit-of-functor} that $\tilde{\qE}_b$ is equivalent to the  quasi-categorical collage $\coll(\iota_n, \qE_{b_n})$, and as in the proof of Proposition \ref{prop:pullback-equivalence-1}, the map $\gamma_b$ factors to define a map 
\[
\begin{tikzcd} \coll(\iota_n, \qE_{b_n}) \arrow[dr, two heads, "\rho"'] \arrow[rr, dashed, "\ell"] & [-.5em]~ & \qE_b \arrow[dl, "\pi_\ell \cdot b"] \\ & \Del^1
\end{tikzcd}
\]
in this case involving the map $\pi_\ell \colon \Del^n \to \Del^1$ that carries every element but the last one to $0$.  Observe that $\pi_\ell$ a cocartesian fibration, and indeed a bifibration, as it is covariantly represented by the functor $! \colon [n-1] \to [0]$, which admits both left and right adjoints; see Corollary \ref{cor:lurie-adjunction}.

Our task, again, is to show that $\ell$ is an equivalence. By Lemma \ref{lem:collage-qcat-cocartesian} and the fact that cocartesian fibrations compose, it is a functor between cocartesian fibrations. Moreover, $\ell$ is bijective on the fibers over $0,1 \in \Del^1$,  the latter being $\qE_{b_n}$ in both cases and the former being $\qE_{b\cdot\face^n}$. As in the proof of Proposition \ref{prop:pullback-equivalence-1}, $\ell$ is a cartesian functor, so Proposition \ref{prop:equivalence-of-fibrations}  implies that $\ell$ is an equivalence, as desired.
\end{proof}

Combining the work in this section, we can finally prove our main result.

\begin{proof}[Proof of Theorem \ref{thm:pullback-equivalence}] 
Our task is to demonstrate that the canonical natural transformation
\[
\begin{tikzcd} \sSet_{/\qB} \arrow[r, bend left, "\tilde{p}^*"] \arrow[r, bend right, "p^*"'] \arrow[r, phantom, "\scriptstyle\Downarrow\gamma"] & \sSet_{/\qE}
\end{tikzcd}
\]
is a componentwise Joyal weak equivalence using the result of Proposition \ref{prop:pullback-equivalence-on-simplices}, which demonstrates that this is the case for the simplices $b \colon \Del^n \to \qB$ of $\qB$.

The category $\sSet_{/\qB}$ is equivalent to the category $\sSet^{\el{\qB}\op}$ of presheaves indexed by the category $\el{\qB}$ of simplices of $\qB$; its objects are simplices $b\colon\Del^n\to \qB$ and a morphism from $b$ to $c \colon \Del^m \to \qB$ is a simplicial operator $\alpha \colon\Del^n \to \Del^m$ so that $c \cdot \alpha = b$.  The representable presheaves generate $\sSet^{\el{\qB}\op}$ under colimits, and such colimits are preserved by both of the functors $p^*$ and $\tilde{p}^*$, the former case because of the right adjoint $\Pi_p$ that exists in the locally cartesian closed category $\sSet$ and the latter case by Lemma \ref{lem:pullback-replacement-colimits}. Under the equivalence $\sSet^{\el{\qB}\op} \cong \sSet_{/\qB}$, these representables correspond to the objects $b \colon \Del^n \to \qB$ whose domain is a simplex. Proposition \ref{prop:pullback-equivalence-on-simplices} verifies that the components of $\gamma$ indexed by such objects are equivalences, which is the moral reason why $\gamma$ is an equivalence at all objects.

To demonstrate this, note that $b \colon X \to \qB$ is a colimit indexed by the category $\el{X}$ of its simplices $\Del^n \xrightarrow{x} X \xrightarrow{b} \qB$, i.e.,
\[ (X \xrightarrow{b} \qB) \cong \colim_{\el{X}}(\Del^n \to \qB).\] The map $\gamma_b$ factors as
\[ \gamma_b \colon \tilde{p}^*(\colim_{\el{X}}\Del^n \to \qB) \cong \colim_{\el{X}}\tilde{p}^*(\Del^n\to \qB) \longrightarrow \colim_{\el{X}} p^*(\Del^n \to \qB) \cong p^*(\colim_{\el{X}}\Del^n \to \qB),\] so it remains only to show that this middle map, the colimit of the equivalences $\gamma_{bx}$ indexed by the simplices of $X$, is itself an equivalence. The indexing category $\el{X}$ is a Reedy category, so if we can show that the two $\el{X}$-indexed diagrams are Reedy cofibrant and that the category $\el{X}$ has fibrant constants, then the pointwise equivalence between the diagrams will induce the desired equivalence between their colimits. To say that the Reedy category $\el{X}$ has fibrant constants means that for each element $x \colon \Del^n \to X$, the category of elements of the covariant representable boundary functor $\partial\el{X}_x$ is either empty or connected. This category is empty just when $x$ is non-degenerate and has a terminal object, and is connected in particular, when $x$ is degenerate. So $\el{X}$ has fibrant constants and the colimit functor $(\sSet_{/\qB})^{\el{X}}\to \sSet_{/\qB}$ carries pointwise weak equivalences between Reedy cofibrant diagrams to weak equivalences.

To verify this Reedy cofibrancy, it suffices to show 
\begin{enumerate}[label=(\roman*)]
\item the canonical diagram $\el{X} \to \sSet_{/\qB}$ is Reedy cofibrant
\item $\tilde{p}^*$ and $p^*$ preserve Reedy cofibrant objects.
\end{enumerate}
Since $\tilde{p}^*$ and $p^*$ preserve colimits, they in particular preserve latching objects, so for this second item it  suffices to show that both functors also preserve monomorphisms. Here, the fact that the pullback functor $p^*$ preserves monomorphisms is standard, and the fact that its replacement $\tilde{p}^*$ preserves monomorphisms was proven in Lemma \ref{lem:pullback-replacement-monos}. 

So it remains only to prove (i), that is, to argue that the functor
\[ \xymatrix@R=1em@C=5em{ \el{X} \ar[r] & \sSet_{/\qB} \\ x \colon \Del^n \to X \ar@{}[r]|-\mapsto & bx \colon \Del^n \to \qB}\] is Reedy cofibrant. The latching object associated to $x \colon \Del^n \to X$ is the composite $\boundary\Del^n\inc\Del^n\xrightarrow{x} X\xrightarrow{b} \qB$ and the latching map is the inclusion $\partial\Del^n\inc\Del^n$ over $\qB$, which is obviously a monomorphism. This completes the proof.
\end{proof}
  
Specializing Theorem \ref{thm:pullback-equivalence} to the identity morphism on $\qB$, we have

\begin{cor}\label{cor:total-space-as-oplax-colimit}
The domain of a cocartesian fibration $p \colon \qE \tfib \qB$ is equivalent to the oplax colimit of the associated comprehension functor $c_p \colon \gC\qB \to \qCat$, with colimit cocone:
            \begin{equation*}
    \xymatrix@=1.5em{
      & {\catone+\gC{\qB}}\ar@{^(->}[dl]\ar[dr]^{\langle \qE,c_p \rangle} & \\
      {\gC[\qB\join\Del^0]}\ar[rr]_-{\ell^\qE} && {\qCat}
    } \qed
  \end{equation*}
  \end{cor}


%% file: exponentiation.tex

\section{Pushforward along a cocartesian fibration}\label{sec:powerful}

In this section, we shall fix a cocartesian fibration
  $p\colon \qE\tfib \qB$ of quasi-categories and prove that the pushforward functor 
  \[ \sSet_{/\qE} \xrightarrow{\Pi_p} \sSet_{/\qB}\]
   has two properties that are relevant to the development of the category of theory of quasi-categories:
  \begin{enumerate}[label=(\roman*)]
  \item\label{itm:pushforward-iso} The pushforward functor preserves \emph{isofibrations}. In model categorical terminology, this implies that the adjunction 
  \[
\begin{tikzcd} \sSet_{/\qE} \arrow[r, bend right, "{\Pi}_p"'] \arrow[r, phantom, "\bot"] & \sSet_{/\qB} \arrow[l, bend right, "{p}^*"'] 
\end{tikzcd}
\]
   is Quillen with respect to slices of the Joyal model structure.
  \end{enumerate}
  Moreover:
\begin{enumerate}[label=(\roman*), resume]
  \item\label{itm:pushforward-cart} The pushforward functor preserves \emph{cartesian fibrations} and \emph{cartesian functors} between them.
  \end{enumerate}
In fact, both pullback and pushforward along $p$ define \emph{cosmological functors}, a property we briefly note for use in future work.  
  
Both of the properties \ref{itm:pushforward-iso} and \ref{itm:pushforward-cart} are more easily established for an alternate model of the pushforward functor defined as a right adjoint to the functor $\tilde{p}^* \colon \sSet_{/\qB} \to \sSet_{/\qE}$ introduced in \S\ref{sec:pullback}.  Theorem \ref{thm:pullback-equivalence} demonstrates that the pullback $E_b \to \qE$ of a functor $b \colon X \to \qB$ along a cocartesian fibration $p \colon \qE \tfib \qB$ is computed, up to equivalence, as the oplax colimit of a particular diagram
\[
\begin{tikzcd} \gC{X}\arrow[r, "\gC{b}"] & \gC{\qB} \arrow[r, "c_p"] & \qCat
\end{tikzcd}
\]
When the oplax colimit is defined strictly as a simplicial set it enjoys the universal property of Definition \ref{defn:simp-weight}: maps $\colim^{\mathrm{oplax}}{(c_p \circ \gC{b})} \to \qF$ correspond to lax cocones under $c_p \circ \gC{b}$ with nadir $\qF$. This correspondence defines a right adjoint
\[
\begin{tikzcd} \sSet_{/\qE} \arrow[r, bend right, "\tilde{\Pi}_p"'] \arrow[r, phantom, "\bot"] & \sSet_{/\qB} \arrow[l, bend right, "\tilde{p}^*"']
\end{tikzcd}
\]
characterized on an object $q \colon \qF \to \qE$ by the bijection
\[
\begin{tikzcd}[column sep={5em,between origins}]
\Del^n \arrow[rr] \arrow[dr, "b"'] & & \tilde{\Pi}_p\qF \arrow[dl, "\tilde{\Pi}_pq"]   &~ & \colim\nolimits^{\mathrm{oplax}}(c_p \circ \gC{b}) \arrow[dr, "{\tilde{p}^*(b) \defeq \ell^E\vert_b}"'] \arrow[rr] & & \qF \arrow[dl, "q"] 
\\ & \qB & &  \arrow[u, phantom, "\leftrightsquigarrow"] & & \qE
\end{tikzcd}
\]
That is, $n$-simplices in $\tilde{\Pi}_p\qF$ over $b \colon \Del^n \to \qB$ correspond to lax cocones under the homotopy coherent $n$-simplex $c_p\circ\gC{b}$ with nadir $\qF$ whose whiskered composite with $q$ recovers the restriction $\ell^\qE\vert_b$ of the lax cocone produced by the comprehension construction of Theorem \ref{thm:comprehension}.

To make this simplex level construction of $\tilde{\Pi}_pq$ precise, we require a simplicial set whose $n$-simp\-lic\-es correspond to lax cocones under a homotopy coherent $n$-simplex in $\qCat$ with nadir $F$. One might think that the slice quasi-category $\slicer{\qqCat}{\qF}$ provides just such a gadget, where $\qqCat$ is the quasi-category of quasi-categories  defined by passing to the maximal Kan complex enriched core and then applying the homotopy coherent nerve,  but this isn't quite correct: since we've passed to the $(\infty,1)$-categorical core of $\qCat$ before taking the homotopy coherent nerve, simplices in $\slicer{\qqCat}{\qF}$ correspond to \emph{pseudo} cocones rather than lax cocones. The solution is to drop the core functor, in which case the homotopy coherent nerve $\qqCat_2 \defeq N\qCat$ is not a quasi-category but rather a \emph{2-complicial set}, a type of marked simplicial set which is introduced in \S\ref{sec:2complicial}.   Definition \ref{defn:join-slice} introduces a slice construction for marked simplicial sets, which does \emph{not} commute with the functor that forgets the markings, but this is a good thing. The marked slice $\laxslicer{\qqCat_2}{\qF}$ has exactly the property we desire, in that its $n$-simplices correspond to lax cocones under a homotopy coherent $n$-simplex in $\qCat$ with nadir $F$. 


In \S\ref{sec:right-adjoint}, we describe $\tilde\Pi_pq$ explicitly as the pullback of a map between lax slices of the homotopy coherent nerve of $\qCat$ defined by ``whiskering with $q$.''  After establishing the properties \ref{itm:pushforward-iso} and \ref{itm:pushforward-cart} for $\tilde{\Pi}_p$, we use the natural Joyal equivalence $\gamma \colon \tilde{p}^* \To p^*$ to transfer these properties to the pushforward functor $\Pi_p$. Having established that pullback and pushforward along a cocartesian fibration both preserve cartesian fibrations, in \S\ref{sec:exponentiation}, we construct a closely related exponentiation operation $(q \colon \qF \tfib \qB)^{p \colon \qE \tfib \qB}$ of a cartesian fibration $q$ by a cocartesian fibration $p$ with the same codomain. These exponentials are used in  \S\ref{sec:monadicity} to establish the comonadicity and monadicity of the quasi-category of cartesian fibrations over $\qB$ over the quasi-category of $\ob\qB$-indexed families of quasi-categories.

\subsection{2-complicial sets}\label{sec:2complicial}

  We know from Cordier and Porter~\cite{Cordier:1982:HtyCoh,
    Cordier:1986:HtyCoh} that the homotopy coherent nerve of a 
   Kan complex enriched category is itself a
  quasi-category. But when we apply the homotopy coherent nerve 
   to a quasi-category enriched category, such as $\qCat$ itself, it is
  \emph{not} the case that these nerves are quasi-categories. Since the
  hom-spaces of a quasi-category enriched category contain $1$-simplices that are not invertible, its homotopy coherent nerve
  contains 2-simplices which are not invertible. A homotopy coherent nerve of this kind is most naturally regarded as
  possessing the structure of a \ifnewcut\else\textit{$2$-complicial set\/} 
\cite{Verity:2007rm, Verity:2007:wcs1}.\fi \ifnewadd \textit{$2$-complicial set\/}.

We leave the precise definition to the original sources \cite{Verity:2007rm, Verity:2007:wcs1} or to more recent expository accounts such as \cite{Riehl:2017cs} or \S X.D.1 and instead present an overview of the main ideas. Extending the terminology used by Lurie in \cite[\S 3.1]{Lurie:2009fk}, a \textit{marked simplicial set} is a simplicial set equipped
with a chosen subset of \emph{marked simplices}, which must be positive-dimensional and contain all degeneracies. Maps of marked simplicial sets preserve the markings. A \textit{complicial set} is a marked simplicial set with the right lifting property with respect to certain marked horn inclusions---including both inner and outer marked horns---as well as certain marking extensions. A complicial set is \emph{saturated} if ``all $n$-equiv\-al\-ences  are marked,'' where the notion of $n$-equivalence is defined relative to the collection of marked $(n+1)$-simplices. An $n$-\emph{complicial set} is a saturated complicial set in which all simplices above dimension $n$ are marked.\footnote{One might think of the $n$-complicial sets as being a model    for the theory of $(\infty,n)$-categories, although we will not pursue that    intuition  here.} We refer to maps between complicial sets of any of the varieties just introduced that have the right lifting property with respect to the marked horn inclusions as \emph{isofibrations}.

For our purposes here, we note the following examples.

\begin{ex}[Kan complexes and quasi-categories as complicial sets]\label{ex:qcat-as-comp}
  The 0-complicial sets are precisely the Kan complexes with all positive-dimensional simplices marked and an isofibration between such is just a Kan fibration. 
  
  The 1-complicial sets are precisely the quasi-categories with their \emph{natural markings}---in which all 1-dimensional isomorphisms and all higher simplices are marked---and the isofibrations between such coincide with the isofibrations between quasi-categories. In this setting, the outer horn lifting property of 1-complicial sets and isofibrations between them is typically referred to as ``special outer horn'' lifting; see \cite[1.3]{Joyal:2002:QuasiCategories} or Proposition X.D.4.6.
\end{ex}

\begin{ex}\label{ex:nerves-as-complicial-sets}
  Suppose that $\eK$ is a quasi-category enriched category. Its homotopy coherent nerve $\hN\eK$ has:  \begin{itemize}
    \item \textbf{$0$-simplices} corresponding to the objects $a$ of $\eK$,
    \item \textbf{$1$-simplices} corresponding to $0$-arrows $f\colon a_0\to a_1$,
    \item \textbf{$2$-simplices} corresponding to diagrams
      \begin{equation*}
        \xymatrix@=1.5em{
          {a_0}\ar[rr]^{f_{02}}_{}="one"\ar[dr]_{f_{01}} && {a_2} \\
          & {a_1}\ar[ur]_{f_{12}} & \ar@{=>}^{\alpha} "one"-<0pt,0.7em>;"2,2"+<0pt,1.1em>
        }
      \end{equation*}
      where $\alpha$ is a $1$-arrow in the hom-space $\Fun_{\eK}(a_0,a_2)$ with source
      $f_{02}$ and target $f_{12}\circ f_{01}$.
    \end{itemize}
    Now we define the \textit{natural marking\/} of the homotopy coherent nerve by marking:
    \begin{enumerate}[label=(\roman*)]
    \item all $n$-simplices with $n>2$,
    \item those $2$-simplices, as depicted above, for which $\alpha$ is an invertible
      arrow in the quasi-category $\Fun_{\eK}(a_0,a_2)$, and
    \item\label{itm:1-equivalence-data} each $1$-simplex $f\colon a_0\to a_1$ which is an \emph{equivalence}, in
      the sense that it possesses an \emph{equivalence inverse} $f'\colon a_1\to
      a_0$ witnessed by a pair of invertible $1$-arrows
      \begin{equation*}
        \vcenter{\xymatrix@=1.5em{
          {a_0}\ar@{=}[rr]^{\id}_{}="one"\ar[dr]_{f} && {a_0} \\
          & {a_1}\ar[ur]_{f'} & \ar@{=>}^{\sim} "one"-<0pt,0.7em>;"2,2"+<0pt,1.1em>
        }}\mkern30mu
        \vcenter{\xymatrix@=1.5em{
          {a_1}\ar@{=}[rr]^{\id}_{}="one"\ar[dr]_{f'} && {a_1} \\
          & {a_0}\ar[ur]_{f} & \ar@{=>}^{\sim} "one"-<0pt,0.7em>;"2,2"+<0pt,1.1em>
        }}
      \end{equation*}
      in the quasi-categories $\Fun_\eK(a_0,a_0)$ and $\Fun_\eK(a_2,a_2)$
      respectively.
    \end{enumerate}
By \cite[Theorem~40]{Verity:2007rm}, the naturally marked homotopy coherent
  nerve $\qK_2\defeq\nrvhc\eK$ is a 2-complicial set. 
\end{ex}

\fi

\ifnewcut\else

\begin{defn}[marked simplicial sets]\label{defn:marked-simp-cat}
  A \textit{marked simplicial set}, called a \textit{stratified simplicial
    set\/} in~\cite{Verity:2007:wcs1} or a \textit{simplicial set with
    thinness\/} in~\cite{Street:1987:Oriental}, is a simplicial set $X$ equipped
  with a chosen subset of its simplices of dimension greater than zero called
  {\em marked simplices} which (at least) contains all of the degenerate
  simplices. Simplicial maps $f\colon X\to Y$ between marked simplicial sets are
  required to preserve the chosen markings. The category of all such marked simplicial
  sets and mark preserving simplicial maps between them is denoted
  $\msSet$.
  \end{defn}

\begin{defn}[marking conventions]\label{defn:marked-basics}
  When working with marked simplicial sets, the following conventions are commonly adopted:
  \begin{enumerate}[label=(\roman*)]
  \item When $X$ is a general simplicial set we silently promote it to a marked
    simplicial set by marking only those simplices that are degenerate. This is
    called the \emph{minimal marking} of $X$.
  \item When $\qA$ is a quasi-category we silently promote it to a marked
    simplicial set by marking all simplices above dimension $1$ and marking
    those $1$-simplices that are isomorphisms. This is called the
    \emph{natural marking}, and any simplicial map of quasi-categories
    preserves this natural marking.
  \item We say that a marked simplicial set is \textit{$n$-marked\/} when all of
    its simplices above dimension $n$ are marked. For example, under its natural
    marking a quasi-category is $1$-marked. A quasi-category is a Kan complex if
    and only if it is $0$-marked under its natural marking.
  \item\label{itm:core-adjunction} When $X$ is marked simplicial set, we let $\sharp_nX$ denote the marked
    simplicial set constructed by extending the marking of $X$ to mark all
    simplices above dimension $n$; $\sharp_0X$ is the \emph{maximal marking} of a simplicial set. 
This operator, called \textit{$n$-sharp}, is
    functorial and gives rise to an adjunction:
    \begin{equation*}
      \adjdisplay \sharp_n -| \core_n : \msSet -> \msSet .
    \end{equation*}
    When $X\in\msSet$ is a marked simplicial set, its \textit{$n$-core\/}
    $\core_n(X)$ is the subset of $X$ comprising those simplices whose
    faces above dimension $n$ are all marked.
\item The abbreviated notation $\sharp\Del^n \defeq \sharp_{n-1}\Del^n$ is used for the \emph{marked $n$-simplex}, which extends the minimal marking by also marking the top-dimensional non-degenerate simplex.
  \end{enumerate}
\end{defn}

Our primary motivation for recalling these notions is to define the natural marking of the homotopy coherent nerve of a quasi-category enriched category.

\begin{defn}\label{defn:natural-marking}
  Suppose that $\eK$ is a quasi-category enriched category. Its homotopy coherent nerve $\hN\eK$ has:  \begin{itemize}
  \item \textbf{$0$-simplices} corresponding to the objects $a$ of $\eK$,
  \item \textbf{$1$-simplices} corresponding to $0$-arrows $f\colon a_0\to a_1$,
  \item \textbf{$2$-simplices} corresponding to diagrams
    \begin{equation*}
      \xymatrix@=1.5em{
        {a_0}\ar[rr]^{f_{02}}_{}="one"\ar[dr]_{f_{01}} && {a_2} \\
        & {a_1}\ar[ur]_{f_{12}} & \ar@{=>}^{\alpha} "one"-<0pt,0.7em>;"2,2"+<0pt,1.1em>
      }
    \end{equation*}
    where $\alpha$ is a $1$-arrow in the hom-space $\Fun_{\eK}(a_0,a_2)$ with source
    $f_{02}$ and target $f_{12}\circ f_{01}$.
  \end{itemize}
  Now we define the \textit{natural marking\/} of the homotopy coherent nerve, a marked simplicial set we denoted by $\qK_2 \defeq \hN\eK$, by marking:
  \begin{enumerate}[label=(\roman*)]
  \item all $n$-simplices with $n>2$,
  \item those $2$-simplices, as depicted above, for which $\alpha$ is an invertible
    arrow in the quasi-category $\Fun_{\eK}(a_0,a_2)$, and
  \item\label{itm:1-equivalence-data} each $1$-simplex $f\colon a_0\to a_1$ which is an \emph{equivalence}, in
    the sense that it possesses an \emph{equivalence inverse} $f'\colon a_1\to
    a_0$ witnessed by a pair of invertible $1$-arrows
    \begin{equation*}
      \vcenter{\xymatrix@=1.5em{
        {a_0}\ar@{=}[rr]^{\id}_{}="one"\ar[dr]_{f} && {a_0} \\
        & {a_1}\ar[ur]_{f'} & \ar@{=>}^{\sim} "one"-<0pt,0.7em>;"2,2"+<0pt,1.1em>
      }}\mkern30mu
      \vcenter{\xymatrix@=1.5em{
        {a_1}\ar@{=}[rr]^{\id}_{}="one"\ar[dr]_{f'} && {a_1} \\
        & {a_0}\ar[ur]_{f} & \ar@{=>}^{\sim} "one"-<0pt,0.7em>;"2,2"+<0pt,1.1em>
      }}
    \end{equation*}
    in the quasi-categories $\Fun_\eK(a_0,a_0)$ and $\Fun_\eK(a_2,a_2)$
    respectively.
  \end{enumerate}
\end{defn}

  Complicial sets~\cite{Verity:2007:wcs1} are certain marked simplicial sets
  satisfying a horn filler condition generalising those that characterise Kan
  complexes and quasi-categories. To describe this notion precisely we introduce
  a couple more standard marked simplicial sets.

\begin{defn}[complicial sets]\label{defn:complicial}
 Suppose that $k$ is an integer
  in some $[n]$ then:
  \begin{itemize}
  \item The \textit{standard $k$-admissible $n$-simplex\/} $\Del^{n,k}$ is
    obtained from $\Del^n$ by marking any face that has all of the elements of
    $\{k-1,k,k+1\}\cap[n]$ among its vertices.
  \item The \textit{standard $(n,k)$-horn} $\Horn^{n,k}$ has the usual $k$-horn as
    its underlying simplicial set and inherits its markings from $\Del^{n,k}$.
  \item The marked simplicial set $\natural\Del^{n,k}$ is constructed from
    $\Del^{n,k}$ by also marking its codimension-one faces $\face^i\colon[n-1]\to[n]$
    for $i\in\{k-1,k+1\}\cap[n]$ and $\sharp\Del^{n,k}$ is constructed from
    $\natural\Del^{n,k}$ by also marking the further codimension-one face
    $\face^k\colon[n-1]\to[n]$. 
  \end{itemize}
  Using these we may make the following definitions:
  \begin{enumerate}[label=(\roman*)]
  \item A \textit{complicial set} $\qA$ is a marked simplicial set that has the
    right lifting property with respect to the horn inclusions $\Horn^{n,k} \inc
    \Del^{n,k}$ and the inclusions $\natural\Del^{n,k}\inc \sharp\Del^{n,k}$
    for $n \geq 1$ and $k\in[n]$.
  \item\label{itm:isofib-defn} A marked simplicial map $p\colon\qE\to\qB$ between complicial sets is an
    \textit{isofibration of complicial sets\/} if it has the right lifting
    property with respect to the horn inclusions $\Horn^{n,k} \inc \Del^{n,k}$ for $n \geq 1$ and $k \in [n]$.
  \item An inclusion $i\colon X\to Y$ of marked simplicial sets is said to be an
    \emph{anodyne extension} if and only if it may be constructed as the union
    of a countable chain of pushouts of coproducts of horns $\Horn^{n,k} \inc
    \Del^{n,k}$ and inclusions $\natural\Del^{n,k}\inc \sharp\Del^{n,k}$. By a
    standard argument, every isofibration of complicial sets has the right
    lifting property with respect to all anodyne extensions.
  \end{enumerate}
\end{defn}

\begin{rmk}[isofibrations of complicial sets]\label{rmk:isofibration-RLP}
  Corollary~55 of~\cite{Verity:2007:wcs1} demonstrates that a map $p\colon \qE\to \qB$  of complicial sets is an isofibration if it:
  \begin{itemize} \item
has the right lifting property with respect to all   {\em inner horns\/} $\Horn^{n,k}\inc\Del^{n,k}$, for $n \geq 2$ and $0 <k<n$
   \item 
  admits lifts of marked $1$-equivalences. More precisely, this latter condition
  asks for $p$ to have the right lifting property with respect to either of the
  inclusions $\Del^0\inc\sharp\Del^1$
  \end{itemize}
  In particular, a map between naturally marked quasi-categories is an isofibration of complicial sets if and only if it is an isofibration of quasi-categories. In the quasi-categorical context, the outer horn lifting property of Definition \ref{defn:complicial}\ref{itm:isofib-defn} is typically referred to as ``special outer horn'' lifting. 
\end{rmk}

\begin{defn}[saturated complicial sets]\label{defn:saturation}
A complicial set $\qA$ is \emph{saturated} if ``all $n$-equiv\-al\-ences in $\qA$ are marked,'' where the notion of $n$-equivalence is defined relative to the collection of marked $(n+1)$-simplices. An $n$-\emph{complicial set} is a saturated complicial set which is $n$-marked.\footnote{One might think of the $n$-complicial sets as being a model    for the theory of $(\infty,n)$-categories, although we will not pursue that    intuition  here.} 

A complicial set is saturated if and only if it satisfies a certain unique right lifting property \cite[3.2.7]{Riehl:2017cs}, so in particular any complicial set that is defined as the limit of a diagram of saturated complicial sets is again saturated. The class of saturated complicial sets is also stable under passing to $n$-cores and forming slices, in the sense of \ref{defn:join-slice} below.
\end{defn}

\begin{ex}[Kan complexes and quasi-categories as complicial sets]\label{ex:qcat-as-comp}
  Under their natural markings, Kan complexes are precisely the 0-complicial sets and
  quasi-categories are precisely the  1-complicial sets. The isofibrations between  1-complicial sets coincide with the usual classes of isofibrations between quasi-categories.
\end{ex}

\begin{ex}\label{ex:nerves-as-complicial-sets}
  Suppose that $\eK$ is a quasi-category enriched category, then Theorem~40
  of~\cite{Verity:2007rm} tells us that its naturally marked homotopy coherent
  nerve $\qK_2\defeq\nrvhc\eK$ is a 2-complicial set. 
\end{ex}

\fi

\begin{defn}[joins and slices of marked simplicial sets]\label{defn:join-slice}
The join operation extends to marked simplicial sets as follows.  Concretely, the join $X\join Y$ of two marked augmented simplicial sets $X$
  and $Y$ has as its simplices pairs $(x,y)$ with $x\in X$ and $y\in Y$ of
  arbitrary dimension with $\dim(x,y)=\dim(x)+\dim(y)+1$, where the convention is to augment a marked simplicial set with a single $-1$-simplex. We declare that a simplex $(x,y)\in X\join
  Y$ is marked if $x$ is marked in $X$ or $y$ is marked in $Y$. 
 
  Now consider a map of marked simplicial sets $f\colon X\to Y$. The slice $Y_{\sslice f}$ is
    the simplical set of whose $n$-simplices are maps $g\colon \Del^n\join X\to
    Y$ which restrict on $X\subseteq \Del^n\join X$ to the fixed map $f\colon
    X\to Y$. Such a simplex $g\colon \Del^n\join X\to Y$ is marked if and only if
    it extends along the inclusion $\Del^n\join X\subseteq
    \sharp\Del^n\join X$\ifnewcut\else,\fi\ifnewadd---where $\sharp\Del^n$ extends the minimally marked $n$-simplex $\Del^n$ by marking the non-degenerate $n$-simplex---\fi and this happens exactly when $g$ maps every simplex
    $(\id_{[n]},x)$ for $x\in X$ to a marked simplex in $Y$.  A dual construction defines $\prescript{f\sslice}{}{Y}$.

  Suppose that $\qA$ is a complicial set and that $f\colon X\to \qA$ is any map
  of marked simplicial sets. As shown in~\cite{Verity:2007:wcs1}, it is then the
  case that $\prescript{f\sslice\!\!}{}{\qA}$ and $\qA_{\sslice f}$ are also
  complicial sets and that the projections $r^f\colon
  \prescript{f\sslice\!\!}{}{\qA}\to\qA$ and $r_f\colon \qA_{\sslice f} \to\qA$ are
  isofibrations of such.
\end{defn}

  \subsection{The right adjoint to pullback}\label{sec:right-adjoint}
  
  We now have all the tools we require to construct an alternate model of the pushforward functor along a cocartesian fibration $p\colon \qE \tfib \qB$ whose value at any isofibration $q \colon \qF \tfib \qB$ will be equivalent to those of the strict pushforward. The alternate model for the pushforward 
\[ \tilde\Pi_pq \colon \tilde\Pi_p\qF \to \qB\] of an isofibration $q \colon \qF \tfib \qE$ along a cocartesian fibration $p \colon \qE \tfib \qB$ is defined as a pullback of a whiskering map for slices of homotopy coherent nerves that we now introduce.. Let $\eK$ be a quasi-category enriched category such as $\qCat$, and write $\qK_2 \defeq \hN\eK$ for its naturally marked homotopy coherent nerve, a 2-complicial set. 
  
  \begin{lem}\label{lem:lax-slice-whiskering-comp} Let $q \colon F \to E$ be a 0-arrow in a quasi-category enriched category $\eK$.
  \begin{enumerate}[label=(\roman*)]
\item\label{itm:lax-slice-whiskering-comp-i} There is a functor of slice 2-complicial sets \[ \begin{tikzcd} \laxslicer{\qK_2}{\qF} \arrow[r, "q\circ  -"] & \laxslicer{\qK_2}{\qE} \end{tikzcd}\] induced from the whiskering operation for lax cocones.
\item\label{itm:lax-slice-whiskering-comp-ii}  If $q \colon F \tfib E$ is a representably-defined isofibration, then $q\circ - \colon \laxslicer{\qK_2}{\qF} \tfib \laxslicer{\qK_2}{\qE}$ is a isofibration of complicial sets.
\end{enumerate}
\end{lem}
\begin{proof}
 By the Yoneda lemma and the natural isomorphisms arising from the slice and homotopy coherent nerve adjunctions
  \[ \sSet(X, {\qK_2}_{\sslice F}) \cong \sSet_{\top\mapsto F}(X \join\Del^0, {\qK_2}) \cong \sCat_{\top\mapsto F}(\gC[X\join\Del^0],\eK),\] to define the map in \ref{itm:lax-slice-whiskering-comp-i}, it suffices to provide a natural operation that converts a lax cocone of shape $X$ with nadir $F$ into a lax cocone with shape $X$ and nadir $E$. 
 The whiskering operation for lax cocones described in Observation \ref{obs:whiskering-cocone} defines such a natural transformation. Since whiskering preserves fibered equivalences and isomorphisms, which correspond to marked 1- and 2-simplices in $\laxslicer{\qK_2}{F}$,  this defines the desired map of 2-complicial sets.

\ifnewadd
For \ref{itm:lax-slice-whiskering-comp-ii}, we must show that the map between the sliced complicial sets has the right lifting property with respect to the marked horn inclusions of \cite[15]{Verity:2007:wcs1}, which define  lifting problems of underlying simplicial sets of the form  
\[
\begin{tikzcd} \Horn^{n,k} \arrow[d, hook] \arrow[r] & \laxslicer{\qK_2}{F} \arrow[d, "q \circ -"] \\ \Del^{n,k} \arrow[r] \arrow[ur, dashed] & \laxslicer{\qK_2}{E}
\end{tikzcd}
\]
for $n \geq 1$ and $0 \leq k \leq n$ with additional marking constraints that we describe below. 
\fi 
\ifnewcut\else 
We must show that the map between the sliced complicial sets has the right lifting property with respect to each of the marked anodyne extensions introduced in Definition \ref{defn:complicial}.
\[
\begin{tikzcd} \Horn^{n,k} \arrow[d, hook] \arrow[r] & \laxslicer{\qK_2}{F} \arrow[d, "q \circ -"] \\ \Del^{n,k} \arrow[r] \arrow[ur, dashed] & \laxslicer{\qK_2}{E}
\end{tikzcd}
\]
\fi
 By \cite[Corollary 49]{Verity:2007:wcs1},    it suffices to consider the case $0 < k \leq n$. Here the bottom horizontal functor is given by a homotopy coherent $n+1$-simplex
     \[ \gC \Del^{n+1}  \to \eK\] that sends the first $n+1$ objects to $E_0,\ldots, E_n$ and the final object to $E \in \eK$ and satisfies one additional condition forced by the markings on $\Del^{n,k}$ and ${\qK_2}_{\sslice E}$. If $k < n$, then this functor must be defined so that the 1-simplex $\alpha \in \Fun(E_{k-1},E_{k+1})$ is invertible.
         \begin{equation*}
      \xymatrix@=1.5em{
        {E_{k-1}}\ar[rr]^{f_{k-1,k+1}}_{}="one"\ar[dr]_{f_{k-1,k}} && {E_{k+1}} \\
        & {E_k}\ar[ur]_{f_{k,k+1}} & \ar@{=>}^{\alpha} "one"-<0pt,0.7em>;"2,2"+<0pt,1.1em>
      }
    \end{equation*}  If $k=n$, then the 1-simplex $\alpha \in \Fun(E_{n-2},E_n)$ must be invertible and $f_{n-1,n} \colon E_{n-1} \to E_n$ must admit an equivalence inverse.
             \begin{equation*}
      \xymatrix@=1.5em{
        {E_{n-2}}\ar[rr]^{f_{n-2,n}}_{}="one"\ar[dr]_{f_{n-2,n-1}} && {E_{n}} \\
        & {E_{n-1}}\ar[ur]_{f_{n-1,n}} & \ar@{=>}^{\alpha} "one"-<0pt,0.7em>;"2,2"+<0pt,1em>
      }
    \end{equation*}  

 The $\face^{n+1}$-face of this homotopy coherent simplex and the top horizontal together define a simplicial functor
    \[ \gC \Horn^{n+1,k} \to \eK\] that carries the $n+1$ objects to $E_0,\ldots, E_n, F$, respectively, and has the property that for each $0 \leq j \leq n$ the diagram of function complexes commutes:
    \[ 
    \xymatrix{ \Fun_{\gC \Horn^{n+1,k}}(j,n+1) \ar@{^(->}[d] \ar[r] & \Fun(E_j,F) \ar@{->>}[d]^{q \circ - } \\ \Fun_{\gC \Del^{n+1}}(j,n+1) \ar[r] & \Fun(E_j,E)}
    \] 
   
   Because $0 < k < n+1$, by \ifnewadd a calculation of homotopy coherent realizations in Example \refVI{ex:subcomputad-of-simplex}\fi\ifnewcut\else Example \ref{ex:subcomputad-of-simplex}\fi, to solve the original lifting problem, it remains only to construct a single lift
       \[ 
    \xymatrix{ \CHorn^{n,k}_1\cong \Fun_{\gC \Horn^{n+1,k}}(0,n+1) \ar@{^(->}[d] \ar[r] & \Fun(E_0,F) \ar@{->>}[d]^{q \circ - } \\\Cube^n \cong \Fun_{\gC \Del^{n+1}} (0,n+1) \ar[r] \ar@{-->}[ur] & \Fun(E_0,E)}
    \]
    the other inclusions being full. The left-hand side is a cubical horn and the right-hand side is an isofibration of quasi-categories. \ifnewcut\else By Remark \ref{rmk:isofibration-RLP},\fi\ifnewadd As noted in Example \ref{ex:qcat-as-comp}, \fi  isofibrations of quasi-categories admit lifts against ``special outer horns'' \ifnewcut\else$\Horn^{m,m} \inc\Del^{m,m}$\fi --- those in which the image of the final edge is invertible. Such extensions solve this lifting problem.
\end{proof}

\begin{prop}\label{prop:alt-right-defn}
There is a right adjoint
\[
\begin{tikzcd} \sSet_{/\qE} \arrow[r, bend right, "\tilde{\Pi}_p"'] \arrow[r, phantom, "\bot"] & \sSet_{/\qB} \arrow[l, bend right, "\tilde{p}^*"']
\end{tikzcd}
\]
to the oplax colimit functor defined at $q \colon F \to E$ by the pullback
\[
\begin{tikzcd} \tilde{\Pi}_pF \arrow[d, two heads, "\tilde{\Pi}_pq"'] \arrow[r] \arrow[dr, phantom, "\lrcorner" very near start] & \laxslicer{\qqCat_2}{\qF} \arrow[d, two heads, "q \circ -"] \\ \qB \arrow[r, "\ell^E"] & \laxslicer{\qqCat_2}{\qE}
\end{tikzcd}
\]
Moreover, when $q \colon F \tfib E$ is an isofibration, $\tilde{\Pi}_pq \colon \tilde{\Pi}_pF \tfib \qB$ is an isofibration between quasi-categories.
\end{prop}
\begin{proof}
Recall from Lemma \ref{lem:lax-slice-whiskering-comp} that $n$-simplices in $\laxslicer{\qqCat_2}{\qF}$ correspond to lax cocones under homotopy coherent simplices with nadir $F$, and observe that the whiskering functor $q\circ - \colon \laxslicer{\qqCat_2}{\qF} \to \laxslicer{\qqCat_2}{\qE}$ does not change the underlying homotopy coherent diagram. By the defining universal property, an $n$-simplex in the pullback over $b \colon \Del^n \to \qB$ corresponds to lax cocone under the homotopy coherent $n$-simplex $c_p\circ \gC{b} \colon \gC\Del^n \to \qCat$ with nadir $F$ that whiskers with $q$ to the lax cocone of \eqref{eq:lax-cocone-under-comprehension}. This recovers the characterization of the right adjoint $\tilde{\Pi}_p$ given above and Lemma \ref{lem:pullback-replacement-colimits} demonstrates that this adjoint correspondence extends to all elements of $\sSet_{/\qB}$. 

The action of $\tilde{\Pi}_p$ on morphisms $u \colon G \to F$ over $E$ is given similarly by the pullback
\[
\begin{tikzcd} \tilde{\Pi}_pG \arrow[d, dashed]  \arrow[r] \arrow[dr, phantom, "\lrcorner" very near start] & \laxslicer{\qqCat_2}{\qG} \arrow[d,  "u \circ -"]  \\ \tilde{\Pi}_pF \arrow[d, two heads, "\tilde{\Pi}_pq"'] \arrow[r]  \arrow[dr, phantom, "\lrcorner" very near start] & \laxslicer{\qqCat_2}{\qF} \arrow[d, two heads, "q \circ -"] \\ \qB \arrow[r, "\ell^E"] & \laxslicer{\qqCat_2}{\qE}
\end{tikzcd}
\]

By Lemma \ref{lem:lax-slice-whiskering-comp} and Example \ref{ex:qcat-as-comp}, it is immediate from the fact that $\qB$ is a quasi-category and $\qqCat_2$ is a  2-complicial set that $\tilde\Pi_p\qF$ is a  2-complicial set. We argue that in fact all 2-simplices are marked: by the defining universal property, a 2-simplex in $\tilde\Pi_p\qF$ corresponds to a pair comprised of a 2-simplex in $\qB$ and a 2-simplex in ${\qqCat_2}_{\sslice \qF}$ and both of these 2-simplices are marked. Thus, $\tilde\Pi_p\qF$ is a  1-complicial set, which Example \ref{ex:qcat-as-comp} tells us is the same thing as a quasi-category, and now the isofibration of complicial sets $\tilde{\Pi}_pq$ becomes an isofibration between quasi-categories.
\end{proof}

\begin{cor}\label{cor:alt-right-isofib} If $p \colon \qE \tfib \qB$ is a cocartesian fibration.:
  \begin{enumerate}[label=(\roman*)]
  \item The functor $\tilde\Pi_p \colon \sSet_{/\qE} \to \sSet_{/\qB}$ carries isofibrations over $\qE$ to isofibrations over $\qB$, restricting to define a functor
  \begin{equation}\label{eq:alt-pi-as-slice-functor} \qCat_{/\qE} \xrightarrow{\tilde\Pi_p} \qCat_{/\qB}.\end{equation}
  \item The functor \eqref{eq:alt-pi-as-slice-functor} preserves isofibrations, now considered as morphisms in these slice categories.
  \end{enumerate}
  \end{cor}
\begin{proof}
Proposition \ref{prop:alt-right-defn} demonstrates that $\tilde{\Pi}_p$ carries isofibrations to isofibrations, restricting to define a functor $\tilde{\Pi}_p \colon \qCat_{/\qE} \to \qCat_{/\qB}$. Moreover this functor preserves isofibrations, now considered as morphisms in these slice categories, since the action of $\tilde{\Pi}_p$ on an isofibration $u \colon \qG \tfib \qF$ over $\qE$ is defined by pulling back the isofibration of complicial sets $u \circ - \colon \laxslicer{\qqCat_2}{\qG} \tfib \laxslicer{\qqCat_2}{\qF}$. 
\end{proof}

We now transfer the properties of the functor $\tilde{\Pi}_p$ to the right adjoint $\Pi_p \colon \sSet_{/\qE} \to \sSet_{/\qB}$ to the strict pullback functor $p^* \colon \sSet_{/\qB} \to \sSet_{/\qE}$.

\begin{prop}\label{prop:quillen-exponentiation-adjunctions} If $p \colon E \tfib \qB$ is a cocartesian fibration, then the adjunctions
\[
\begin{tikzcd} \sSet_{/\qE} \arrow[r, bend right, "{\Pi}_p"'] \arrow[r, phantom, "\bot"] & \sSet_{/\qB} \arrow[l, bend right, "{p}^*"'] & \text{and} &  \sSet_{/\qE} \arrow[r, bend right, "\tilde{\Pi}_p"'] \arrow[r, phantom, "\bot"] & \sSet_{/\qB} \arrow[l, bend right, "\tilde{p}^*"']
\end{tikzcd}
\] are Quillen with respect to the sliced Joyal model structure. 
\end{prop}

In particular, $\Pi_p$ preserves both fibrant objects and the fibrations between, and thus has the properties enumerated for $\tilde\Pi_p$ in Corollary \ref{cor:alt-right-isofib}. Consequently, the natural Joyal equivalence $\gamma \colon \tilde{p}^* \To p$ of Theorem \ref{thm:pullback-equivalence}, which defines a natural isomorphism of total left derived functors, transposes to a natural equivalence $\hat{\gamma} \colon \Pi_p  \To \tilde{\Pi_p}$, which defines a natural isomorphism of total right derived functors.

\begin{proof}
By an observation of Joyal and Tierney \cite[7.15]{Joyal:2007kk}, to show that $\tilde{p}^*\dashv \tilde{\Pi}_p$ is Quillen it suffices to show that the left adjoint preserves cofibrations and the right adjoint preserves fibrations between fibrant objects. Lemma \ref{lem:pullback-replacement-monos} demonstrates the first of these and Corollary \ref{cor:alt-right-isofib}(ii) proves the second.

To prove that $p^* \dashv \Pi_p$ is Quillen, we prove that $p^*$ is left Quillen. Thus functor preserves cofibrations because pullbacks preserve monomorphisms. By Theorem \ref{thm:pullback-equivalence}, $p^*$ is naturally weakly equivalent to the left Quillen functor $\tilde{p}^*$. Since all objects in $\sSet_{/\qE}$ are cofibrant, the left Quillen functor $\tilde{p}^*$ preserves all Joyal weak equivalences, and hence by the 2-of-3 property $p^*$ does as well.
\end{proof}

We now consider the actions of the pushforward functors $\tilde{\Pi}_p$ and $\Pi_p$ along a cocartesian fibration $p \colon \qE \tfib \qB$ when applied to a cartesian fibration $q \colon \qF \tfib \qE$. As before, we demonstrate directly that $\tilde{\Pi}_pq \colon \tilde{\Pi}_p\qF \tfib \qB$ is a then a cartesian fibration and then use Theorem \ref{thm:pullback-equivalence} to conclude the same for $\Pi_p$.

\begin{lem}\label{lem:lax-slice-whiskering-cart} Let $q \colon \qF \tfib \qE$ between quasi-categories.  Then the corresponding map $q \circ - \colon \laxslicer{\qqCat_2}{\qF} \tfib \laxslicer{\qqCat_2}{\qE}$ of 2-complicial sets has the right lifting property with respect to any outer horn inclusion
\[
\begin{tikzcd} \Del^{\fbv{n-1,n}} \arrow[r] \arrow[rr, bend left=10, "\chi"] & \Horn^{n,n} \arrow[r] \arrow[d, hook] & \laxslicer{\qqCat_2}{\qF} \arrow[d, two heads, "q\circ-"] \\&  {\Del^n} \arrow[ur, dashed] \arrow[r] & \laxslicer{\qqCat_2}{\qE}
\end{tikzcd}
\]
whose final edge defines a cartesian 1-arrow
\[
\begin{tikzcd}[sep=small] E_{n-1} \arrow[rr, "f_{n-1}"] \arrow[dr, "f_{n-1,n-2}"'] & \arrow[d, phantom, "\scriptstyle\Downarrow\chi"] & F  \\ & E_{n} \arrow[ur, "f_{n}"'] 
\end{tikzcd}
\]
for the cartesian fibration $q \circ - \colon \Fun(E_{n-1},F) \to \Fun(E_{n-1},E)$.
\end{lem}
\begin{proof}
As in the proof of Lemma \ref{lem:lax-slice-whiskering-comp}, the bottom horizontal functor is given by a homotopy coherent $n+1$-simplex $\gC\Del^{n+1} \to \qCat$ that sends the first $n+1$ objects to $E_0, \ldots E_n$ and the final object to $E$, while the $\face^{n+1}$-face of this homotopy coherent simplex and the top horizontal functor together define a simplicial functor $\gC\Horn^{n+1,n} \to \qCat$ that caries the $n+2$ objects to $E_0,\ldots, E_n,F$ and has the property that for each $0 \leq j \leq n$ the diagram of function complexes commutes:
\[
\begin{tikzcd}
\gC\Horn^{n+1,n}(j,n+1) \arrow[r] \arrow[d, hook] & \Fun(E_j,E) \arrow[d, two heads, "q \circ -"] \\ \gC\Del^{n+1}(j,n+1) \arrow[r] & \Fun(E_j,F)
\end{tikzcd}
\]

By Example
\ifnewadd \refVI{ex:subcomputad-of-simplex}\fi\ifnewcut\else \ref{ex:subcomputad-of-simplex}\fi, to solve the original lifting problem, we need only construct a single lift
\[
\begin{tikzcd}
\sqcap^{n,n}_1 \cong \gC\Horn^{n+1,n}(0,n+1) \arrow[r] \arrow[d, hook] & \Fun(E_0,E) \arrow[d, two heads, "q \circ -"] \\ \square^n\cong \gC\Del^{n+1}(0,n+1) \arrow[r] \arrow[ur,dashed] & \Fun(E_0,F)
\end{tikzcd}
\]
This extension problem can be solved by filling inner horns and ``special outer horns''
 $\Horn^{m,m} \to \Del^m$, those whose final edges are complies of the 1-simplex $\chi \in \Fun(E_{n-1},F)$ pre-composed with some functor $E_0 \to E_{n-1}$. Such 1-simplices represent $(q\circ -)$-cartesian cells so these ``special outer horn'' lifting problems also admit solutions by Lemma \ifnewadd\refVI{lem:cocart-cylinder-extensions}\fi\ifnewcut\else\ref{lem:cocart-cylinder-extensions}(i)\fi.
 \end{proof}

\begin{prop}\label{prop:alt-right-cart-fun} If $p \colon \qE \tfib \qB$ is a cocartesian fibration and $q \colon \qF \tfib \qE$ is a cartesian fibration between quasi-categories, then
\[ \tilde{\Pi}_pq \colon \tilde{\Pi}_p\qF \tfib \qB
\]
is a cartesian fibration between quasi-categories. Moreover, $\tilde{\Pi}_p$ preserves cartesian functors, restricting to define a functor
\[ \tilde{\Pi}_p \colon \Cart(\qCat)_{/\qE} \to \Cart(\qCat)_{/\qB}.\]
\end{prop}
\begin{proof}
  By Proposition \ref{prop:alt-right-defn}, $ \tilde\Pi_p\qF \tfib \qB$ defines an isofibration between quasi-categories. Lemma \ref{lem:lax-slice-whiskering-cart}  identifies a class of cartesian 1-arrows in $\tilde\Pi_p\qF$ in the sense of Definition \ref{defn:qcat-cocart}, which we now describe explicitly.  Recall from the construction of Proposition \ref{prop:alt-right-defn}, that a 1-simplex $\chi \colon \Del^1 \to \tilde{\Pi}_p\qF$ in the fiber over $b \colon \Del^1 \to \qB$ corresponds to a 1-arrow 
\[
\begin{tikzcd}[sep=small] E_{0} \arrow[rr, "f_{0}"] \arrow[dr, "e_b"'] & \arrow[d, phantom, "\scriptstyle\Downarrow\chi"] & \qF  \\ & E_{1} \arrow[ur, "f_{1}"'] 
\end{tikzcd}
\]
in $\Fun(E_0,\qF)$ that whiskers with $q \colon \qF \to E$ to define the lax cocone that restricts the lax cocone associated to the comprehension construction along $b$. As observed previously, this compatibility condition tells us that the quasi-categories $E_0$ and $E_1$ are the fibers of $p \colon E \tfib \qB$ over the vertices in $b$ and the functor $e_b \colon E_0 \to E_1$ is the comprehension of $b$. To form such a lift with codomain $f_1 \colon E_1 \to \qF$, start by lifting $b$ to the lax cocone 
\[
\begin{tikzcd}[sep=small] E_{0} \arrow[rr, "\ell^E_{0}"] \arrow[dr, "e_b"'] & \arrow[d, phantom, "\scriptstyle\Downarrow\epsilon"] & E  \\ & E_{1} \arrow[ur, "\ell^E_{1}"'] 
\end{tikzcd}
\]
under $e_b \colon E_0 \to E_1$ with nadir $E$ associated with the comprehension construction, as displayed in \eqref{eq:comprehension-over-arrow}. Since $f_1$ is in the fiber of $\tilde{\Pi}_pq \colon \tilde{\Pi}_p\qF \tfib \qB$ over the codomain of $b$, we must have $q f_1 = \ell^E_1$. Now we can lift $\epsilon$ along the cartesian fibration $q$ to a $q$-cartesian cell with codomain $f \circ e_b$. This defines the $(q\circ -)$-cartesian cell $\chi$.

Now if $u \colon G \to \qF$ is a cartesian functor from $r \colon G \tfib E$ to $q \colon \qF \tfib E$, then $u$ is representabily cartesian in the sense that $u \circ - \colon \Fun(X,G) \to \Fun(X,\qF)$ carries $(r\circ -)$-cartesian 1-arrows to $(q\circ -)$-cartesian 1-arrows. Since $u \circ - \colon \laxslicer{\qqCat_2}{\qG} \to \laxslicer{\qqCat_2}{\qF}$ preserves the cartesian 1-arrows just identified, proving that $\tilde{\Pi}_p$ carries this map to a cartesian functor between cartesian fibrations over $\qB$.\end{proof}

\begin{cor}\label{cor:cartesian-workhorse}  If $p \colon \qE \tfib \qB$ is a cocartesian fibration and $q \colon \qF \tfib E$ is a cartesian fibration between quasi-categories, then
\[ {\Pi}_pq \colon {\Pi}_p\qF \tfib \qB
\]
is a cartesian fibration between quasi-categories. Moreover, ${\Pi}_p$ preserves cartesian functors, restricting to define a functor
\[ \begin{tikzcd} \Cart(\qCat)_{/\qE} \arrow[r, "\Pi_p"] & \Cart(\qCat)_{/\qB}. \end{tikzcd} \]
\end{cor}
\begin{proof} By Proposition \ref{prop:quillen-exponentiation-adjunctions}, the components
\[
\begin{tikzcd}[sep=small] \Pi_p\qF \arrow[rr, "\hat{\gamma}_q", "\sim"'] \arrow[dr, two heads, "\Pi_pq"'] & & \tilde{\Pi}_p\qF \arrow[dl, two heads, "\tilde{\Pi}_pq"] \\ & \qB
\end{tikzcd}
\]
at an isofibration $q \colon \qF \tfib E$ of the transpose $\hat\gamma \colon \Pi_p \To \tilde{\Pi}_p$ of the natural weak equivalence of Theorem \ref{thm:pullback-equivalence} are equivalences of isofibrations over $\qB$. If $q$ is a cartesian fibration, then Proposition \ref{prop:alt-right-cart-fun} proves that $\tilde{\Pi}_pq$ is a cartesian fibration, and since the notion of cartesian fibration is equivalence-invariant, $\Pi_pq$ must be as well. 
\end{proof}

\begin{thm}\label{thm:cosmological-pushforward} For a cocartesian fibration $p \colon E \tfib \qB$ between quasi-categories, the pullback-pushforward adjunction restricts to define an adjunction
\[
\begin{tikzcd}[row sep=large] \qCat_{/\qE} \arrow[r, bend right=20, "\Pi_p"'] \arrow[r, phantom, "\bot"] & \qCat_{/\qB} \arrow[l, bend right=20, "p^*"'] \\  \Cart(\qCat)_{/\qE} \arrow[r, bend right=20, "\Pi_p"'] \arrow[r, phantom, "\bot"] \arrow[u, hook] & \Cart(\qCat)_{/\qB} \arrow[l, bend right=20, "p^*"'] \arrow[u, hook]
\end{tikzcd}
\]
\end{thm}
\begin{proof}
By Proposition \ref{prop:quillen-exponentiation-adjunctions}, the adjoint functors $p^* \dashv \Pi_p$ define an adjunction 
\[
\begin{tikzcd}[row sep=large] \qCat_{/\qE} \arrow[r, bend right=20, "\Pi_p"'] \arrow[r, phantom, "\bot"] & \qCat_{/\qB} \arrow[l, bend right=20, "p^*"'] 
\end{tikzcd}
\]
By Proposition \ref{prop:pullback-stability}, the left adjoint restricts to a define a functor \[p^* \colon \Cart(\qCat)_{/\qB} \to \Cat(\qCat)_{/\qE}.\] By Corollary \ref{cor:cartesian-workhorse}, the right adjoint also restricts to a functor $\Pi_p \colon \Cart(\qCat)_{/\qE} \to \Cart(\qCat)_{/\qB}$. Since the inclusion $\Cart(\qCat)_{/\qB} \inc \qCat_{/\qB}$ is not full, this is not quite enough to demonstrate adjointness of the restricted adjunction: it remains to argue that the adjoint transpose of a cartesian functor is a cartesian functor.

To that end, let $q \colon \qF \tfib E$ and $r \colon G \tfib \qB$ be cartesian fibrations. A functor $f \colon G \to \tilde{\Pi}_pq$ over $\qB$ is cartesian if and only if the square
\[
\begin{tikzcd} G  \arrow[r, "f"] \arrow[d, two heads, "r"'] & \laxslicer{\qqCat_2}{\qF}  \arrow[d, two heads, "{q \circ -}"] \\ \qB \arrow[r, "\ell^E"'] & \laxslicer{\qqCat_2}{\qE}
\end{tikzcd}
\]
carries $r$-cartesian arrows to representably $q$-cartesian arrows in $\laxslicer{\qqCat_2}{\qF}$, as described in Lemma \ref{lem:lax-slice-whiskering-cart}. Fixing an 1-arrow arrow $\zeta \colon \Del^1 \to G$ over $b \colon \Del^1 \to \qB$ as below-left, the arrow $f\zeta$ transposes to the functor over $E$ displayed below-right
\[
\begin{tikzcd}[column sep={5em,between origins}] \Del^1  \arrow[r, "f\zeta"] \arrow[d, "b"'] & \laxslicer{\qqCat_2}{\qF}  \arrow[d, two heads, "{q \circ -}"] & ~&  \colim\nolimits^{\mathrm{oplax}}{(c_p\circ\gC{b})} \arrow[rr, "\widehat{f\zeta}"] \arrow[dr, "\ell^E\vert_b"'] & & \qF \arrow[dl, two heads, "q"] \\ \qB \arrow[r, "\ell^E"'] & \laxslicer{\qqCat_2}{\qE} &\arrow[u, phantom, "\leftrightsquigarrow"] & & E
\end{tikzcd}
\]
By Proposition \ref{prop:pullback-equivalence-1}, the oplax colimit is equivalent to the fiber $E_b$ and the functor $\widehat{f\zeta}$ represents the whiskered lax cocone 
\[
\begin{tikzcd}[row sep=tiny] 
 E_{b_0} \arrow[rrrd, bend left=10, ""{name=1,right}, "\ell^E_0"] \arrow[ddr, dotted, "e_b"'] \arrow[ddd, "p_0"'] &&& & \qF \arrow[dr, two heads, "q"]  \\  & & & E_b \arrow[ddd, two heads, "p_b"] \arrow[r] \arrow[dddrr, phantom, "\lrcorner" very near start] \arrow[ur, "\widehat{f\zeta}"] & E_r \arrow[r] \arrow[u, "\hat{f}"']  \arrow[dddr, phantom, "\lrcorner" very near start] \arrow[ddd, two heads] & E \arrow[ddd, two heads, "p"]  \\ 
& E_{b_1} \arrow[urr, bend right=10, "\ell^E_1"'] \arrow[from=1,Rightarrow, dotted, shorten >= 0.5em, shorten <= 0.5em, "\chi"]\arrow[ddrr, phantom, "\lrcorner" pos=.05]  \\
1 \arrow[rrrd, bend left=10, ""{name=2,right}, "0"] \arrow[ddr, equals]  \\ & & & \Del^1 \arrow[rr, "b"', bend right] \arrow[r, "\zeta"] & G \arrow[r, two heads, "r"] & \qB\\ 
& 1 \arrow[urr, bend right=10, "1"'] \arrow[from=uuu, crossing over, two heads, "p_1"' pos=.6] \arrow[from=2,Rightarrow,  shorten >= 0.5em, shorten <= 0.5em, "\kappa"] 
\end{tikzcd}
\]
Now $f$ is a cartesian functor if and only if the whiskered composite $\widehat{f\zeta}\chi$ is $q$-cartesian whenever $\zeta$ is $r$-cartesian.  Since Proposition \ref{prop:pullback-stability} demonstrates that cartesian arrows are created by pullbacks, this proves that $f$ is a cartesian functor if and only if the transposed functor is cartesian:
\[
\begin{tikzcd}[sep=small]  E_r \arrow[rr, "\hat{f}"] \arrow[dr, two heads, "p^*r"'] & & \qF \arrow[dl, two heads, "q"] \\ & E
\end{tikzcd} \qedhere
\]
\end{proof}

For use in \ifnewcut\else\cite[\S12.3]{RiehlVerity:2022eo},\fi \ifnewadd sequels to this work,\fi ~we note that the pushforward is a cosmological functor between the $\infty$-cosmoi established in Proposition \refVIII{prop:cartesian-cosmoi}.

\begin{cor}\label{cor:cosmological-pushforward} Let $p \colon \qE \to \qB$ be a cocartesian fibration. Then the pushforward construction defines cosmological functors
\[ \Pi_p \colon \qCat_{/\qE} \to \qCat_{/\qB} \qquad \mathrm{and} \qquad \Pi_p \colon \Cart(\qCat)_{/\qE} \to \Cart(\qCat)_{/\qB},\]
which is to say that they are simplically enriched, preserve all simplicially enriched limits with flexible weights, and preserve the isofibrations, considered as morphisms in the slice category.  
\end{cor}
\begin{proof}
The functor $\Pi_p \colon \qCat_{/\qE} \to \qCat_{/\qB}$ is the restriction of a Quillen right adjoint $\Pi_p \colon \sSet_{/\qE} \to \sSet_{/\qB}$. To prove that this defines a cosmological, it remains only to show that the adjunction $p^* \dashv \Pi_p$ is simplicially enriched. This follows from Observation \ref{obs:tensor-pullback}, which notes that the left adjoint preserves tensors with simplicial sets. Lemma \ref{lem:cart-functor-tensor} observes that the simplicial enrichment descends to the subcosmoi of cartesian fibrations.
\end{proof}

The argument given in the proof of Theorem \ref{thm:cosmological-pushforward} provides a characterization of the $\Pi_pq$-cartesian 1-arrows in the cartesian fibration constructed from a cocartesian fibration $p \colon E \tfib \qB$ and a cartesian fibration $q \colon \qF \tfib E$ between quasi-categories that lift a specified arrow $\beta \colon \Del^1 \to \qB$. 

\begin{lem}\label{lem:cart-arrow-in-pushforward}
If $p \colon \qE \tfib \qB$ is a cocartesian fibration and $q \colon \qF \tfib \qE$ is a cartesian fibration between quasi-categories, then the cartesian 1-arrows $\chi$ in $\Pi_pq \colon \Pi_p\qF \tfib \qB$ are those maps that transpose to define functors  that carry $p$-cocartesian lifts of $\beta$ to $q$-cartesian lifts.
\[
\begin{tikzcd}  \Del^1 \arrow[dr, "\beta"'] \arrow[rr, "\chi"] && \Pi_p\qF  \arrow[dl, two heads, "\Pi_pq"] \\ & \qB
\end{tikzcd}
\qquad\leftrightsquigarrow\qquad \begin{tikzcd} & \qF \arrow[d, two heads, "q"] \\ \qE_\beta \arrow[ur, "\hat{\chi}"] \arrow[r, "\ell_\beta"] \arrow[d, two heads, "p_\beta"'] \arrow[dr, phantom, "\lrcorner" very near start] & \qE \arrow[d, two heads, "p"] \\ \Del^1 \arrow[r, "\beta"'] & \qB
\end{tikzcd}
\]
\end{lem}
\begin{proof}
By Theorem \ref{thm:pullback-equivalence}, $\qE_\beta$ may be identified with the oplax colimit of the canonical lax cocone formed by taking a $p$-cocartesian lift of $\beta$. From this perspective, the transposed functor $\hat\chi \colon \qE_\beta \to \qF$ acts by whiskering this $p$-cocartesian arrow. By the construction in the proof of Theorem \ref{thm:cosmological-pushforward}, the $\Pi_pq$-cartesian lifts of $\beta$ are those arrows for which this whiskered composite is $q$-cartesian, as claimed.
\end{proof}

\subsection{Exponentiation}\label{sec:exponentiation}

As is familiar in any locally cartesian closed category, the pullback and pushforward functors can be used to  construct \emph{exponentials} in $\Cart(\qCat)_{/\qB}$, where the exponent is given by a cocartesian fibration.

\begin{defn}[exponentials]\label{defn:exponential}
For $p \colon \qE \tfib \qB$ either a cartesian or cocartesian fibration and  $q \colon \qF \tfib \qB$ an isofibration, define 
\begin{equation}\label{eq:exponential} (q \colon \qF \tfib \qB)^{p \colon \qE \tfib \qB}\qquad \in \qCat_{/\qB}\end{equation} to be the image of $q$ under the composite functor
\[\xymatrix{\qCat_{/\qB} \ar[r]^-{p^*} &\qCat_{/\qE} \ar[r]^-{\Pi_p} & \qCat_{/\qB}}\]
 
Note that by the adjunctions $\Sigma_p \dashv p^* \dashv \Pi_p$, if $r \colon \qG \tfib \qB$ is also a cartesian or cocartesian fibration, there are natural isomorphisms
\[ \Fun_\qB(\qG \tfib \qB,(\qF \tfib \qB)^{\qE \tfib \qB}) \cong \Fun_{\qB}( \qE \times_{\qB} \qG \tfib \qB, \qF \tfib \qB) \cong \Fun_\qB(\qE \tfib \qB,(\qF \tfib \qB)^{\qG \tfib \qB}), \] and the left-hand isomorphism still holds in the case where $r \colon \qG \to\qB$ is a mere functor, whose domain need not even be a quasi-category.
\end{defn}

\begin{prop}\label{prop:cartesian-workhorse} 
 If $p \colon \qE \tfib \qB$ is a cocartesian fibration and $q \colon \qF \tfib \qB$ is a cartesian fibration, then  \eqref{eq:exponential} is a cartesian fibration whose cartesian 1-arrows are those maps that transpose to define functors that carry $p$-cocartesian lifts of $b$ to $q$-cartesian lifts of $b$. 
\[
\xymatrix{ \Del^1 \ar[dr]_b \ar[rr]^-\chi & & (\qF \xtfib{q} \qB)^{\qE\xtfib{p}\qB} \ar@{->>}[dl]  \ar@{}[dr]|{\displaystyle\leftrightsquigarrow} & 
\qE_b \ar@{->>}[d]_{p_b} \ar[r] \pbexcursion \ar@/^2ex/[rr]^{\widehat{\chi}} & \qE \ar@{->>}[d]^p  & \qF \ar@{->>}[dl]^q  \ar@{}[dr]|{\displaystyle\leftrightsquigarrow} & \qE_b \ar@{->>}[dr]_{p_b} \ar[rr]^{\widehat{\chi}} & & \qF_b \ar@{->>}[dl]^{q_b} \\ & \qB & & \Del^1 \ar[r]_b & \qB & &  & \Del^1}\] 
\end{prop}
\begin{proof} The first statement follows from Corollary \ref{cor:cartesian-workhorse}, while the characterization of cartesian cells is given in Lemma \ref{lem:cart-arrow-in-pushforward}.
\end{proof}

Recall the function complexes constructed in Definitions \ref{defn:fun-over-B} and \ref{defn:cart-fun-over-B}.

\begin{lem}\label{lem:cartesian-functor-transpose} Let $q \colon \qF \tfib \qB$ be a cartesian fibration, let $p \colon \qE \tfib \qB$ be a cocartesian fibration, and let $\pi \colon \qA \times \qB \tfib \qB$ denote the projection, a bifibration. Then the isomorphism
\[ \Fun_{\qB}(\qF \xtfib{q}\qB, (\qA \times \qB \tfib \qB)^{\qE \xtfib{p} \qB}) \cong  \Fun_{\qB}(\qE \xtfib{p}\qB, (\qA \times \qB \tfib \qB)^{\qF \xtfib{q} \qB})\] restricts to an isomorphism
\[ \Fun_{\qB}^c(\qF \xtfib{q}\qB, (\qA \times \qB \tfib \qB)^{\qE \xtfib{p} \qB}) \cong  \Fun_{\qB}^c(\qE \xtfib{p}\qB, (\qA \times \qB \tfib \qB)^{\qF \xtfib{q} \qB})\] between the function complexes in the quasi-category enriched categories $\Cart(\qCat)_{/\qB}$ and $\coCart(\qCat)_{/\qB}$, which is to say,  cartesian functors between the cartesian fibrations on the left transpose to define cartesian functors between the cocartesian fibrations on the right.
\end{lem}
\begin{proof}
By adjunction, the data of a map over $\qB$ 
\[
\begin{tikzcd}[column sep={4em,between origins}] \qF \arrow[dr, two heads, "q"'] \arrow[rr, "f"] & & (\qA \times \qB \tfib \qB)^{\qE \xtfib{p} \qB} \arrow[dl, "\pi^p"] \\ & \qB
\end{tikzcd}
\]
is given by a single functor $\hat{f} \colon \qF \times_\qB \qE \to \qA$. By Proposition \ref{prop:cartesian-workhorse}, $f$ is a cartesian functor if and only if for each $q$-cocartesian cell $\chi \colon \Del^1 \to \qF$ over $b \colon \Del^1 \to \qB$, the induced functor
\[
\begin{tikzcd} \qE_b \arrow[d, two heads, "p_b"'] \arrow[r, "\hat{\chi}"] & \qA \times \qB \arrow[d, two heads, "\pi"] \\ \Del^1 \arrow[r, "b"'] & \qB
\end{tikzcd}
\]
carries $p$-cartesian cells $\gamma \colon \Del^1 \to \qE$ over $b$ to $\pi$-cocartesian ones, these being those maps $\Del^1 \to \qA \times \qB$ whose component along the other projection $\qA \times \qB \to \qA$ is invertible.  In summary, the functor $f$ is cartesian if and only if for every cocartesian lift $\chi$ and cartesian lift $\gamma$ of $b$, the composite morphism
\[
\begin{tikzcd} \Del^1 \arrow[r, "{\langle \chi,\gamma \rangle}"] & \qF_b \times_{\Del^1} \qE_b \arrow[r, hook] & \qF \times_\qB\qE \arrow[r, "\hat{f}"] & \qA
\end{tikzcd} \qedhere
\]\end{proof}


%% file: comonadicity.tex

\section{Monadicity and comonadicity of cartesian fibrations}\label{sec:monadicity}

The 0-skeleton of a quasi-category $\qB$ defines the ``underlying set of objects'' $\ob{\qB}$, together with a canonical inclusion $\ob{\qB}\inc \qB$.  By Proposition \ref{prop:pullback-stability}, pulling back along the inclusion $\ob\qB\inc\qB$ induces a forgetful functor
\[ \xymatrix@R=0.5em{\Cart(\qCat)_{/\qB} \ar[r] & \Cart(\qCat)_{/\ob\qB}\cong \prod\limits_{\ob\qB} \qCat \\ (\qE \xtfib{p} \qB) \ar@{}[r]|-\mapsto & (\qE_b)_{b \in \ob\qB}} \]
whose codomain is isomorphic to the product of the quasi-categorically enriched categories of quasi-categories, a cartesian fibration over a set being simply an indexed family of 
quasi-categories. Our aim in this section is to construct left and right adjoints and prove that this functor is monadic and comonadic in a suitable sense.

The adjoint functors are constructed as what we refer to as \emph{biadjoint} functors of quasi-categorically enriched categories: that is, we construct quasi-categorically enriched functors
\[ L,R \colon \Cart(\qCat)_{/\ob\qB} \longrightarrow \Cart(\qCat)_{/\qB}\] together with natural equivalences of function complexes  that encode the adjoint transpose correspondence. The right adjoint makes use of the exponentiation construction of \S\ref{sec:exponentiation} and the Yoneda lemma is used to prove biadjointness. These tasks occupy  \S\ref{sec:adjoint-construction}.

Any quasi-categorically enriched category $\eK$ has a (typically large) \emph{quasi-categorical core} $\qK \defeq \hN g_*\eK$ defined by passing to the maximal Kan complex enriched core and then applying the homotopy coherent nerve. For example, the quasi-category $\qqCat$ of quasi-categories and functors is the quasi-categorical core of $\qCat$. In  \S\ref{sec:adjoint-construction} we also prove that biadjoint functors of quasi-categorically enriched categories descend to adjoint functors between their quasi-categorical cores.

In particular, this implies that the map of large quasi-categories of cartesian fibrations and cartesian functors
\[ \Cart_{/\qB} \longrightarrow \Cart_{/\ob\qB} \cong \prod\limits_{\ob\qB}\qqCat\] admits both left and right adjoints.  In \S\ref{sec:comonadicity}, we prove first that this forgetful functor is comonadic and then use comonadicity to prove that it is also monadic. To do so, we appeal to the comonadicity theorem proven in \S\refII{sec:monadicity} and recalled as Theorem \ref{thm:comonadicity} below. The monadicity of this forgetful functor will be used in \S\ref{sec:groupoidal-reflection} to construct a ``groupoidal reflection'' functor for cartesian fibrations.

\subsection{Adjoint functors}\label{sec:adjoint-construction}

A functor $U \colon \eK \to \eL$ between quasi-categorically enriched categories gives rise to a functor between the large quasi-categories defined by passing to the Kan complex enriched cores of $\eK$ and $\eL$ and applying the homotopy coherent nerve construction. We frequently find it convenient to construct adjoints to this functor of quasi-categories ``at the point-set level'' by producing the structures axiomatized in the following definition:

\begin{defn} A \emph{biadjunction of quasi-categorically enriched categories} consists of:
\begin{itemize}
\item a pair of quasi-categorically enriched categories $\eK$ and $\eL$;
\item a pair of simplicial functors $F \colon \eL \to \eK$ and $U \colon \eK \to \eL$; and
\item a simplicially-enriched natural equivalence
\[ \Fun_{\eK}(FL,K) \simeq \Fun_{\eL}(L,UK)\] of function complexes .
\end{itemize}
\end{defn}

\begin{prop} If $F \colon \eL \to \eK$ and $U \colon \eK \to \eL$
    define a biadjunction of quasi-categorically enriched categories, then the induced functors between the quasi-categorical cores $\qK \defeq \hN g_*\eK$ and $\qL \defeq \hN g_*\eL$ define an adjunction:
\[ \adjdisplay F -| U : \qK -> \qL.\] 
\end{prop}
\begin{proof}
Reprising a construction from the proof of Theorem \refI{thm:simplicial-Quillen-adjunction}, we define simplicial categories $\coll(F, \eK)$ and $\coll(\eL,U)$ whose objects are $\ob\eK +\ob\eL$ and which include $\eK$ and $\eL$ as full subcategories. Each of the function complexes from an object of $\eK$ to an object of $\eL$ are empty, while for $L \in \eL$ and $K \in \eK$, we define
\[ \Fun_{\coll(F,\eK)}(L,K) \defeq \Fun_{\eK}(FL,K) \qquad \mathrm{and} \qquad \Fun_{\coll(\eL,U)}(L,K) \defeq \Fun_{\eL}(L,UK).\] The natural equivalence $\Fun_{\eK}(FL,K) \we \Fun_{\eL}(L,UK)$ of the biadjunction gives rise to a simplicial functor $\coll(F,\eK) \to \coll(\eL,U)$ under $\eK + \eL$ that is bijective on objects and a local equivalence of quasi-categories. 

Passing to Kan complex enriched cores, this functor defines a Dwyer-Kan equivalence, and thus yields an equivalence of quasi-categories upon passing to homotopy coherent nerves. Note that the homotopy coherent nerve of the groupoid core of $\coll(F,\eK)$ is isomorphic to the quasi-categorical collage of the underlying functor $F \colon \qL \to \qK$ constructed in Definition \ref{defn:qcat-collage}. Now Corollary \ref{cor:lurie-adjunction} demonstrates that $F \dashv U$ as functors between the quasi-categories $\qK$ and $\qL$.
\end{proof}

\begin{prop}\label{prop:right-adjoint-construction}
The functor $\Cart(\qCat)_{/\qB} \longrightarrow \Cart(\qCat)_{/\ob\qB}$ admits a quasi-cat\-e\-gor\-ically enriched right biadjoint defined by
\[ 
\xymatrix@R=1em{ \Cart(\qCat)_{/\ob\qB} \ar[r] & \Cart(\qCat)_{/\qB} \\  (\qE_b)_{b \in\ob\qB} \ar@{}[r]|-\mapsto & \prod_{b} (\qE_b \times \qB \tfib \qB)^{ b \comma \qB \tfib \qB} }
\]
that is, an $\ob\qB$-indexed family of quasi-categories $(\qE_b)_{b\in\ob\qB}$ is sent to the product in $\Cart(\qCat)_{/\qB}$ of the cartesian fibrations $(\qE_b \times \qB \tfib \qB)^{b \comma \qB \tfib \qB}$.
\end{prop}

\begin{proof}
    \ifnewadd Here $b\comma \qB \tfib \qB$ is the (groupoidal) cocartesian fibration represented by the vertex $b \in \qB$, sending an arrow in $\qB$ with domain $b$ to its codomain; see Example \refIV{ex:comma-fibrations}. The fiber over a vertex $x \in \qB$, is the Kan complex $b \comma x$ of maps from $b$ to $x$ in $\qB$. \fi  
    \ifnewcut\else
Recall from Example \ref{ex:comma-qcat} that $b\comma \qB \tfib \qB$ is the cocartesian fibration represented by the vertex $b \in \qB$. \fi Since the product projection $\pi \colon \qE_b \times \qB \tfib \qB$ is a bifibration, \ifnewadd both cartesian and cocartesian, \fi Proposition \ref{prop:cartesian-workhorse} implies that  $(\qE_b \times \qB \tfib \qB)^{b \comma \qB \tfib \qB}$ is a cartesian fibration.

For another cartesian fibration $q\colon \qF \tfib\qB$ with fibers $(\qF_b)_{b\in\ob\qB}$, we will define a natural equivalence of function complexes 
\[ \Fun^c_{\qB}(\qF \xtfib{q}\qB,  \prod\limits_{b \in \ob\qB}(\qE_b \times \qB \tfib \qB)^{ b \comma \qB\tfib \qB}) \we \prod\limits_{b\in\ob\qB} \Fun(\qF_b,\qE_b)\] and so establish the claimed adjoint correspondence.

To begin, the universal properties of the product and exponential provide isomorphisms
\begin{align*} \Fun_{\qB}(\qF \xtfib{q}\qB,  \prod_{b\in\ob\qB} (\qE_b \times \qB \tfib \qB)^{ b \comma \qB \tfib \qB})  &\cong \prod\limits_{b \in \ob\qB}\Fun_{\qB}(\qF \xtfib{q}\qB, (\qE_b \times \qB \tfib \qB)^{ b \comma \qB\tfib \qB}) \\ &\cong \prod\limits_{b\in\ob\qB} \Fun_{\qB}(b \comma \qB\tfib \qB, (\qE_b \times \qB\tfib \qB)^{\qF \xtfib{q} \qB}).\end{align*} By Lemma \ref{lem:cartesian-functor-transpose}, these isomorphisms restrict to the full sub quasi-categories spanned by the cartesian functors. Hence, 
\[ \Fun^c_{\qB}(\qF \xtfib{q}\qB,  \prod\limits_{b \in \ob\qB}(\qE_b \times \qB \tfib \qB)^{ b \comma \qB \tfib \qB})  \cong\prod\limits_{b \in \ob\qB}\Fun^c_{\qB}(b \comma \qB\tfib \qB, (\qE_b \times \qB\tfib \qB)^{\qF \xtfib{q} \qB}).\] By the dual of the Yoneda lemma, proven as Theorem \refIV{thm:yoneda}, restriction along the element $1 \to b \comma \qB$ corresponding to the identity at $b \colon 1 \to \qB$ defines an equivalence
\[ \Fun^c_{\qB}(b \comma \qB\tfib \qB, (\qE_b \times \qB\tfib \qB)^{\qF \xtfib{q} \qB})\we \Fun_{\qB}(1\xrightarrow{b} \qB,(\qE_b \times \qB\tfib \qB)^{\qF \xtfib{q} \qB}). \] The proof is finished by the isomorphisms
\[\Fun_{\qB}(1\xrightarrow{b} \qB,(\qE_b \times \qB\tfib \qB)^{\qF \xtfib{q} \qB}) \cong \Fun_{\qB}(\qF_b \to 1 \xrightarrow{b}\qB, \qE_b \times \qB \to \qB) \cong \Fun(\qF_b, \qE_b). \qedhere \]
\end{proof}

The construction of the left adjoint to $\Cart(\qCat)_{/\qB} \longrightarrow \Cart(\qCat)_{/\ob\qB}$ will make use of the tensor of a cartesian fibration, namely $\qB\comma b \tfib \qB$, with a quasi-category, namely $\qE_b$, described in Observation \ref{obs:tensor-pullback} and Lemma \ref{lem:cart-functor-tensor}.

\begin{prop}\label{prop:left-adjoint-construction}
The functor $\Cart(\qCat)_{/\qB} \longrightarrow \Cart(\qCat)_{/\ob\qB}$ admits a quasi-cat\-e\-gor\-ically  enriched left biadjoint defined by
\[ 
\xymatrix@R=1em{ \Cart(\qCat)_{/\ob\qB} \ar[r] & \Cart(\qCat)_{/\qB} \\(\qE_b)_{b \in\ob\qB} \ar@{}[r]|-\mapsto & \coprod_{b} \qE_b \times \qB\comma b\tfib \qB }
\]
that is, an $\ob\qB$-indexed family of quasi-categories $(\qE_b)_{b\in\ob\qB}$ is sent to the coproduct in $\Cart_{/\qB}$ of the cartesian fibrations $\qE_b \times \qB\comma b\xtfib{\pi} \qB\comma b\xtfib{p_0} \qB$.
\end{prop}

\begin{proof}
To make sense of this construction, note that the coproduct of cartesian fibrations over $\qB$ is again a cartesian fibration over $\qB$: since the horns are connected, each lifting problem of Definition \ref{defn:qcat-cocart} is supported in a single component. It follows that a functor out of the coproduct of cartesian fibrations is cartesian if and only if each of its legs is a cartesian functor.

For another cartesian fibration $q\colon \qF \tfib\qB$ with fibers $(\qF_b)_{b\in\ob\qB}$, we will define a natural equivalence of function complexes 
\[ \Fun^c_{\qB}( \coprod_{b\in\ob\qB} \qE_b \times \qB\comma b \tfib \qB, \qF \xtfib{q}\qB) \we \prod\limits_{b\in\ob\qB} \Fun(\qE_b,\qF_b)\] and so establish the claimed adjoint correspondence.

To begin, the universal property of the coproduct provides the first isomorphism, while the universal property of the tensor proven as Lemma \ref{lem:cart-functor-tensor} provides the second
\begin{align*}  \Fun^c_{\qB}( \coprod_{b\in\ob\qB} \qE_b \times \qB\comma b \tfib \qB, \qF \xtfib{q}\qB)  &\cong \prod\limits_{b \in \ob\qB}\Fun^c_{\qB}( \qE_b \times \qB\comma b \tfib \qB,\qF \xtfib{q}\qB) \\ &\cong \prod\limits_{b \in\ob\qB} \Fun^c_{\qB}(\qB\comma b \tfib \qB, \qF \xtfib{q}\qB)^{E_b}.\end{align*} 

 By the Yoneda lemma, proven as Theorem \refIV{thm:yoneda}, restriction along the element $1 \to \qB\comma b$ corresponding to the identity at $b \colon 1 \to \qB$ defines an equivalence
\[ \Fun^c_{\qB}(\qB\comma b\tfib \qB,  \qF \xtfib{q}\qB) \we \Fun_{\qB}(1\xrightarrow{b} \qB,  \qF \xtfib{q}\qB)  \cong \qF_b. \] This equivalence is respected by the cotensor $(-)^{\qE_b}$ and the product, so we have the desired equivalence
\[  \Fun^c_{\qB}( \coprod_{b\in\ob\qB} \qE_b \times \qB\comma b \tfib \qB, \qF \xtfib{q}\qB)  \cong \prod\limits_{b \in\ob\qB} \Fun^c_{\qB}(\qB\comma b \tfib \qB, \qF \xtfib{q}\qB)^{E_b} \simeq  \prod\limits_{b \in\ob\qB}\Fun(\qE_b,\qF_b).\qedhere\]
\end{proof}

\subsection{Monadicity and comonadicity}\label{sec:comonadicity}

A functor $u \colon \qA \to \qB$ between quasi-categories is \emph{comonadic} if it is the left adjoint part of a comonadic adjunction. This means that $\qA$ is equivalent to the \emph{quasi-category of coalgebras} for the \emph{homotopy coherent comonad} on $\qB$ underlying the corresponding \emph{homotopy coherent adjunction} derived from $u$ and its right adjoint. The quasi-category of coalgebras is defined as a particular flexible weighted limit of the homotopy coherent comonad. A few specific details of this construction are needed in the proofs in \S\ref{sec:groupoidal-hocoh-monad} and these are reviewed there. To avoid an unnecessary digression, we refer the reader \S\refII{sec:monadicity} for the definition of these notions and omit them from the current presentation. 

Recall Theorem \refII{thm:monadicity}, presented here in the dual:

\begin{thm}[{comonadicity \refII{thm:monadicity}}]\label{thm:comonadicity} A functor $u \colon \qA \to \qB$ between quasi-categories  is comonadic if and only if:
\begin{enumerate}[label=(\roman*)]
\item\label{itm:comonadicity-adj} $u$ admits a right adjoint,
\item\label{itm:comonadicity-split} $\qA$ admits and $u$ preserves limits of $u$-split cosimplicial objects, and
\item\label{itm:comonadicity-conserv} $u$ is conservative, reflecting isomorphisms.
\end{enumerate}
\end{thm}

Conservative functors between quasi-categories might arise as follows:

\begin{lem}\label{lem:conservative} Let $U \colon \eK \to \eL$ be a functor of quasi-categorically enriched categories that reflects equivalences in the sense that any 0-arrow $f \colon A \to B$ in $\eK$ whose image is an equivalence in $\eL$ is an equivalence in $\eK$. Then the corresponding functor $U \colon \qK \to \qL$ between quasi-categorical cores is conservative.
\end{lem}
\begin{proof} We show that an equivalence $f\colon A \to B$ in a quasi-categorically enriched category $\eK$ corresponds to an isomorphism in its quasi-categorical core. An equivalence in $\eK$ is comprised of the data enumerated in \ifnewadd Example \ref{ex:nerves-as-complicial-sets}\fi\ifnewcut\else\ref{defn:natural-marking}\fi\ref{itm:1-equivalence-data}, which is contained in the subcategory $g_*\eK\subset\eK$. In the homotopy coherent nerve $\qK=\hN g_*\eK$ this data gives rise to a pair of objects $A$ and $B$, a pair of 1-simplices $f \colon A \to B$ and $f' \colon B \to A$, and a pair of 2-simplices witnessing that $f$ and $f'$ compose to identities. This is the data that defines an isomorphism in a quasi-category.\end{proof}

Comonadic functors have the following property:

\begin{thm}[{comonadicity and colimit creation \refIII{thm:monadic-completeness}}]\label{thm:comonadic-cocompleteness} Let $u \colon \qA \to \qB$ be a comonadic functor between quasi-categories. Then $u$ creates any colimits that $\qB$ admits. \qed
\end{thm}

In particular any quasi-category admits colimits of split simplicial objects\footnote{Furthermore, such colimits are \emph{absolute}, that is preserved by any functor.}: a simplicial object $\Del\op \to \qB$ is \emph{split} if it extends along the inclusion $\Del\op\inc\Del+\op\inc\Del[t]$ that augments it with an terminal object --- note that $\Del+\op\cong\Del\op\join\catone$ --- and then adds an ``extra degeneracy'' map in each dimension (see \S\refI{sec:geo-realization}).

\begin{thm}[{split simplicial objects define colimits \refI{thm:splitgeorealizations}}]\label{thm:splitgeorealizations} Any split simplicial object $\Del\op \to \qB$ admits a colimit, whose colimit cone is given by the augmented diagram  $\Del+\op \to \qB$:
\[ \xymatrix{&  {~}\Del\op \ar[dr] \ar@{_(->}[dl]  \ar@{^(->}[d] & \\ \Del+\op \ar@{^(->}[r] &  \Del[t] \ar@{-->}[r] & \qB}
\]
\end{thm}

Combining these results, it follows that if a functor admits both left and right adjoints, its monadicity can be leveraged to help establish its comonadicity, or conversely:

\begin{prop}\label{prop:monadic-from-comonadic} A functor that admits both left and right adjoints is monadic if and only if it is comonadic.
\[
\xymatrix@C=3em{\qA\ar[r]|-u & \qB \ar@/_3ex/[l]_{\ell}^\perp \ar@/^3ex/[l]^r_\perp}\] 
\end{prop}
\begin{proof}
If $u$ is comonadic, then $u$ is conservative, verifying condition \ref{itm:comonadicity-conserv} of the dual Monadicity Theorem \ref{thm:comonadicity}. We have already assumed that the left adjoint required by \ref{itm:comonadicity-adj} exists. Finally, Theorem \ref{thm:splitgeorealizations} implies that $\qB$ admits colimits of $u$-split simplicial objects, and then comonadicity of $u$ together with Theorem \ref{thm:comonadic-cocompleteness} then implies that $\qA$ admits them as well and these are preserved by $u$. This verifies \ref{itm:comonadicity-split}, and Theorem \ref{thm:comonadicity} then implies that $u$ is also monadic. A dual argument proves the converse implication.
\end{proof}

\begin{thm}\label{thm:comonadic-cart} The forgetful functor 
\[ u\colon \Cart(\qqCat)_{/\qB} \longrightarrow \Cart(\qqCat)_{/\ob\qB} \cong \prod_{\ob\qB}\qqCat\]
is comonadic and hence also monadic.
\end{thm}
\begin{proof}
We use Theorem \ref{thm:comonadicity} to prove comonadicity and then deduce monadicity from Proposition \ref{prop:left-adjoint-construction} and Proposition \ref{prop:monadic-from-comonadic}. The right adjoint to $u$ is constructed in Proposition \ref{prop:right-adjoint-construction} proving \ref{itm:comonadicity-adj}. Example \refVIII{ex:comp-fibs} and Remark \refVIII{rmk:cosmos-pres-qcat-limits} combine to prove that $\Cart_{/\qB}$ admits and $\Cart_{/\qB} \inc \qqCat_{/\qB}$ preserves all limits. The functor $u$ is the composite of this inclusion with the projection functor $\qqCat_{/\qB} \to \qqCat_{/\ob\qB}$, which also preserves all limits, being the homotopy coherent nerve of a functor of $\infty$-cosmoi $\qCat_{/\qB} \to \qCat_{/\ob\qB}$.  This proves \ref{itm:comonadicity-split}. Proposition \ref{prop:equivalence-of-fibrations} and Lemma \ref{lem:conservative} assert that $u \colon \Cart_{/\qB} \to \Cart_{/\ob\qB}$ is conservative, proving \ref{itm:comonadicity-conserv}. 
\end{proof}

%% file: reflection.tex
\section{Groupoidal reflection}\label{sec:groupoidal-reflection}

In this section, we give a first application of the monadicity and comonadicity results of the previous section. A cartesian fibration between quasi-categories is \emph{groupoidal} if its fibers are Kan complexes, rather than quasi-categories. In this section we construct a reflection to the fully faithful inclusion of the quasi-category of groupoidal cartesian fibrations into the quasi-category of cartesian fibrations:
\[ \adjdisplay \mathrm{invert} -| {} : \Cart\gr_{/\qB} -> \Cart_{/\qB}.\] A different proof of this result will be appear in a sequel, making use of explicit fiberwise coinverters.

In \S\ref{sec:qcat-into-kan}, we study the relationship between the quasi-categorically enriched category of cartesian fibrations and its subcategory of groupoidal cartesian fibrations and establish a groupoidal reflection functor in the ``base case,'' reflecting quasi-categories into Kan complexes. In \S\ref{sec:groupoidal-monadicity}, we prove that the monadic and comonadic adjunctions of Theorem \ref{thm:comonadic-cart} restrict to define analogous monadic and comonadic adjunctions for groupoidal cartesian fibrations. It follows that the large quasi-categories $ \Cart\gr_{/\qB}$ and $\Cart_{/\qB}$ can be understood as quasi-categories of algebras for closely related homotopy coherent monads acting on $\prod_{b \in\ob\qB}\qKan$ and $\prod_{b \in\ob\qB}\qqCat$ respectively. In \S\ref{sec:groupoidal-hocoh-monad}, we exploit this presentation to construct an adjunction defining the reflection of a cartesian fibration into a groupoidal cartesian fibration:
\[ \adjdisplay \mathrm{invert} -| {~} : \Cart\gr_{/\qB} -> \Cart_{/\qB}.\] 

\subsection{Reflecting quasi-categories into Kan complexes}\label{sec:qcat-into-kan}

Before we begin, we note that the notion of groupoidal cartesian fibration of quasi-categories just defined agrees with the definition given in \S\refIV{ssec:groupoidal}, which declares that a cartesian fibration $p \colon \qE \tfib \qB$ is \emph{groupoidal} just when the quasi-category of functors from any $f \colon X \to \qB$ to $p$ is a Kan complex\ifnewadd; see X.12.2.3.\fi\ifnewcut\else. The reader who is content to work with the present definition may safely skip this bit of bookkeeping.

\begin{lem}\label{lem:groupoidal-quasi-is-fiberwise} Let $p \colon \qE \tfib \qB$ be an isofibration between quasi-categories. Then $p$ is groupoidal as an object of $\qCat_{/\qB}$ if and only if the fibers of $p$ are Kan complexes.\footnote{Warning: this result does not hold for generic $\infty$-cosmoi, having to do with the fact that the terminal quasi-category is a generating object in $\qCat$ in a suitable sense.}
\end{lem}
\begin{proof}
The quasi-category of functors in $\qCat_{/\qB}$ from $b \colon 1 \to \qB$ to $p \colon \qE \to \qB$ is the fiber $\qE_b$ of $p$ over $b$, so if $p$ is groupoidal as an object of $\qCat_{/\qB}$, then $p$ necessarily has Kan complex fibers. For the converse, we must argue that the pullback
\[ 
\xymatrix{ \Fun_{\qB}(f \colon X \to \qB, p \colon \qE \tfib \qB) \pbexcursion \ar@{->>}[d] \ar[r] & \Fun(X,\qE) \ar@{->>}[d]^{\Fun(X,p)} \\ 1 \ar[r]^-f & \Fun(X,\qB)}
\] is a Kan complex supposing that $p$ has groupoidal fibers. Unraveling the definition, we must show that any map $h$ that lies over $f$ in the sense of the following commutative diagram
\[ \xymatrix{ X \times \Del^1 \ar[d]_{\pi} \ar[r]^-h & \qE \ar@{->>}[d]^p \\ X \ar[r]^-f & \qB}\] defines an isomorphism $\qE^X$. Corollary \refI{cor:pointwise-equiv} observes that such 1-simplices are invertible if and only if each component $h(x,-) \colon \Del^1 \to \qE$ indexed by a vertex $x \in X$ defines an isomorphism in $\qE$. As this 1-simplex lives in the fiber over $px$, it is invertible, so the result follows.\end{proof}\fi

At the level of simplicially-enriched categories, the subcategory of groupoidal cartesian fibrations is defined by the pullback:
\begin{equation}\label{eq:groupoidal-fibrations-as-pullback}
\xymatrix{ \Cart\gr(\qCat)_{/\qB} \ar@{^(->}[r] \ar[d] \pbexcursion & \Cart(\qCat)_{/\qB} \ar[d]  \\ \Cart\gr(\qCat)_{/\ob\qB} \cong \prod\limits_{\ob\qB} \Kan \ar@{^(->}[r] & \prod\limits_{\ob\qB}\qCat \cong \Cart(\qCat)_{/\ob\qB}}
\end{equation}
In  \S\ref{sec:groupoidal-hocoh-monad}, we construct a \emph{groupoidal reflection functor}, by which we mean a left adjoint to the inclusion \[\Cart\gr(\qqCat)_{/\qB}\inc\Cart(\qqCat)_{/\qB}\] as a functor between large quasi-categories. We begin by describing groupoidal reflection in the case where $\qB=1$.

\begin{thm}\label{thm:basic-groupoidal-reflection} The inclusion $\qKan\inc\qqCat$ admits both left and right adjoints 
\[ 
\xymatrix@C=3em{\qKan \ar@{^(->}[r]^-\perp_\perp & \qqCat \ar@/_3ex/[l]_{\mathrm{invert}} \ar@/^3ex/[l]^\core}\]
and is monadic and comonadic.
\end{thm}
\begin{proof}
The inclusion $\Kan\inc\qCat$ is left adjoint to a functor $\core \colon \qCat\to\Kan$ that carries a quasi-category to the maximal sub Kan complex spanned by the isomorphisms\ifnewcut\else; this is a restriction of the 0-sharp $0$-core adjunction of \ref{defn:marked-basics}\ref{itm:core-adjunction} to the subcategories of saturated 1-complicial and saturated 0-complicial sets\fi. This adjunction is simplicial with respect to the Kan complex enrichments of both $\Kan$ and $\qCat$, the latter obtained by applying the $\core$ functor to the function complexes, so this simplicially enriched adjunction descends to provide a right adjoint to $\Kan\inc\qCat$.

The left adjoint can also be modeled at the point-set level. The quasi-category $\qqCat$ is isomorphic to the homotopy coherent nerve of the Kan-complex enriched category of naturally marked quasi-categories in the sense of Example \ref{ex:qcat-as-comp}. There is a simplicial Quillen adjunction
\[ \adjdisplay U -| (-)^\sharp : \sSet -> \msSet.\] connecting this simplicial model structure for quasi-categories to the Quillen model structure for Kan complexes on simplicial sets. Applying Theorem \refI{thm:simplicial-Quillen-adjunction}, this provides the left adjoint to the inclusion $\qKan\inc\qqCat$. From this vantage point, we may apply Proposition \refVII{prop:gpdal-infty-cosmos} to see that $\Kan$ is closed in $\qCat$ under flexible weighted limits, so we conclude that $\qKan\inc\qqCat$ creates all limits.

The functor $\Kan\inc\qCat$ reflects equivalences, so by Lemma \ref{lem:conservative} the inclusion $\qKan\inc\qqCat$ is conservative. Now Theorem \ref{thm:comonadicity} implies that $\qKan\to\qqCat$ is comonadic, and Proposition \ref{prop:monadic-from-comonadic} then implies that  $\qKan\to\qqCat$ is also monadic.
\end{proof}

\subsection{(Co)monadicity of groupoidal cartesian fibrations}\label{sec:groupoidal-monadicity}

We now argue that 
\[ \xymatrix@R=0.5em{\Cart\gr(\qCat)_{/\qB} \ar[r] & \Cart\gr(\qCat)_{/\ob\qB}\cong \prod\limits_{\ob\qB} \eop{Kan} \\ (\qE \xtfib{p} \qB) \ar@{}[r]|-\mapsto & (\qE_b)_{b \in \ob\qB}} \]
admits left and right quasi-categorically enriched biadjoints, given by restricting those from the non-groupoidal case, and that moreover the restricted adjunction is both monadic and comonadic at the level of functors between underlying quasi-categories.

\begin{thm}\label{thm:groupoidal-cartesian-monadicity} The quasi-categorical biadjoints to $\Cart(\qCat)_{/\qB} \to \Cart(\qCat)_{/\ob\qB}$ restrict to groupoidal cartesian fibrations
\[
\xymatrix{ \Cart\gr(\qCat)_{/\qB} \ar@{^(->}[r] \ar[d] & \Cart(\qCat)_{/\qB} \ar[d]  \\ \Cart\gr(\qCat)_{/\ob\qB}  \ar@{^(->}[r] \ar@/_3ex/@{-->}[u]_R^\dashv \ar@/^3ex/@{-->}[u]^L_\dashv& \Cart(\qCat)_{/\ob\qB}\ar@/_3ex/[u]_R^\dashv \ar@/^3ex/[u]^L_\dashv}
\]
 and moreover these restricted adjunctions display the functor between the quasi-categorical cores $\Cart\gr_{/\qB} \to \Cart\gr_{/\ob\qB}\cong \prod_{/\ob\qB}\qKan$ as both monadic and comonadic.
\end{thm}
\begin{proof}
Proposition \ref{prop:left-adjoint-construction} defines $L \colon \Cart(\qCat)_{/\ob\qB} \to \Cart(\qCat)_{/\qB}$ to be the functor that carries a family $(\qE_b)_{b \in \ob\qB}$ to 
\[ L( (\qE_b)_{b \in \ob\qB}) \defeq \coprod_{b} \qE_b \times \qB\comma b \tfib \qB.\] The fiber over $x \in \ob\qB$ is
$\coprod_{b} \qE_b \times x\comma b$. Since $x\comma b$ is a Kan complex, it is clear that this fiber is groupoidal if each $\qE_b$ is a Kan complex. Thus, we see immediately that $L$ restricts to groupoidal cartesian fibrations.

Proposition \ref{prop:right-adjoint-construction} defines $R \colon \Cart(\qCat)_{/\ob\qB} \to \Cart(\qCat)_{/\qB}$ to be the functor that carries a family $(\qE_b)_{b \in \ob\qB}$ to
\[ R((\qE_b)_{b \in \ob\qB}) \defeq  \prod_{b} (\qE_b \times \qB \tfib \qB)^{ b \comma \qB \tfib \qB}.\] The fiber over $x \in \ob\qB$ of this product of fibrations is isomorphic to the product of the fibers of each $(\qE_b\times\qB \tfib \qB)^{b\comma \qB\tfib \qB}$ over $x$, so it suffices to show that each of these fibers is a Kan complex if $\qE_b$ is a Kan complex. By the bijection of Definition \ref{defn:exponential}, a 1-simplex in the fiber of $(\qE_b\times\qB \tfib \qB)^{b\comma \qB\tfib \qB}$ over $x \colon 1 \to \qB$ corresponds to the displayed dashed map
\[
\xymatrix{ \Del^1 \times b \comma x \ar@{->>}[d] \ar[r] \pbexcursion \ar@/^3ex/@{-->}[rrr] &b \comma x \ar@{->>}[d] \ar[r] \pbexcursion & b\comma \qB \ar@{->>}[d] & \qE_b \times \qB \ar@{->>}[dl]^{\pi} \\ \Del^1 \ar[r] & \Del^0 \ar[r]_x & \qB}\]
i.e., to a map $\Del^1 \times b\comma x \to \qE_b$. If each $\qE_b$ is a Kan complex, this map can be extended along the inclusion $\Del^1 \inc\iso$ from the 1-simplex into the free-living isomorphism. This proves that every 1-simplex in the fiber of $(\qE_b\times\qB \tfib \qB)^{b\comma \qB\tfib \qB}$ is an isomorphism, which tells us that $R$ restricts to groupoidal cartesian fibrations.

By what is now a familiar line of argument, we apply Theorem \ref{thm:comonadicity} to prove that the restricted adjunctions are comonadic, and then deduce monadicity from Proposition \ref{prop:monadic-from-comonadic}. The required adjoints have already been constructed and Proposition \ref{prop:equivalence-of-fibrations} and Lemma \ref{lem:conservative} imply that $\Cart\gr_{/\qB} \to \Cart\gr_{/\ob\qB}$ is conservative, so it remains only to establish condition \ref{itm:comonadicity-split} of Theorem \ref{thm:comonadicity}.

The pullback of simplicial categories \eqref{eq:groupoidal-fibrations-as-pullback} is preserved by passing to the level of quasi-categories:
\[
\xymatrix{ \Cart\gr_{/\qB} \ar@{^(->}[r] \ar[d] \pbexcursion & \Cart_{/\qB} \ar[d]  \\  \prod\limits_{\ob\qB} \qKan \ar@{^(->}[r] & \prod\limits_{\ob\qB}\qqCat}
\]
By the monadicity established in Theorems \ref{thm:comonadic-cart} and \ref{thm:basic-groupoidal-reflection} and the dual of Theorem \ref{thm:comonadic-cocompleteness}, both the left-hand vertical and lower horizontal functors create all limits present in $\prod_{\ob\qB}\qqCat$, which is to say all limits, since by Proposition \refVII{prop:qcat-of-cosmos-complete} $\qqCat$ is complete. The lower inclusion is replete up to isomorphism, so $\qKan\inc\qqCat$ defines an isofibration of large quasi-categories. Thus, Lemma \ifnewcut\else\ref{lem:pullback-colimits} below\fi\ifnewadd X.6.3.17 \fi applies to tell us that $\Cart\gr(\qqCat)_{/\qB}$ is also complete and all limits are preserved by the left-hand vertical functor. 

Now Theorem \ref{thm:comonadicity} implies that $\Cart\gr_{/\qB} \to \Cart\gr_{/\ob\qB} \cong \prod_{\ob\qB}\qKan$  is comonadic, and monadicity follows from Proposition \ref{prop:monadic-from-comonadic}.
\end{proof}

\ifnewcut\else
\begin{lem}\label{lem:pullback-colimits} Consider a pullback diagram of quasi-categories  whose vertical morphisms are isofibrations
\[ \xymatrix{ \qE \pbexcursion \ar@{->>}[d]_p \ar[r]^f & \qF \ar@{->>}[d]^q \\ \qB \ar[r]_g & \qA}\] Then $p$ creates and $f$ preserves  any class of limits or colimits that $g$ preserves and $q$ creates.
\end{lem}
\begin{proof}
This can be proven directly or by appealing to Theorem \refIII{thm:limits-in-limits}, which states that the subcategory of $\infty$-categories admitting and functors preserving limits or colimits of a particular variety is closed under flexible weighted limits. Since $q \colon \qF \tfib \qA$ is an isofibration, the strict limit of the cospan $\qB \xrightarrow{g} \qA \xfibt{q} \qF$ is equivalent to the limit weighted by the projective cofibrant weight $\pbshape \to \sSet$ with image $\Del^0 \xrightarrow{\face^0} \iso \xleftarrow{\face^1} \Del^0$. Theorem  \refIII{thm:limits-in-limits} now implies that $\qE$ admits and $f$ and $p$ preserve any limits or colimits present in $\qB$, $\qA$, and $\qF$ and preserved by $g$ and $q$. 

It remains only to argue that (co)limits are created by $p \colon \qE \tfib \qB$. Consider a family of diagrams $d \colon \qD \to \qE^X$ so that $pd \colon \qD \to \qB^X$ has a colimit $b \colon \qD \to \qB$. Then $gb \colon \qD \to \qA$ is a colimit for $gpd = qfd\colon \qD \to \qA^X$, and the fact that $q$ creates such colimits implies that there is an colimit object $c \colon \qD \to \qB$ for $fd$ with $qc$ isomorphic to $gb$ in $\Fun(\qD,\qA)$. In fact, since 
$\Fun(\qD,q) \colon \Fun(\qD,F) \tfib \Fun(\qD,\qA)$ is an isofibration, we may assume that $qc = gb$, so this pair induces an object $e \colon \qD \to \qE$ with $pe = b$. By construction, there is a pullback of contravariant represented modules displayed below left
\[ \xymatrix{ e\comma \qE \ar@{->>}[d]_p \ar[r]^f \pbexcursion & c \comma \qB \ar@{->>}[d]^q & & d \comma \Delta \pbexcursion \ar[r]^f \ar@{->>}[d]_p & fd \comma \Delta \ar@{->>}[d]^q \\ b \comma \qB \ar[r]_g & gb \comma \qA & & pd \comma \Delta \ar[r]_g & gpd \comma \Delta}\]
whose underlying cospan is equivalent to the cospan in the pullback of quasi-categories of cones displayed above right. It follows that $e \comma \qE \simeq d \comma \Delta$, which says that $e$ is a colimit for $d$.
\end{proof}
\fi
\subsection{Groupoidal reflection}\label{sec:groupoidal-hocoh-monad}

Our monadicity results, Theorems \ref{thm:groupoidal-cartesian-monadicity} and \ref{thm:comonadic-cart}, tell us that the quasi-categories $\Cart_{/\qB}$ and $\Cart\gr_{/\qB}$ are equivalent to the quasi-categories of algebras associated to closely related homotopy coherent monads acting on $\Cart_{/\ob\qB} \cong \prod_{\ob\qB}\qqCat$ and $\Cart\gr_{/\ob\qB}\cong\prod_{\ob\qB}\qKan$ respectively. In this section, we will use this result to lift the reflection functor $\mathrm{invert} \colon \qqCat \to \qKan$ from quasi-categories to Kan complexes to a groupoidal reflection functor $\mathrm{invert} \colon \Cart_{/\qB} \to \Cart\gr_{/\qB}$ that is left adjoint to the inclusion.

To do this we make use of a convenient representation for adjoint functors that can be expressed in any 2-category, dual to the more familiar representation of the unit of an adjunction as an absolute left extension diagram:

\begin{lem}[{\refI{ex:adjasabslifting}}]\label{lem:adj-as-abs-lifting} To define a left adjoint to a functor $u \colon A \to B$ is to define an absolute left lifting of $\id_B$ along $u$, in which case $f \dashv u$ with unit $\eta \colon \id_B \To uf$.
\[
\xymatrix{ \ar@{}[dr]|(.7){\Uparrow\eta} & A \ar[d]^u \\ B \ar@{=}[r]  \ar@{-->}[ur]^f & B} \qed
\] 
\end{lem}

Let $\pbshape$ denote the category indexing a cospan and write $\qCat^\pbshape$ for the simplicially enriched category of cospans of quasi-categories, whose objects are cospans and whose 0-arrows are natural transformations
  \begin{equation}\label{eq:obj-K-pbshape-trans}
    \xymatrix@=1.5em{
      {} & {\qB}\ar[d]_{f}\ar[dr]^{v} & {} \\
      {\qC}\ar[r]^{g}\ar[dr]_{w} & {\qA}\ar[dr]|*+<1pt>{\scriptstyle u} & 
      {\qB'}\ar[d]^{f'} \\
      {} & {\qC'}\ar[r]_{g'} & {\qA'}
    }
  \end{equation}
  
  \begin{defn}\label{defn:abs-lift-trans-left-exactness}
Transformations of the kind depicted in~\eqref{eq:obj-K-pbshape-trans} between diagrams which admit absolute left liftings give rise to the following diagram
  \begin{equation}\label{eq:induced-mate}
    \vcenter{
      \xymatrix@=1.5em{
        {} \ar@{}[dr]|(.7){\Uparrow\lambda}
        & {\qB}\ar[d]^{f}\ar[dr]^{v} & {} \\
        {\qC}\ar[ur]^{\ell}\ar[r]_{g}\ar[dr]_{w} & 
        {\qA}\ar[dr]|*+<1pt>{\scriptstyle u} & 
        {\qB'}\ar[d]^{f'} \\
        {} & {\qC'}\ar[r]_{g'} & {\qA'}
      }
    } 
    \mkern40mu = \mkern40mu
    \vcenter{
      \xymatrix@=1.5em{
        {} & {\qB}\ar[dr]^{v} & {} \\
        {\qC}\ar[ur]^{\ell}\ar[dr]_{w} & 
        {\scriptstyle\Uparrow\tau}\ar@{}[dr]|(.7){\Uparrow\lambda'} & 
        {\qB'}\ar[d]^{f'} \\
        {} & {\qC'}\ar[ur]^{\ell'}\ar[r]_{g'} & {\qA'}
      }
    } 
  \end{equation}
  in which the triangles are absolute left liftings and the 2-cell $\tau$ is induced by the universal property of the triangle on the right. We say that the transformation \eqref{eq:induced-mate} is {\em left exact\/} if and only if the induced 2-cell $\tau$ is an isomorphism. This left exactness condition holds if and only if, in the diagram on the left, the whiskered 2-cell $u\lambda$ displays $v\ell$ as the absolute left lifting of $g'w$ through $f'$.
\end{defn}

Our interest in these notions is on account of the following result:

\begin{prop}[{\refIII{prop:abslifts-in-limits}}]\label{prop:abslifts-in-limits} Consider any simplicial functor $T \colon \eA \to \qCat^\pbshape$ and any flexible weight $W \colon \eA \to \SSet$. If each of the objects in the image of $T$ admits an absolute left lifting and each of the 0-arrows in the image of $T$ is left exact, then the weighted limit $\lim^W\!T \in \qCat^\pbshape$ admits an absolute left lifting and the legs of the limit cone are left exact transformations. \qed
\end{prop}

 By Lemma \refII{lem:projcof1}, the quasi-category of algebras construction, introduced in Definition \refII{defn:EMobject}, is an instance of a flexible weighted limit. We will use Proposition \ref{prop:abslifts-in-limits} applied in a larger Grothendieck universe to the quasi-categorically enriched category $\eop{QCAT}$ of large quasi-categories to lift the absolute left lifting diagram
\[ 
\xymatrix{ \ar@{}[dr]|(.7){\Uparrow\eta} & \qKan \ar@{^(->}[d] \\ \qqCat \ar[ur]^{\mathrm{invert}} \ar@{=}[r] & \qqCat}\]
whose left adjoint part defines the groupoidal reflection associated to the inclusion
\[ \Cart\gr_{/\ob\qB}\cong\prod\limits_{\ob\qB}\qKan \inc \prod\limits_{\ob\qB}\qqCat\cong \Cart_{/\ob\qB}\]
to these flexible weighted limits, defining an absolute left lifting diagram
\[ 
\xymatrix{ \ar@{}[dr]|(.7){\Uparrow\eta} & \Cart\gr_{/\qB} \ar@{^(->}[d] \\ \Cart_{/\qB} \ar@{-->}[ur]^{\mathrm{invert}} \ar@{=}[r] &  \Cart_{/\qB}}\]
the left adjoint being the desired groupoidal reflection functor.

\begin{thm}\label{thm:groupoidal-reflection-cartesian}
There is a left adjoint to the inclusion 
\[
\begin{tikzcd}[sep=large] \Cart\gr_{/\qB} \arrow[r, bend right] \arrow[r, phantom, "\bot"] & \Cart_{/\qB} \arrow[l, bend right, "\mathrm{invert}"']
\end{tikzcd}
\]defining the reflection of a cartesian fibration into a groupoidal cartesian fibration.
\end{thm}
\begin{proof}
By Theorem \refII{thm:hty-coherence-exist-I}, the adjunction 
\begin{equation}\label{eq:this-adjunction}
\begin{tikzcd}[sep=large]\Cart_{/\qB} \arrow[r, bend right] \arrow[r, phantom, "\bot"] & \prod_{\ob\qB}\qqCat \arrow[l, bend right, "L"']
\end{tikzcd}
\end{equation}
defined in Proposition \ref{prop:left-adjoint-construction} extends to a \emph{homotopy coherent adjunction}: a simplicial functor $\Adj \to \eop{QCAT}$ valued in the quasi-categorically enriched category of large quasi-categories whose domain is the simplicial computad obtained by applying the nerve functor to the hom-categories of the free 2-category containing an adjunction; see \S\refII{sec:generic-adj}.  The two objects of $\Adj$, called ``$-$'' and ``$+$'' are mapped to the quasi-categories $\Cart_{/\qB}$ and $\prod_{\ob\qB}\qqCat$ respectively. The full subcategory $\Mnd\inc\Adj$ on the object ``$+$'' defines the free \emph{homotopy coherent monad}. In this case, the data of the underlying homotopy coherent monad $T \colon \Mnd\inc\Adj \to \eop{QCAT}$ on the object $\prod_{\ob\qB}\qqCat$ is given by a the map of objects $+ \mapsto \prod_{\ob\qB}\qqCat$ and the map of function complexes
\[T\colon \Del+ = \Fun_{\Mnd}(+,+) \to \Fun_{\eop{QCAT}}(\prod_{\ob\qB}\qqCat,\prod_{\ob\qB}\qqCat),\] satisfying an appropriate simplicial functoriality condition. This functoriality condition implies that the 0-arrows of $\Mnd$, are indexed by the objects $[n] \in \Del+$, are all finite composites of the object $[0]$, whose image
\[ t \defeq T[0] \colon \prod_{\ob\qB}\qqCat \to \prod_{\ob\qB}\qqCat\] is the monad endofunctor defined by composing the left and right adjoints of \eqref{eq:this-adjunction}.

To apply Proposition \ref{prop:abslifts-in-limits} we must extend the homotopy coherent monad $T$ to a homotopy coherent monad on $\eop{QCAT}^\pbshape$. To do so, we argue that this homotopy coherent monad restricts along $\prod_{\ob\qB}\qKan\inc\prod_{\ob\qB}\qqCat$ to define a homotopy coherent monad $T\gr \colon \Mnd \to \eop{QCAT}$ on $\prod_{\ob\qB}\qKan$ in such a way that this map will define the component of a simplicial natural transformation $T\gr\Rightarrow T$. To see this, note that $\prod_{\ob\qB}\qKan$, the nerve of the simplicially enriched category spanned by $\ob\qB$-indexed families of Kan complexes, is a full sub quasi-category of $\prod_{\ob\qB}\qqCat$, the nerve of the simplicially enriched category spanned by $\ob\qB$-indexed families of quasi-categories, in the sense that it contains all of the $n$-simplices whose vertices are Kan complexes, not mere quasi-categories. So to check that the data of the homotopy coherent monad restricts to define a simplicial functor given by $+ \mapsto \prod_{\ob\qB}\qKan$ and
\[ T\gr\colon\Del+ = \Fun_{\Mnd}(+,+) \to \Fun_{\eop{QCAT}}(\prod_{\ob\qB}\qKan,\prod_{\ob\qB}\qKan),\] it suffices to check this at the level of vertices $[n] \in\Del_+$, which amounts to checking that the monad $t$ of \eqref{eq:this-adjunction} restricts to define a monad $t\gr$ on groupoidal cartesian fibrations; this was done in Theorem \ref{thm:groupoidal-cartesian-monadicity}. 

In this way we obtain a simplicial functor $\Mnd \to \eop{QCAT}^\cattwo$ sending the object ``$+$'' to the arrow $\prod_{\ob\qB}\qKan \inc\prod_{\ob\qB}\qqCat$. Pairing this with the identity simplicial natural transformation we obtain a simplicial functor $\Mnd \to \eop{QCAT}^\pbshape$ sending the object ``$+$'' to the cospan displayed below-left
\begin{equation}\label{eq:hocoh-monad-cospan} \xymatrix{ & \prod_{\ob\qB}\qKan \ar@{^(->}[d]  & \ar@{}[dr]|(.7){\Uparrow\eta} &  \prod\limits_{\ob\qB}\qKan \ar@{^(->}[d] \\ \prod_{\ob\qB}\qqCat \ar@{=}[r] & \prod_{\ob\qB}\qqCat &  \prod\limits_{\ob\qB} \qqCat \ar[ur]^{\prod_{\ob\qB}\mathrm{invert}} \ar@{=}[r] &   \prod\limits_{\ob\qB} \qqCat }\end{equation} 
This cospan admits an absolute left lifting displayed above right, defining the left adjoint and unit of an adjunction whose counit is invertible. In fact the entire diagram $\Mnd \to \eop{QCAT}^\pbshape$ restricts to the subcategory  spanned by those cospans that admit absolute left liftings and those 0-arrows that define left exact transformations between them. To see this, we need only argue that the generating 0-arrow $[0]$ in $\Mnd$, the endofunctor of the monad, defines a left exact transformation. That is, we must show that the endotransformation of \eqref{eq:hocoh-monad-cospan} whose components are the functors
\[ 
\xymatrix@R=0.5em{ \prod_{\ob\qB}\qKan \ar[r]^{t\gr} & \prod_{\ob\qB}\qKan & \prod_{\ob\qB}\qqCat \ar[r]^t & \prod_{\ob\qB}\qqCat \\ (\qE_b)_{b \in \qB} \ar@{}[r]|-\mapsto & \left(\coprod\limits_{b\in\ob\qB} \qE_b \times x\comma b \right)_{x \in\ob\qB} & (\qE_b)_{b \in \qB} \ar@{}[r]|-\mapsto & \left(\coprod\limits_{b\in\ob\qB} \qE_b \times x\comma b \right)_{x \in\ob\qB} 
}
\] is left exact. This amounts to showing that the whiskered 2-cell
\[
\xymatrix{  \ar@{}[dr]|(.7){\Uparrow\eta} &  \prod\limits_{\ob\qB}\qKan\ar@{^(->}[d] \ar[r]^{t\gr} &  \prod\limits_{\ob\qB}\qKan\ar@{^(->}[d]  \\   \prod\limits_{\ob\qB} \qqCat \ar[ur]^{\prod_{\ob\qB}\mathrm{invert}} \ar@{=}[r] &   \prod\limits_{\ob\qB} \qqCat \ar[r]_{t} &  \prod\limits_{\ob\qB} \qqCat }
\]
is invertible. This is the case because the process of freely inverting a family of quasi-categories commutes up to equivalence with forming the product with the Kan complex $x\comma b$ and with the coproduct $\coprod_{b\in\ob\qB}$.

In this way we obtain a homotopy coherent monad $\Mnd \to\eop{QCAT}^\pbshape$ valued in the subcategory of cospans admitting absolute left liftings and left exact transformations between them. There is a flexible  weight $W_- \colon \Mnd \to \SSet$ introduced in Definition \refII{defn:EMobject} --- the precise details of which are not relevant here --- so that the $W_-$-weighted limit of a homotopy coherent monad define its quasi-category of algebras, as characterized up to equivalence by the Monadicity Theorem \ref{thm:comonadicity}. By Theorems \ref{thm:groupoidal-cartesian-monadicity} and  \ref{thm:comonadic-cart} the $W_-$-weighted limit of the composite diagram $\Mnd\to\eop{QCAT}^\pbshape$ defines the cospan displayed below and by Proposition \ref{prop:abslifts-in-limits} it therefore admits an absolute left lifting:
\[ 
\xymatrix{ \ar@{}[dr]|(.7){\Uparrow\eta} & \Cart\gr_{/\qB} \ar@{^(->}[d] \\ \Cart_{/\qB} \ar@{-->}[ur]^{\mathrm{invert}} \ar@{=}[r] &  \Cart_{/\qB}}\]
By Lemma \ref{lem:adj-as-abs-lifting}, this absolute left lifting diagram defines the adjunction that constructs the groupoidal reflection of a cartesian fibration.
\end{proof}
